\documentclass[11pt,twoside]{article}

\usepackage{fancyhdr}
\usepackage[colorlinks,citecolor=blue,urlcolor=blue,linkcolor=blue,bookmarks=false]{hyperref}
\usepackage{epsf,epsfig,graphics}
\graphicspath{ {./figures/} }
\usepackage{subcaption}
\usepackage{float}
\usepackage{afterpage}
\usepackage{fullpage}
\usepackage[T1]{fontenc}
\usepackage{psfrag,xspace}
\usepackage{color,etoolbox}
\usepackage{bbm}
\usepackage{scalerel}
\usepackage{multirow}
\usepackage{import}
\usepackage{todonotes}
\usepackage{multicol}
\usepackage{doi}
\usepackage[font=small,skip=0pt,labelfont=bf]{caption}
\usepackage[export]{adjustbox}

\usepackage{tcolorbox}

\newtcolorbox{customcolorbox}[3][]
{ colframe = #2!30, colback  = #2!8, coltitle = #2!30!black, title    = {#3},
  #1, }

\usepackage{csquotes}
\newenvironment{itquote}
{\begin{quote}\itshape} {\end{quote}}

\usepackage[square,comma,authoryear,sort&compress]{natbib}

\usepackage{./Latex/macros_01_proof_envs}
\usepackage{./Latex/macros_02_operators_commands}
\usepackage{./Latex/macros_03_symbols}
\usepackage{./Latex/macros_04_propernames}
\usepackage{./Latex/macros_05_reviewer}

\newcommand{\Prbb}[2]{\mathop{\bf Pr\/}\left(#1 \,\middle|\, #2 \right)}
\newcommand{\resp}{R}
\newcommand{\absresp}{T}
\newcommand{\univupperboundmu}{C}

\newcommand{\univupperboundsigma}{\Psi}
\newcommand{\univlowerboundsigma}{\psi}

\newcommand{\muhatascifit}{\boldsymbol{\widehat{\mu}_{\textnormal{ascifit}}} }
\newcommand{\muhatnaive}{\boldsymbol{\widehat{\mu}_{\textnormal{naive}}} }

\newtoggle{neurips2022sub} 
\newtoggle{arxivsub} 
\newtoggle{anonymoussub} 

\usepackage{tabulary}
\usepackage{booktabs}

\usepackage{inconsolata}

\usepackage{enumitem}
\usepackage{float}
\usepackage{color,etoolbox}

\newtcolorbox{mycolorbox}[3][]
{ colframe = #2!25, colback  = #2!10, coltitle = #2!20!black, title    = {#3},
  #1, }

\usepackage{algorithm}
\usepackage{algpseudocode}

\usepackage[toc,page,header]{appendix}
\usepackage{minitoc}

\toggletrue{arxivsub} 

\begin{document}
\doparttoc 
\faketableofcontents 

\begin{center} {\LARGE{\bf{Adversarial Sign-Corrupted Isotonic Regression}}} \\

    \vspace*{.3in}

    {\large {
            \begin{tabular}{ccccc}
                Shamindra Shrotriya & Matey Neykov
            \end{tabular}

            \vspace*{.1in}

            \begin{tabular}{ccc}
                Department of Statistics \& Data Science \\
            \end{tabular}

            \begin{tabular}{c}
                Carnegie Mellon University \\
                Pittsburgh, PA 15213
            \end{tabular}

            \vspace*{.2in}

            \begin{tabular}{c}
                {\texttt{\{sshrotri, mneykov\}@andrew.cmu.edu}}
            \end{tabular}
        }

    }

    \vspace*{.2in}

    \today

    \vspace*{.2in}

    \begin{abstract}
    \vspace*{.2in}
    Classical univariate isotonic regression involves nonparametric estimation
    under a monotonicity constraint of the true signal. We consider a variation
    of this generating process, which we term adversarial sign-corrupted
    isotonic (\texttt{ASCI}) regression. Under this \texttt{ASCI} setting, the
    adversary has full access to the true isotonic responses, and is free to
    sign-corrupt them. Estimating the true monotonic signal given these
    sign-corrupted responses is a highly challenging task. Notably, the
    sign-corruptions are designed to violate monotonicity, and possibly induce
    heavy dependence between the corrupted response terms. In this sense,
    \texttt{ASCI} regression may be viewed as an adversarial stress test for
    isotonic regression. Our motivation is driven by understanding whether
    efficient robust estimation of the monotone signal is feasible under this
    adversarial setting. We develop \texttt{ASCIFIT}, a three-step estimation
    procedure under the \texttt{ASCI} setting. The \texttt{ASCIFIT} procedure is
    conceptually simple, easy to implement with existing software, and consists
    of applying the \texttt{PAVA} with crucial pre- and post-processing
    corrections. We formalize this procedure, and demonstrate its theoretical
    guarantees in the form of sharp high probability upper bounds and minimax
    lower bounds. We illustrate our findings with detailed simulations.
\end{abstract}

\end{center}

\part{} 
\parttoc 

\clearpage
\section{Introduction}\label{sec:intro}

Isotonic regression is a classically studied nonparametric regression problem in
which the underlying signal satisfies a monotonicity constraint. In the
univariate case, this regression setup provides a flexible nonparametric
generalization to simple linear regression. That is, the underlying signal may
be non-linear, but still satisfies monotonicity as in the simple linear model.
The classically studied isotonic regression generating process is formally
described in \Cref{ndefn:isotonic-classical-model}: \bnumdefn[Classical isotonic
regression]\label{ndefn:isotonic-classical-model} We consider $n$ observations,
$\thesetb{Y_{i}}{i \in [n]}$, where each observation $Y_{i}$ is generated from
the following model:
\begin{align}
    Y_{i}
     & = \mu_{i} + \varepsilon_{i}
    \label{neqn:classical-iso-reg-1}         \\
    \text{s.t. } \mu_1
     & \leq \mu_{2} \leq \ldots \leq \mu_{n}
    \label{neqn:classical-iso-reg-2}         \\
    \text{ and } \varepsilon_{i}
     & \distiid \distNorm\parens{0,\sigma^2}
    \label{neqn:classical-iso-reg-3}
\end{align}
\enumdefn
The statistical goal under this classical setup is to estimate the underlying
signal vector $\boldsymbol{\mu} \defined \parens{\mu_1, \ldots,
\mu_{n}}^{\top}$, subject the monotonicity constraint in
\Cref{neqn:classical-iso-reg-2}, while $\sigma$ is an unknown (nuisance)
parameter. Throughout this paper we will adopt the convention, without loss of
generality, that the signal vector is monotonically increasing (as per
\Cref{neqn:classical-iso-reg-2}). Additionally we will assume that all
estimation errors are computed under square loss (in Euclidean metric), in high
probability.

\subsection{Adversarial sign-corrupted isotonic (\texttt{ASCI})
    regression}\label{subsec:asci-regression}

Our work in this paper is motivated by a variation of the classical isotonic
regression estimation problem, per \Cref{ndefn:isotonic-classical-model}. We
refer to this newly proposed model as \textit{adversarial sign-corrupted
isotonic} (\texttt{ASCI}) \textit{regression}. The generating process for this
\texttt{ASCI} estimation problem is formalized in
\Cref{ndefn:isotonic-adversarial-sign-model}.
\begin{restatable}[Adversarial sign-corrupted isotonic (\texttt{ASCI})
        regression]{ndefn}{ndefnisotonicadversarialsignmodel}\label{ndefn:isotonic-adversarial-sign-model}
        We consider $n$ observations, $\thesetb{\resp_{i}}{i \in [n]}$, where
        each observation $\resp_{i}$ is generated from the following model:
    \begin{align}
        \resp_{i}
         & = \xi_{i}(\mu_{i} + \varepsilon_{i})
        \label{neqn:adversarial-iso-reg-1}                         \\
        \text{s.t. } 0
         & < \eta \leq \mu_1 \leq \mu_{2} \leq \ldots \leq \mu_{n}
        \label{neqn:adversarial-iso-reg-2}                         \\
        \text{ and } \varepsilon_{i}
         & \distiid \distNorm\parens{0,\sigma^2}
        \label{neqn:adversarial-iso-reg-3}                         \\
        \text{ and } \xi_{i}
         & \in \theseta{-1, 1}
        \label{neqn:adversarial-iso-reg-4}
    \end{align}
\end{restatable}

\bnrmk\label{nrmk:eta-arbitrary} We note that the constant $\eta > 0$ is known
to the observer and adversary as part of the generating process. This means that
the true signal is positive, since it is uniformly bounded away from $\eta$. It
is an artefact of our method but as we will see later in
\Cref{nexa:asci-rademacher-mixture,nexa:equiv-conds}, highly non-trivial
estimation tasks are still contained under this constraint.
\enrmk

\bnrmk\label{nrmk:asci-defn-ref-name} Throughout this paper we will
interchangeably use the terms \texttt{ASCI} \textit{regression},
\textit{setting}, \textit{setup}, \textit{model}, and \textit{generating
process} to refer to \Cref{ndefn:isotonic-adversarial-sign-model}.
\enrmk

By comparing
\Cref{ndefn:isotonic-classical-model,ndefn:isotonic-adversarial-sign-model},
this \texttt{ASCI} regression generating process is a partial generalization of
the classical isotonic regression. It can be briefly described as follows. Here
the classical isotonic regression responses, $\mu_{i} + \varepsilon_{i}$ in
\Cref{neqn:classical-iso-reg-1}, are \emph{sign-corrupted} in a manner chosen by
an adversary, as captured by the multiplicative $\xi_{i}$ terms. Here the
$\xi_{i} \in \theseta{-1, 1}$ are \textit{sign-corruptions} for the
\textit{true} data generating process, \ie, $Y_{i} \defined \mu_{i} +
\varepsilon_{i}$. It is important to note that the $\xi_{i} \in \theseta{-1, 1}$
for each $i \in [n]$, are chosen given that the adversary has full access to the
true responses, \ie, $\thesetb{\mu_{i} + \varepsilon_{i}}{i \in [n]}$. As such,
\Cref{neqn:classical-iso-reg-1} in the classical isotonic regression setup
represents a special case of
\Cref{neqn:adversarial-iso-reg-1,neqn:adversarial-iso-reg-4} by taking $\xi_{i}
\eqas 1$ for each $i \in [n]$. However, we note that in this \texttt{ASCI}
setting, in \Cref{neqn:adversarial-iso-reg-2} the monotonically increasing
signal vector $\boldsymbol{\mu} \defined \parens{\mu_1, \ldots, \mu_{n}}^{\top}$
is bounded below by $\eta$, which is some fixed and known positive constant. As
such this represents a restriction of the classical isotonic regression
condition described in \Cref{neqn:classical-iso-reg-2}. In summary,
\texttt{ASCI} regression represents both a restriction and relaxation of the
classical isotonic regression generating process. We will see why the restriction is necessary in
this work, but we will later suggest possible ways in which it can be relaxed in
future work.

\subsection{Interesting special cases of \texttt{ASCI}
    regression}\label{subsec:asci-regression-special-cases}

Interestingly, we note that even some special cases of this \texttt{ASCI}
regression setup can result in highly non-trivial estimation tasks. Two
particular \texttt{ASCI} regression special cases are formalized in
\Cref{nexa:asci-rademacher-mixture,nexa:equiv-conds,nexa:equiv-conds}.

\begin{restatable}[Two-component Gaussian mixture \texttt{ASCI}  regression
        special
        case]{nexa}{nexaascisimulation}\label{nexa:asci-rademacher-mixture} We
        consider $n$ observations, $\thesetb{\resp_{i}}{i \in [n]}$, where each
        observation $\resp_{i}$ is generated from the following model:
    \begin{align}
        \resp_{i}
         & = \xi_{i}\mu_{i} + \varepsilon_{i}
        \label{neqn:simulations-01}                                             \\
        \text{s.t. } 0
         & < \eta \leq \mu_{1} \leq \mu_{2} \leq \ldots \leq \mu_{n}
        \label{neqn:simulations-02}                                             \\
        \text{ and } \varepsilon_{i}
         & \distiid \distNorm\parens{0, \sigma^2}
        \label{neqn:simulations-03}                                             \\
        \text{ and } \xi_{i}
         & \distiid \distRademacher{(p)}, \, p \in (0, 1), \text{ and } \xi_{i}
        \indep \varepsilon_{i}
        \label{neqn:simulations-04}
    \end{align}
\end{restatable}

\bnrmk\label{nrmk:asci-rademacher-mixture} Note that $\xi_{i} \distiid
    \distRademacher{(p)}$ for each $i \in [n]$, means that $\xi_{i} = +1$ with
    probability $p$, and $\xi_{i} = -1$ with probability $1 - p$. We note that
    the model defined in \Cref{nexa:asci-rademacher-mixture} is a special case
    of \Cref{ndefn:isotonic-adversarial-sign-model}. This is formally proved in
    \Cref{subsec:asci-rademacher-mixture-proofs}.
\enrmk

We note that in the univariate setting, \Cref{nexa:asci-rademacher-mixture}
represents a generalization of the two-component Gaussian mixture model studied
in detail in \citet[Section~3.2.1]{balakrishnan2017emalgostatanalysis}. Our
model generalizes their setting in the sense that we allow a different mean,
\ie, $\mu_{i}$, for each of the $n$ univariate observations. Interestingly, in
this more general univariate mixture setting, our proposed \texttt{ASCIFIT}
estimator (see \Cref{sec:ascifit-three-step-est-proc}) provides an efficient
alternative to the EM algorithm \citep{dempster1977emalgomleincompletedata}.
Such models have extensive applications, \eg, community detection
\citep{royer2017adaptclustsemidefinite,giraud2018partialrecboundsclust}.

\begin{restatable}[Non-convex \texttt{ASCI} regression special
        case]{nexa}{nexaequivconds}\label{nexa:equiv-conds} We consider $n$
        observations, $\thesetb{\resp_{i}}{i \in [n]}$, where each observation
        $\resp_{i}$ is generated from the following model:
    \begin{align}
        \resp_{i}
         & = \gamma_{i} + \varepsilon_{i}
        \label{neqn:equiv-conds-reg-1}                                                             \\
        \text{s.t. } 0
         & < \eta \leq \absa{\gamma_{1}} \leq \absa{\gamma_{2}} \leq \ldots \leq \absa{\gamma_{n}}
        \label{neqn:equiv-conds-reg-2}                                                             \\
        \text{ and } \varepsilon_{i}
         & \distiid \distNorm\parens{0,\sigma^2}
        \label{neqn:equiv-conds-reg-3}
    \end{align}
\end{restatable}
\bnrmk\label{nrmk:equiv-conds} Verifying that \Cref{nexa:equiv-conds} is a
special case of \Cref{ndefn:isotonic-adversarial-sign-model} is not a priori
obvious, and is formally proved in \Cref{subsec:equiv-conds-proofs}.
\enrmk

In \Cref{nexa:equiv-conds}, one can think of this model being generated from the
\texttt{ASCI} model as per \Cref{ndefn:isotonic-adversarial-sign-model}. In this
special case, the adversary randomly chooses sign-corruptions
\textit{independently} of the error terms, \ie, $\xi_{i} \indep \varepsilon_{i}$
for each $i \in [n]$. Under this setup, the resulting response term in
\Cref{neqn:equiv-conds-reg-1} is the same as the classical isotonic regression
response, as seen by comparing to \Cref{neqn:classical-iso-reg-1}. However,
interestingly the adversarial sign-corruption is now absorbed into the revised
monotonicity constraint $0 < \eta \leq \absa{\gamma_1} \leq \absa{\gamma_{2}}
\leq \ldots \leq \absa{\gamma_{n}}$, per \Cref{neqn:equiv-conds-reg-2}. As a
result, this generating process is a highly non-convex constrained estimation
problem. In this case \texttt{ASCIFIT} will allow one to recover $\vert
\gamma_{i} \vert$, where as \texttt{PAVA} will not provide \textit{any}
information on the $\vert \gamma_{i} \vert$ (or $\gamma_{i}$), given the
non-convex constraint.

\subsection{Motivation and focus of our
    work}\label{subsec:asci-motivation-focus-of-our-work} With the \texttt{ASCI}
    regression setup clearly defined, we turn our attention to describing the
    focus of our analysis in this work. This is summarized by the following core
    question of interest:
\begin{tcolorbox}
    \begin{itquote}
        \textbf{Core question:} Given the adversarial sign-corrupted isotonic
        \textnormal{(\texttt{ASCI})} regression setup in
        \Cref{ndefn:isotonic-adversarial-sign-model}, can we find a
        computationally efficient estimator for $\boldsymbol{\mu}$, and
        demonstrate its precise (non-asymptotic) statistical optimality?
    \end{itquote}
\end{tcolorbox}

To the best of our knowledge the \texttt{ASCI} model, and our core question of
interest, have not been explicitly studied before in the literature. We note
that this \texttt{ASCI} estimation problem is inherently challenging, and thus
interesting, for three main reasons, \ie, \ref{itm:asci-challenge-1} --
\ref{itm:asci-challenge-3}:
\benum[wide=0pt, label={\color{blue}\textbf{Challenge \Roman*}}, align=left, start=1]
\item \label{itm:asci-challenge-1} \textbf{(Dependent responses):} in this
estimation problem the adversary is free to choose the sign-corruption terms
$\xi_{i}$, \textit{after} observing all samples, possibly resulting in a strong
dependence between the original isotonic responses. As such, any \texttt{ASCI}
estimator must be able to handle arbitrary dependence structure between the
sign-corrupted responses.
\item \label{itm:asci-challenge-2} \textbf{(Violating signal monotonicity):}
qualitatively speaking, the sign-corruptions are in a sense `extreme' in that by
selectively changing the sign of the observations the adversary fundamentally
`attacks' the isotonic monotonicity constraint directly. It is this convex
monotone constraint which classical isotonic estimators, \ie, \texttt{PAVA}, are
\textit{designed} to exploit.
\item \label{itm:asci-challenge-3} \textbf{(Interesting special cases):} The
\texttt{ASCI} setting contains interesting non-trivial special cases as
described in \Cref{nexa:asci-rademacher-mixture,nexa:equiv-conds}. Naively
applying typical least squares estimation techniques will be unable to provide
any relevant information on the estimated quantity of interest.
\eenum

Given these three formidable challenges posed by \texttt{ASCI} regression, any
computationally and statistically efficient estimator here needs to utilize new
techniques to exploit the potential non-convex structure in our setting. Our
motivations here are thus driven by understanding the robustness of existing
isotonic regression estimators under such adversarial settings. Moreover, for
the \texttt{ASCI} setting to be worth studying, we wanted ensure practical
algorithms for estimation under this adversarial setting, with sharp minimax
(worst-case) statistical guarantees, both of which we were able to provide. We
thus view the \texttt{ASCI} setting as \textit{stimulating prototype} for more
such research into adversarial robustness in isotonic regression.


\subsection{Prior and related work}\label{sec:prior-related-work}

As noted, to the best of knowledge our core question of interest, \ie, isotonic
regression under the proposed \texttt{ASCI} setup, has not been previously
studied. Our work however builds on and utilizes known estimators from the
classical isotonic regression literature. As such we limit our prior and related work
summary on known risk bounds (and rates) for such isotonic regression
estimators, and the efficient algorithms (\ie, the \texttt{PAVA}) and practical
implementations thereof.

\vspace{4pt}
\nit \underline{\textbf{Isotonic regression (classical):}}

\vspace{4pt}
A lively historical overview of isotonic regression estimation from a
computational lens is given in \citet[Section~1]{deleeuw2009isotoneoptimpava}.
In brief, we note that the origins of isotonic regression can be traced back to
a number of independently written papers in the 1950s. In particular it was
studied by \citet{brunk1955mlemonotone,ayer1955empdistfnincompleteinfo}. Such
estimators for ``ordered parameters'' were also analyzed in the series of papers
\citet{vaneeden1956mlepartiallyord,vaneeden1957leastsquaresordparams,vaneeden1957mlepartiallyordparta,vaneeden1957mlepartiallyordpartb}
which culminated into a PhD thesis in by the same author
\citep{vaneeden1958estorderedparamsdistns}. Shortly thereafter the articles
\citep{bartholomew1959homogeneityparta,bartholomew1959homogeneitypartb} also
investigated the related idea of hypothesis testing under monotonicity
constraints. We refer the interested reader to the classical comprehensive
references
\citet{barlow1972statinfunderordrestrictions,robertson1988orderreststatinf}, for
further reading.

The classical isotonic regression setup per
\Cref{ndefn:isotonic-classical-model} under square loss is a convex optimization
problem. As such, it has a unique solution, \ie, the Euclidean projection onto
the closed convex monotone cone given by the constraint in
\Cref{neqn:classical-iso-reg-2}. In this case, the non-asymptotic risk bounds
for the least squares estimator (LSE) of the monotone parameters $\mu_{i}$ are
of the order $n^{- 2 / 3}$ in sample complexity. This convergence rate has been
established across a number of papers including
\citet{vandegeer1990estregressionfunc,vandegeer1993hellingerconsistencynpmles,
donoho1990gelfandnwidthsleastsquares,birge1993ratesconvmincontrastests,
wang1996l2riskisotonicestimate,meyer2000ondofshaperestregression,
zhang2002riskboundsisotonicreg,chatterjee2015onriskboundsisotonicreg}. Broadly
speaking, these results typically vary in the generality of their underlying
assumptions on the normality or independence of the error terms in classical
isotonic regression. As noted in the excellent recent survey
\cite{guntuboyina2018nonparshprstreg}, the same risk rate for this (and for more
general) LSEs was established using an alternative approach in
\cite{chatterjee2014newpersplstsqrsconvexconst}. Moreover, in the case of
minimax lower bounds, the matching risk rate (up to constant terms) for isotonic
regression was established in \citet{chatterjee2015onriskboundsisotonicreg} and
also in \citet{bellec2015sharporaclemonoconvexreg}, in both high probability and
expectation terms.

\vspace{4pt}
\nit \underline{\textbf{Pool Adjacent Violators Algorithm (\texttt{PAVA}):}}

\vspace{4pt}
Rather remarkably, despite the nonparametric setup of classical isotonic
regression, the LSE in this case has an explicit `max-min' formulation
\citep[Equation~(1.9)]{barlow1972statinfunderordrestrictions}. However, in
practice it is efficiently computed using the  pool adjacent violators algorithm
(\texttt{PAVA}). As described in \citet{tibshirani2011nearlyisoreg} the
\texttt{PAVA} in effect estimates the ordered parameters by scanning through the
(sorted) observations. For each adjacent pair of observations, the monotonicity
constraint is checked. If the constraint is `violated' by a given observation,
the average of the observations is used as the estimate, with appropriate
(minimal) backtracking to ensure that any restrospectively incurred  violations
are similarly corrected for. Efficient \texttt{PAVA} implementations, \eg, as
described in
\cite{grotzinger1984projontoordersimplexes,best1990activesetisoregunifframe},
have a computational complexity of $\bigoh{n}$, where $n$ is the sample size.
Since we will use the \texttt{PAVA} in just one step in our proposed three-step estimator
for the \texttt{ASCI} regression parameter $\boldsymbol{\mu}$, we will not
detail it further here. However, such open-source \texttt{PAVA} implementation
details can be found in
\citet{deleeuw2009isotoneoptimpava,pedregosa2011sklearnpython}.


\subsection{Main contributions}\label{sec:contributions} Our contributions in
this paper are twofold and are summarized as follows:

\bitems
\item \textbf{Computable estimators with non-asymptotic upper bounds:} We
propose a computationally efficient three-step algorithm \texttt{ASCIFIT}, to
estimate the required parameter $\boldsymbol{\mu}$, under the \texttt{ASCI}
setting. Our \texttt{ASCIFIT} estimator converges at a $n^{-2/3}$ rate, with
high probability. We illustrate our findings with extensive numerical
simulations.

\item \textbf{Sharp minimax lower bounds:} we provide matching high probability
lower bounds (up to constant and log factors) under square loss, and thus
demonstrate that our estimators are minimax optimal in this sense.
\eitems

In particular, our upper bound proofs involve rather subtle theoretical details
about the \texttt{PAVA}, and our use of method of moment techniques is quite
unique in this literature. We believe these proof techniques will be of
independent interest to researchers in isotonic regression. In particular, for
similar adversarial estimation tasks, where traditional convex M-estimation
techniques are infeasible.


\subsection{Organization of the paper}\label{subsec:organization}

The rest of this paper is organized as follows. In
\Cref{sec:ascifit-three-step-est-proc} we introduce \texttt{ASCIFIT}, our
three-step estimation procedure for $\boldsymbol{\mu}$. In
\Cref{sec:step-two-unique-root} we provide high probability upper bounds on
estimation rates using \texttt{ASCIFIT}. In \Cref{sec:lower-bounds} we establish
sharp minimax lower bounds for the parameter estimation in our \texttt{ASCI}
setting. In \Cref{sec:simulations} we provide extensive numerical \texttt{ASCI}
simulations, to illustrate our findings. In \Cref{sec:conclusion} we summarize
our results and describe exciting future research directions.

\subsection{Notation}\label{sec:notation}

Throughout this paper, we typically use lowercase for scalars in $\reals$, \eg,
$(x, y, z, \ldots)$, bold lowercase for vectors, \eg, $(\bfx, \bfy, \bfz,
\ldots)$, and bold uppercase for matrices, \eg, $(\bfX, \bfY, \bfZ, \ldots)$. We
use $\lesssim$ and $\gtrsim$ to mean $\leq$ and $\geq$, respectively, up to
positive universal constants. We denote $a \vee b \defined \max\braces{a, b}$
for each $a, b \in \reals$. We say that a sequence $a_{n} \defined \bigoh{1}$ if
there exists $C > 0, N \in \nats$ such that $\absa{a_{n}} < C$ for each $n > N$.
Similarly, $a_{n}=\bigoh{b_{n}}$ iff $\frac{a_{n}}{b_{n}} = \bigoh{1}$. We say
that a sequence $a_{n} = \litoh{1}$ if $a_{n} \rightarrow 0$ as $n \rightarrow
\infty$. Similarly, $a_{n}=\litoh{b_{n}}$ iff $\frac{a_{n}}{b_{n}} = \litoh{1}$.
We denote the finite set $\theseta{1, \ldots, n}$ by $[n]$. We define the
\textit{indicator function} $\Indb{\Omega}{\bfx}$ to take the value $1$ when
$\bfx \in \Omega \subseteq \reals^{d}$, and $0$ otherwise. We say that a
function $f : \Omega \to \reals$ is increasing, if for all $u, v \in \Omega
\subset \reals$ such that $u \leq v$, implies $f(u) \leq f(v)$. We use
\textit{strictly} increasing in the case where these inequalities are strict.
Similarly we note that $f$ is decreasing (or \textit{strictly} decreasing) when
these respective inequalities are reversed. We provide a useful notation summary
table in \Cref{subsec:app-notation-summary}.

\section{\texttt{ASCIFIT:} A three-step estimation procedure for
  \texorpdfstring{$\boldsymbol{\mu}$}{mu}}\label{sec:ascifit-three-step-est-proc}

As per our core question of interest, we now turn our attention to
\texttt{ASCIFIT}, \ie our proposed estimation procedure for $\boldsymbol{\mu}$,
under the \texttt{ASCI} setup. The \textit{Folded Normal} distribution, and in
particular its mean and variance, will be fundamental to \texttt{ASCIFIT}. As
such, we first formalize the key properties of the Folded Normal distribution in \Cref{ndefn:folded-normal}.
\begin{restatable}[Folded Normal distribution]{ndefn}{ndefnfoldednormal}
    \label{ndefn:folded-normal}
    Suppose $\resp \sim \distNorm(\mu, \sigma^{2})$, and let $\absresp \defined
        \absa{\resp}$. We then say that $\absresp \sim \distFoldNorm(\mu,
        \sigma)$, is a \emph{Folded Normal} distribution. We denote the mean and
    variance of $\absresp$, by $f(\mu, \sigma)$ and $g(\mu, \sigma)$,
    respectively. They are given as follows:
    \begin{align}
        \iftoggle{neurips2022sub}{f(\mu, \sigma)
            \defined \Exp{\absresp}
         & = \sigma \sqrt{2/\pi} \exp(- \mu^2/(2 \sigma^2))
            - \mu(1 - 2\Phi(\mu/\sigma))}{f(\mu, \sigma)
            \defined \Exp{\absresp}
         & = \sigma \sqrt{2/\pi} \exp(- \mu^2/(2 \sigma^2))
            - \mu(1 - 2\Phi(\mu/\sigma))}.
        \label{neqn:folded-normal-exp-1}                    \\
        g(\mu, \sigma)
        \defined \Var{\absresp}
         & = \mu^{2} + \sigma^{2} - f(\mu, \sigma)^{2}.
        \label{neqn:folded-normal-var}
    \end{align}
\end{restatable}
\bnrmk\label{nrmk:folded-normal-mean} We refer the reader to
\cite{tsagris2014foldednormdistn,elandt1961foldednormalestparmom} for more
details. We only consider $\mu > \eta > 0$ per
\Cref{neqn:adversarial-iso-reg-2}, and we use the shorthand notation $f(\mu,
    \sigma)^{2} \defined \parens{f(\mu, \sigma)}^{2}$.
\enrmk

We now describe \texttt{ASCIFIT}, our three-step procedure to estimate
$\boldsymbol{\mu}$ under the \texttt{ASCI} setting, as follows:

\vspace{8pt}
\begin{mycolorbox}{blue}{\textbf{\texttt{ASCIFIT}}: Three-step procedure to
        estimate $\boldsymbol{\mu}$ under the \texttt{ASCI} setting}
    \benum[wide=0pt, label={\color{blue} \textbf{Step \Roman*}}, align=left, start=1]

    \item \label{itm:asci1} \textbf{(\textit{Pre-processing and
        }\texttt{PAVA}):} \nit \vspace{1mm} \newline Obtain an initial
    \textit{naive} estimate of $\boldsymbol{\mu} \defined \parens{\mu_1,
            \ldots, \mu_{n}}^{\top}$ by fitting isotonic regression (using the
    \texttt{PAVA}) on $\absresp_{i} \defined \absa{\resp_{i}}$. Denote these
    estimates by $\boldsymbol{\widehat{\mu}_{\text{naive}}} \defined
        (\widehat{\absresp}_{1}, \ldots, \widehat{\absresp}_{n})^{\top}$.
    \vspace{2mm}
    \item \label{itm:asci2} \textbf{(\textit{Second moment matching}):} \nit
    \vspace{1mm}
    \newline
    Estimate $\sigma$ in the following way. Pick the $\sigma$ solving the
    following equation, and denote the corresponding solution as
    $\widehat{\sigma}$:
    \begin{equation}\label{neqn:second-moment-matching-step}
        G(\sigma)
        \defined \sigma^2 +
        \frac{1}{n} \sum_{i = 1}^{n} (f^{-1}(\widehat{\absresp_{i}} \vee f(\eta, \sigma), \sigma))^2
        = \frac{1}{n} \sum_{i = 1}^{n} \absresp_{i}^2.
    \end{equation}
    Here $f^{-1}(\cdot, \sigma)$, denotes the inverse function of $f(\mu,
        \sigma)$ with respect to $\mu$, when we hold $\sigma$ fixed to the value
    $\sigma$.
    \vspace{2mm}
    \item \label{itm:asci3} \textbf{(\textit{Post-processing via plug-in}):}
    \nit \vspace{1mm}
    \newline
    From $\boldsymbol{\widehat{\mu}_{\text{naive}}}$ in \ref{itm:asci1}, and
    $\widehat{\sigma}$ in \ref{itm:asci2}, compute
    $\muhatascifit \defined
        \parens{\widehat{\mu}_{1}, \ldots, \widehat{\mu}_{n}}^{\top}$ as follows:
    \begin{equation}\label{neqn:correction-step}
        \widehat{\mu}_{i} \defined f^{-1}(\widehat{\absresp_{i}} \vee f(\eta, \widehat{\sigma}), \widehat{\sigma}), \, \text{for each $i \in [n]$}.
    \end{equation}
    \eenum
\end{mycolorbox}

\subsection{Intuition for the three \texttt{ASCIFIT}
    steps}\label{sec:ascifit-three-step-est-procs}

We now provide more precise intuition for each of the three \texttt{ASCIFIT}
steps, \ie, \ref{itm:asci1} -- \ref{itm:asci3}.

\vspace{4pt}
\nit \underline{\textbf{Intuition for} \ref{itm:asci1}:}

\vspace{4pt}
Here, we begin with the
pre-processing operation $\absresp_{i} \defined \absa{\resp_{i}}$. This serves
the critical dual purpose of removing the effect of the sign-corruptions
$\xi_{i}$, and also induces \emph{independence} of the resulting observations
$(\absresp_{1}, \ldots, \absresp_{n})^{\top}$. This helps directly address
\ref{itm:asci-challenge-1} and \ref{itm:asci-challenge-2} under the
\texttt{ASCI} setup. To better understand this dual effect, note that in the
\texttt{ASCI} setup, the $\xi_{i} \in \theseta{-1, 1}$ may be arbitrarily chosen
by the adversary (without a precise distributional assumption). However, the
critical information in our model is given by pre-processing each observation,
$R_{i}$, as $\absresp_{i} \defined \absa{\resp_{i}}$. More specifically we have
that $\absresp_{i} = \absa{\xi_{i}\parens{\mu_{i} + \varepsilon_{i}}} =
    \absa{\mu_{i} + \varepsilon_{i}}$. Since $\mu_{i} + \varepsilon_{i} \distinid
    \distNorm\parens{\mu_{i}, \sigma^{2}}$, per \Cref{ndefn:folded-normal} we have
that $\absresp_{i} \distinid \distFoldNorm\parens{\mu_{i}, \sigma}$, per \Cref{ndefn:folded-normal}. We note
that our pre-processed observations $\theset{\absresp_{1}, \ldots \absresp_{n}}$
are all \inid\footnote{\ie, independent but not identically distributed.}, since
they have a common variance $\sigma^{2}$ but varying means $\mu_{i}$ for each
observation $i \in [n]$. Moreover, fitting an isotonic regression to
$\absresp_{i}$ intuitively estimates $f(\mu_{i}, \sigma)$ which are the mean of
each $\absresp_{i}$. This step is formally justified by the results of
\citet{zhang2002riskboundsisotonicreg}.

\vspace{4pt}
\nit \underline{\textbf{Intuition for} \ref{itm:asci2}:}

\vspace{4pt}
This is motivated by
second moment matching to estimate $\sigma$. Specifically, using the fact that
the expected value of $\frac{1}{n} \sum_{i = 1}^n \absresp_{i}^2$, is $\sigma^2
    + \frac{1}{n} \sum_{i = 1}^n \mu_{i}^2$. The left hand side of
\Cref{neqn:second-moment-matching-step} directly estimates the term $\sigma^{2}
    + \frac{1}{n} \sum_{i = 1}^n \mu_{i}^2$. In \ref{itm:asci2} it is not clear
a priori whether such an inverse function $f^{-1}(\cdot, \sigma)$ exists, or
whether there exists a unique positive solution for $\sigma$ in
\Cref{neqn:second-moment-matching-step}. We will demonstrate that both
assertions are true, and that the unique solution $\widehat{\sigma}$, to
estimate $\sigma$, can be computed efficiently with a binary search approach. We
would like to note here that estimating $\sigma$ is not an easy problem (it is
not by coincidence that in classical isotonic regression that $\sigma$ is viewed
as a nuisance parameter). This difficulty explains why we need to impose some
additional assumptions on the vector $\mu$ and on $\sigma$ later on. Next, we
provide the intuition on why we use the factor $\widehat{\absresp_{i}} \vee
    f(\eta, \sigma)$ in \ref{itm:asci2}, for each $i \in [n]$. This is summarized in
\Cref{nprop:intuition-for-correction}.
\begin{restatable}[Reason for the ``$\vee f(\eta, \sigma)$''-correction in \ref{itm:asci2}]{nprop}{npropintuitionforcorrection}
    \label{nprop:intuition-for-correction}
    The need for defining the $\vee f(\eta, \sigma)$ in
    \Cref{neqn:second-moment-matching-step} in \textnormal{\ref{itm:asci2}} in
    \texttt{ASCIFIT}, is that the solution to the problem
    \begin{align}
        \argmin_{\widetilde{\absresp}_{1}, \ldots, \widetilde{\absresp}_{n}} \sum_{i = 1}^{n} (\absresp_{i} - \widetilde{\absresp}_{i})^2
         & \mbox{ s.t. }
        f(\eta, \sigma) \leq \widetilde{\absresp}_{1} \leq \ldots \leq \widetilde{\absresp}_{n},
        \label{neqn:intuition-for-correction-01}
    \end{align}
    is related to the solution to
    \begin{align}
        \argmin_{\widehat{\absresp}_{1}, \ldots, \widehat{\absresp}_{n}} \sum_{i = 1}^{n} (\absresp_{i} -  \widehat{\absresp}_{i})^2
         & \mbox{ s.t. }
        \widehat{\absresp}_{1} \leq \ldots \leq \widehat{\absresp}_{n},
        \label{neqn:intuition-for-correction-02}
    \end{align}
    as $\widetilde{\absresp}_{i} \defined \widehat{\absresp}_{i} \vee f(\eta,
        \sigma)$.
\end{restatable}

To understand the significance of \Cref{nprop:intuition-for-correction}, first
note that we apply the \texttt{PAVA} to the $\absresp_{i}$ values in
\ref{itm:asci1}. As such, the corresponding least squares \texttt{PAVA}
estimates, $\widehat{\absresp}_{i}$, actually project onto the
\textit{unconstrained} monotone cone, as per
\Cref{neqn:intuition-for-correction-02}. However, in our setup we actually want
to solve the \textit{constrained non-negative} monotone means, as per
\Cref{neqn:intuition-for-correction-01}. Fortunately, this is not an issue since
we can simply post hoc correct each of the fitted unconstrained \texttt{PAVA}
solutions as $\widetilde{\absresp}_{i} \defined \widehat{\absresp}_{i} \vee
    f(\eta, \sigma)$, for each $i \in [n]$. This follows by adapting
\citet[Corollary~1]{nemeth2012projectmonotonenonnegconepava} to our
\texttt{ASCIFIT} setup. From all of the above discussion, intuitively it follows
that the term $\sigma^2 + \frac{1}{n} \sum_{i = 1}^{n}
    (f^{-1}(\widehat{\absresp_{i}} \vee f(\eta, \sigma), \sigma))^2$ in
\eqref{neqn:second-moment-matching-step} also estimates $\sigma^{2} +
    \frac{1}{n} \sum_{i = 1}^n \mu_{i}^2$, which explains why in \ref{itm:asci2} we
equate that term to $\frac{1}{n} \sum_{i = 1}^n \absresp_{i}^2$.

\vspace{4pt}
\nit \underline{\textbf{Intuition for} \ref{itm:asci3}:}

\vspace{4pt}
To understand the need
for this step, one needs to realize that the \texttt{PAVA} will estimate the
means of $\absresp_{i}$ which are $f(\mu_{i}, \sigma)$. Hence in order for us to
go back at the original $\mu_{i}$ scale, we need to invert the value of the PAVA
estimates $\widehat{\absresp}_{i}$. Ideally we would use the true value of
$\sigma$ in the inversion, but since it is unavailable to us, we use the plug-in
estimate $\widehat{\sigma}$ as computed in \ref{itm:asci2}. In addition the term
``$\vee f(\eta, \widehat{\sigma})$'' in \Cref{neqn:correction-step} is present,
since by assumption the value of each $\mu_{i}$ (or sufficiently just $\mu_{1}$)
must be at least $\eta$, after inverting.

\section{Analysis of \texttt{ASCIFIT}: Upper
  bounds}\label{sec:step-two-unique-root}

We have now described the details and key intuition behind our three-step
\texttt{ASCIFIT} estimator $\muhatascifit$, for $\boldsymbol{\mu}$. We now turn
our attention to formalizing this intuition into least squares estimation risk
bounds. More specifically, our end goal in this section is to describe our high
probability non-asymptotic upper risk bound for $\muhatascifit$, and understand
its dependence on the sample complexity, and other \texttt{ASCI} parameters. We
also provide summary sketch behind the main proof techniques used and what
insight they offer for estimation purposes. Before we state the results we will
define the rate of convergence $r_{n,2}(\mu_{n}, \mu_{1}, \sigma)$, which plays
an important role in all of the Theorems to follow. For an absolute constant
$C_2 > 0$ define
\begin{align}
    r_{n,2}(\mu_{n}, \mu_{1}, \sigma)
    \defined \min \bigg[2 \sigma^{2} C_{2}^{2},
        \frac{27}{4}\parens{\frac{\mu_{n} - \mu_{1}}{n}}^{\frac{2}{3}} (\sigma C_{2})^{\frac{4}{3}}
        + \frac{2 \sigma^{2} C_{2}^{2}}{n} \parens{1 + \log{n}}\bigg]. \label{neqn:upp-bound-Rn2-squared-main}
\end{align}

Importantly, assuming that $\mu_{n} - \mu_{1}, \sigma$ are constants not scaling
with the sample size $n$, we have that $r_{n,2}(\mu_{n}, \mu_{1}, \sigma)
\lesssim \max\big\{\parens{\frac{\sigma^{2} V}{n}}^{\frac{2}{3}},
\frac{\sigma^{2} \log{n}}{n}\}$, where $V \defined \mu_{n} - \mu_{1}$, is the
total variation of the underlying monotone signal. With this essential
background, we are ready to state our first result in
\Cref{nthm:eqn-has-unique-root}.

\begin{restatable}[\Cref{neqn:second-moment-matching-step} has a unique
        root]{nthm}{nthmeqnhasuniqueroot}\label{nthm:eqn-has-unique-root} Assume
        that there exist constants $\univlowerboundsigma, \univupperboundsigma,
        \univupperboundmu > 0$ such that $\univlowerboundsigma \leq \sigma \leq
        \univupperboundsigma$ and $\frac{1}{n}\sum_{i = 1}^n \mu_i^2 \leq
        \univupperboundmu$, for each $n \in \nats$. In addition let
        $r_{n,2}(\mu_{n}, \mu_{1}, \sigma) = o(1)$, where the quantity
        $r_{n,2}(\mu_{n}, \mu_{1}, \sigma)$ is defined in
        \eqref{neqn:upp-bound-Rn2-squared-main}. Then for sufficiently large
        $n$, $\delta = \litoh{\parens{r_{n,2}(\mu_n, \mu_1, \sigma)}^{-1}}$, and
        $\gamma = \litoh{n^{1 / 2}}$, under the \texttt{ASCI} setup,
        \Cref{neqn:second-moment-matching-step} in \texttt{ASCIFIT} has a unique
        root $\sigma^{*} \in \brackets{0,\sqrt{\frac{1}{n} \sum_{i = 1}^{n}
        \absresp_{i}^2}}$ for $\sigma$ with probability at least $1 -
        \delta^{-1} - 2\gamma^{-2}$.
\end{restatable}

The key insight of \Cref{nthm:eqn-has-unique-root} from a statistical
perspective, is that our second moment matching approach in \ref{itm:asci2} will
ensure that our proposed estimator $\widehat{\sigma}$, for $\sigma$, will be
unique with high probability. The core idea behind the proof of
\Cref{nthm:eqn-has-unique-root} is that the map $\sigma \mapsto G(\sigma)
\defined \sigma^2 + \frac{1}{n} \sum_{i = 1}^n (f^{-1}(\widehat{\absresp_{i}}
\vee f(\eta, \sigma), \sigma))^2$ is monotone increasing over $\sigma \geq 0$.
This enables the use of the intermediate value theorem to check that two
endpoints of $G(\sigma) - \frac{1}{n}\sum_{i = 1}^n T_i^2$, evaluated at $\sigma
\in \theset{0,\sqrt{\frac{1}{n} \sum_{i = 1}^{n} \absresp_{i}^2}}$ have opposite
sign with high probability. This has important practical implications for
estimation purposes. In effect, it means that we can efficiently
\textit{compute}
$\widehat{\sigma}$, by solving $G(\sigma) = \frac{1}{n} \sum_{i = 1}^{n}
    \absresp_{i}^2$ (per \Cref{neqn:second-moment-matching-step}) using a binary
    search approach between the two identified endpoints. We would like to
    mention that while using the intermediate value theorem sounds like an easy
    task, it turns out that it is extremely hard to verify that $G(0) \leq
    \frac{1}{n}\sum_{i = 1}^n T_i^2$, for which the bulk of the proof of Theorem
    \ref{nthm:eqn-has-unique-root} is dedicated to.

Although \Cref{nthm:eqn-has-unique-root} gives us a high probability bound on
estimating $\widehat{\sigma}$ uniquely, it is important to next understand how
efficiently $\widehat{\sigma}$ estimates $\sigma$. This is summarized in
\Cref{nthm:sigma-hat-close-sigma}.

\begin{restatable}[$\widehat{\sigma}$ is close to
        $\sigma$]{nthm}{nthmsigmahatclosesigma}\label{nthm:sigma-hat-close-sigma}
        Under the assumptions of \Cref{nthm:eqn-has-unique-root}, we have that $
        \absa{\sigma - \hat \sigma} \lesssim (\delta r_{n,2}(\mu_{n}, \mu_{1},
        \sigma) )^{1/2} + \gamma n^{-1/2}$ with probability at least $1-
        \delta^{-1}- 2\gamma^{-2}$, where $\delta^{-1}, \gamma^{-2}\in (0,1)$.
\end{restatable}

From \Cref{nthm:sigma-hat-close-sigma} we see that $\widehat{\sigma}$ converges
to $\sigma$ roughly at a $n^{-1 / 3}$ rate. In both
\Cref{nthm:eqn-has-unique-root,nthm:sigma-hat-close-sigma}, we require that
there exist constants $\univlowerboundsigma, \univupperboundsigma,
\univupperboundmu > 0$ such that $\univlowerboundsigma \leq \sigma \leq
\univupperboundsigma$ and $\frac{1}{n}\sum_{i = 1}^n \mu_i^2 \leq
\univupperboundmu$, for each $n \in \nats$. For transparency, we note that such
assumptions are an artefact of our methodology and ensure that our risk bounds
can be tightly controlled using the second moment matching approach. Given the
highly adversarial corruptions and non-convex constraints that can arise under
\texttt{ASCI} estimation, \eg, in \Cref{nexa:equiv-conds}, these are slightly
stronger assumptions required for classical convex isotonic regression setup.
They effectively represent a trade-off for the flexibility, and simplicity of
using \texttt{ASCIFIT} under these adversarial settings, whilst still ensuring
precise control in the parameter estimation risk bounds.

$\widehat{\sigma}$ in our post-processing correction for $\muhatascifit \defined
\parens{\widehat{\mu}_{1}, \ldots, \widehat{\mu}_{n}}^{\top}$. That is, our
final estimate for each $\mu_{i}$, is given by $\widehat{\mu}_{i} \defined
f^{-1}(\widehat{\absresp_{i}} \vee f(\eta, \widehat{\sigma}),
\widehat{\sigma})$. With this explicit form and our tightly controlled bounds in
\Cref{nthm:eqn-has-unique-root} and \Cref{neqn:second-moment-matching-step} we
are finally able derive the required least squares risk rate of $\muhatascifit$.
This is summarized in \Cref{nthm:mu-hat-close-mu}. We will shortly discuss this
result further in \Cref{sec:lower-bounds} when we derive high probability
minimax lower bounds.

\begin{restatable}[$\muhatascifit$ is close to
        $\boldsymbol{\mu}$]{nthm}{nthmmuhatclosemu}\label{nthm:mu-hat-close-mu}
        Under the assumptions of \Cref{nthm:eqn-has-unique-root} and
        \Cref{nthm:sigma-hat-close-sigma}, we have that
    \begin{align}\label{neqn:mu-hat-close-mu-main-01}
        \frac{1}{n} \sum_{i = 1}^{n}(f^{-1} (\widehat{\absresp_{i}} \vee f(\eta, \hat \sigma), \hat \sigma) - \mu_{i})^2
         & \lesssim \delta r_{n,2}(\mu_{n}, \mu_{1}, \sigma)  + \gamma^2 n^{-1},
    \end{align}
    with probability at least $1 - \delta^{-1} - 2\gamma^{-2}$.
\end{restatable}

\bnrmk
We note that $\eta$ is currently absorbed in our constants in
\Cref{nthm:eqn-has-unique-root,nthm:sigma-hat-close-sigma,nthm:mu-hat-close-mu}.
The exact form is complicated (but the smaller the $\eta$ the bigger the
constants). For more details, please refer to
\Cref{sec:step-two-unique-root-proofs}.
\enrmk

\section{Lower bounds}\label{sec:lower-bounds}

We now derive high probability minimax lower bounds under the \texttt{ASCI}
setting. We accordingly first introduce the relevant related notation and
definitions here for classes of monotonic sequences. We denote $\mclS^{\uparrow}
  \defined \thesetb{\boldsymbol{\mu} \defined (\mu_1, \ldots,
    \mu_{n})^{\top}}{\mu_1 \leq \ldots \leq \mu_{n}}$ to be the set of all
non-decreasing sequences. We define $k(\boldsymbol{\mu}) \geq 1$, for
$\boldsymbol{\mu} \in \mclS^{\uparrow}$, to be the integer such that
$k(\boldsymbol{\mu})-1$ is the number of inequalities $\mu_{i} \leq \mu_{i+1}$
that are strict for $i \in [n-1]$ (\ie, the number of `jumps' of
$\boldsymbol{\mu}$). The class of bounded monotone functions are
$\mclS^{\uparrow}(V^{*}) \defined \thesetb{\boldsymbol{\mu} \in
    \mclS^{\uparrow}}{V(\boldsymbol{\mu}) \leq V^{*}}$, for some fixed $V^{*} \geq
  0$, and $V(\boldsymbol{\mu}) \defined \mu_{n} - \mu_{1}$, is the total variation
of any $\boldsymbol{\mu} \in \mclS^{\uparrow}$. We focus on the
\texttt{ASCI}-restricted class of monotone sequences, \ie,
$\mclS^{\uparrow}(V^{*}, \eta, \univupperboundmu) \defined
  \thesetb{\boldsymbol{\mu} \in \mclS^{\uparrow}(V^{*})}{\frac{1}{n}\sum_{i =
      1}^{n} \mu_{i}^{2} \leq \univupperboundmu, \mu_{1} > \eta > 0}$.

We closely follow the approach of
\citet[Proposition~4]{bellec2015sharporaclemonoconvexreg} but non-trivially adapt it
to our \texttt{ASCI} setting by ensuring the monotonicity constraint in
\Cref{neqn:adversarial-iso-reg-2} is satisfied in the lower bound construction.
The proof uses well established techniques including the Varshamov-Gilbert bound
\citet[Lemma~2.9]{tsybakov2009intrononparmestimation}, and Fano's Lemma
arguments using \citet[Theorem~2.7]{tsybakov2009intrononparmestimation}. This
leads to our minimax lower bound result in \Cref{nprop:minimax-lower-bounds-final}.

\begin{restatable}[Minimax lower
    bounds]{nprop}{ncorminimaxlowerbounds}\label{nprop:minimax-lower-bounds-final} Let
  $n \geq 2, V^{*} > 0$ and $\sigma > 0$, and define
  $\widetilde{r}_{n,2}(V^{*}, \sigma) \defined
    \max\big\{\parens{\frac{\sigma^{2} V^{*}}{n}}^{\frac{2}{3}},
    \frac{\sigma^{2}}{n}\}$. Then, there exist absolute constants $c,
    c^{\prime}>0$ such that:
  \begin{equation}\label{neqn:minimax-lower-bounds-02}
    \inf_{\hat{\boldsymbol{\mu}}}
    \sup_{\substack{\text{$\mclS^{\uparrow}(V^{*}, \eta, \univupperboundmu)$}}}
    \Prba{\boldsymbol{\mu}}{\frac{1}{n} \norma{\widehat{\boldsymbol{\mu}} -
        \boldsymbol{\mu}}^{2} \geq \\ c \widetilde{r}_{n,2}(V^{*}, \sigma)} > c^{\prime}
  \end{equation}
\end{restatable}

Crucially, \Cref{nprop:minimax-lower-bounds-final} demonstrates that our high probability upper bounds for
\texttt{ASCIFIT} in \Cref{nthm:mu-hat-close-mu} are sharp in the minimax sense, up to constants and $\log$
factors. This is evident by directly comparing $\widetilde{r}_{n,2}(V^{*},
  \sigma)$ to $r_{n,2}(\mu_{n}, \mu_{1}, \sigma)$ per
\Cref{neqn:upp-bound-Rn2-squared-main}.

\section{Simulations}\label{sec:simulations}

\begin{figure}[!h] \centering
    \includegraphics[width=0.8\linewidth]{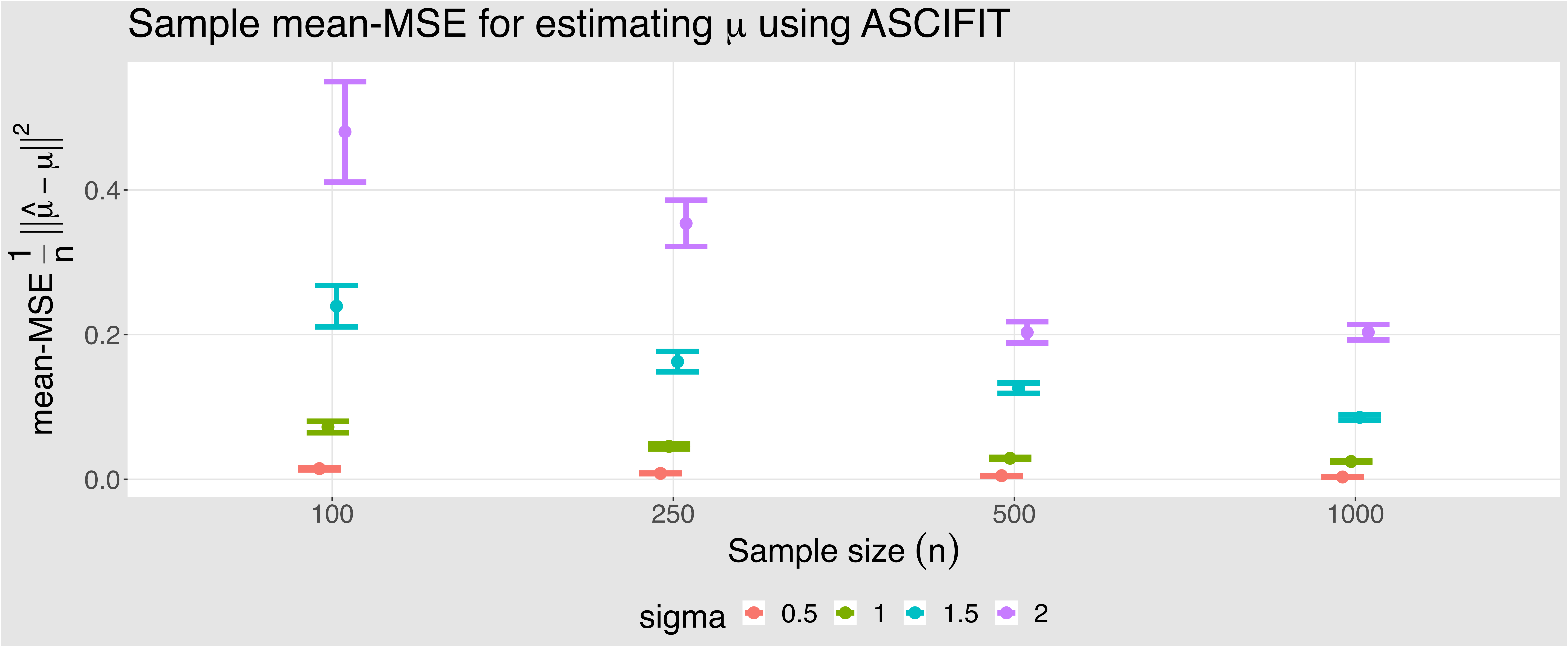}
    \caption{Mean (sample) MSE of \texttt{ASCIFIT} as a function of $n,
    \sigma$.} \label{fig:mean-mse-plot-01}
\end{figure}

We now demonstrate our \texttt{ASCIFIT} estimation algorithm in action through a
variety of simulations\iftoggle{anonymoussub}{}{\footnote{Reproducible code for
all figures in this paper is found at:
\url{https://github.com/shamindras/ascifit}. All of the simulation results in
this section were run on a personal Macbook laptop with macOS, Intel Core i9
CPU, and 64GB RAM. The total runtime for a single run of all simulations is
approximately 90 minutes.}}. Specifically, for simulation purposes we consider
$n$ observations, $\thesetb{\resp_{i}}{i \in [n]}$, where each observation
$\resp_{i}$ is generated from \Cref{nexa:asci-rademacher-mixture}. For
sufficiently large $n$, the \texttt{ASCI} model in
\Cref{nexa:asci-rademacher-mixture} roughly translates to $(1 - p)$-proportion
of observations being independently sign-corrupted by the adversary. Moreover we
assume per \Cref{neqn:simulations-04} that the adversarial sign-corruptions,
$\xi_{i}$, are chosen independently of all true errors, $\varepsilon_{i}$, for
each $i \in [n]$. The true monotone signal is defined to be $\mu_{i} \defined
\eta + \parens{1 - \eta} \frac{i - 1}{n}$, for each $i \in [n]$. We run this
generating process over the following parameter grid: $\eta \defined 0.2$, $p
\defined 0.5$, $\sigma \in \theset{0.5, 1, 1.5, 2}$, $n \in \theset{100, 250,
500, 1000}$. We perform 50 replications for each combination of simulation grid
parameters. In each replication of this generating process we fit the
\texttt{ASCIFIT} estimator $\boldsymbol{\widehat{\mu}_{\text{ascifit}}}$, for
$\boldsymbol{\mu}$. The main summary result from running our simulation, is
shown in \Cref{fig:mean-mse-plot-01}.

To clarify, given $\eta = 0.2, p = 0.5$, \Cref{fig:mean-mse-plot-01} plots the
sample mean-MSE, $\frac{1}{n} \norma{\muhatascifit - \boldsymbol{\mu}}^{2}$,
over 50 \texttt{ASCIFIT} replications for each value of $n \in \theset{100, 250,
500, 1000}$. Here the sample mean-MSE is a useful simulation proxy for the least
squares error, our core theoretical risk measure of interest. This sample
mean-MSE is plotted separately for each of the four sigma values, $\sigma \in
\theset{0.5, 1, 1.5, 2}$. The mean-MSE value of each replication ($\pm$ 2
standard errors) are shown using error bars in an effort to quantify replication
uncertainty. The plot in \Cref{fig:mean-mse-plot-01} is as expected in that all
of the sample mean-MSE values show a steady decreasing trend in $n$. Importantly
the relative uncertainty in sample mean-MSE reduces in $n$, as seen by the
smaller error bars to the right of \Cref{fig:mean-mse-plot-01}. For smaller
$\sigma$ values, \ie, smaller variance in the underlying generating model, we
see a much lower sample mean-MSE on average compared to higher $\sigma$-valued
simulations. That is, our \texttt{ASCIFIT} estimator achieves better accuracy,
with smaller underlying variability in the model, on average when other factors
are held constant.

\begin{figure}[!h] \centering
    \includegraphics[width=0.8\linewidth]{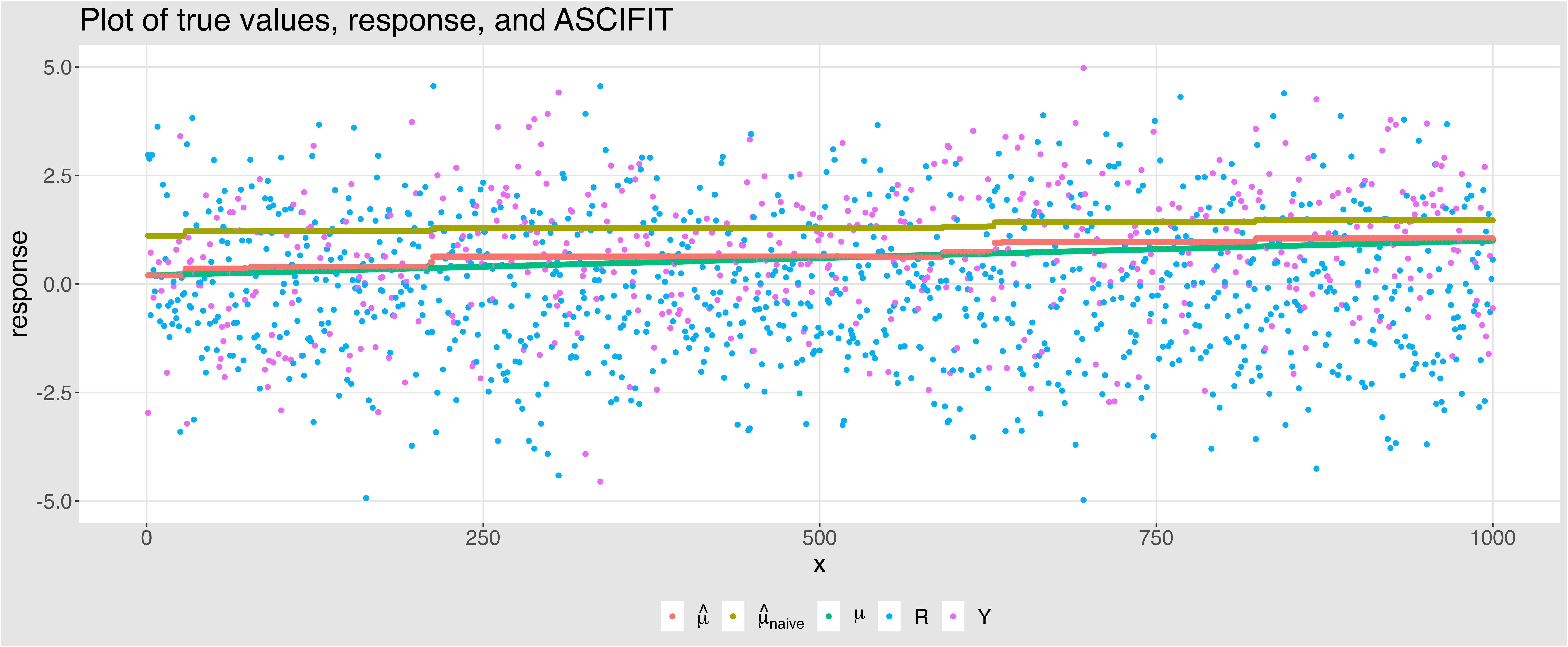}
    \caption{Mean (sample) MSE of \texttt{ASCIFIT} as a function of $n,
    \sigma$.} \label{fig:mean-mse-plot-02}
\end{figure}

Finally, in order to precisely gauge how well the \texttt{ASCIFIT} estimator
$\boldsymbol{\widehat{\mu}_{\text{ascifit}}}$, actually fits the true signal
$\boldsymbol{\mu}$, it is instructive to plot both directly on the original
generating sample data. This is seen for one instance over our parameter grid of
simulations in \Cref{fig:mean-mse-plot-02}.

Specifically, for $\eta = 0.2, p = 0.5, n = 1000, \sigma = 1.5$,
\Cref{fig:mean-mse-plot-02} plots the simulated true generating process,
$\boldsymbol{\mu}$, against the \texttt{ASCIFIT} estimator, $\muhatascifit$.
Additionally both the original and sign-corrupted individual observations are
plotted to emphasize the difficulty of this estimation task. Moreover, for
comparison purposes we also plot the naive estimator, \ie, $\muhatnaive$. Here
$\muhatnaive$ represents the estimator by stopping at \ref{itm:asci1} in
\texttt{ASCIFIT}. That is, estimating $\boldsymbol{\mu}$, by simply fitting
isotonic regression (using the \texttt{PAVA}) on $\absresp_{i} \defined
\absa{\resp_{i}}$. Furthermore, since $p = 0.5$, as expected, on average roughly
half of the true responses are adversarially sign-corrupted. Despite this, one
can see that \texttt{ASCIFIT} is relatively stable and reasonably recovers the
true signal. This shows more directly (in such an instance), the robustness of
\texttt{ASCIFIT} under such randomized adversarial sign-corruptions. Moreover
since $n = 1000$, we can see that \texttt{ASCIFIT} indeed fits well with
increasing sample complexity. In addition it highlights the importance of
\ref{itm:asci2} and \ref{itm:asci3} in \texttt{ASCIFIT}.



\section{Conclusion}\label{sec:conclusion}

In this paper we have considered a variation of the original isotonic regression
problem in which the observations can be adversarially corrupted in their sign
value. In this \texttt{ASCI} setting, adversarially refers to the fact that the
sign-corruptions can be chosen to have strong dependence with the error terms in
the original model. Our simple three-step estimation procedure,
\texttt{ASCIFIT}, is easy to implement with existing software and has sharp
non-asymptotic minimax guarantees on the estimation error, under square loss.
For future directions we note that that true signal is required to be strictly
positive for our guarantees to hold. We believe this restriction can be lifted
if one uses \textit{unimodal} regression instead of isotonic regression in
\ref{itm:asci1}. However, sharp risk guarantees need to first be proven similar
to \cite{zhang2002riskboundsisotonicreg} under this unimodal setting. It would
also be interesting to see if the moment matching technique could be extended
subgaussian error terms. We leave these exciting directions for future
work.\looseness=-1

\section{Acknowledgments}\label{sec:acknowledgments}
We would like to thank Arun Kumar Kuchibhotla, Alex Reinhart, Alessandro
Rinaldo, Larry Wasserman from the Carnegie Mellon University (CMU) Department of Statistics \& Data Science,
and Yang Ning from the Cornell Department of Statistics \& Data Science. Their
encouragement, and extensive feedback throughout this work greatly shaped the
final outcome. This paper extensively utilizes the \texttt{R} statistical
software \citep{rcoreteam2021cran} for conducting simulations and plots. In
particular, we relied primarily on the \texttt{tidyverse}
\citep{wickham2019tidyversejoss} collection of \texttt{R} packages.

\clearpage
\bibliographystyle{apalike}
\bibliography{refs}

\clearpage
\appendix
\addcontentsline{toc}{section}{Appendix} 
\part{Appendix} 
\parttoc 
\section{Mathematical Preliminaries}\label{sec:app-mathematical-preliminaries}

In this appendix we provide detailed proofs of all key statements from the main
paper. Since our work relies a variety of core ideas from isotonic regression we
first introduce some common definitions which will be referred to in subsequent
proofs.

\subsection{Notation Summary}\label{subsec:app-notation-summary}

To ensure that the Appendix is can be read in a standalone manner, we
consolidate key notation used in the paper in \Cref{tab:notation}. Unless stated
otherwise $K \subseteq \reals^{d}$ is a closed, non-empty convex set, and
$\Omega \subseteq \reals^{d}$.

\begin{table}[htbp]\caption{Notation and conventions used in our paper}
    \centering 
    \begin{tabular}{r l p{10cm} }
        \toprule
        \multicolumn{2}{c}{\underline{\textbf{Variables and inequalities}}}
        \\
        \multicolumn{2}{c}{}
        \\
        $a \wedge b$                                       & $\min\braces{a, b}$
        for each $a, b \in \reals$
        \\
        $a \vee b$                                         & $\max\braces{a, b}$
        for each $a, b \in \reals$
        \\
        scalars                                            & $x, y, z \in
        \reals$
        \\
        vectors                                            & $\bfx, \bfy, \bfz
        \in \reals^{d}$
        \\
        matrices                                           & $\bfX, \bfY, \bfZ
        \in \reals^{d \times m}$
        \\
        $\lesssim$                                         & $\leq$ up to
        positive universal constants
        \\
        $\gtrsim$                                          & $\geq$ up to
        positive universal constants
        \\
        $a_{n} = \bigoh{1}$                                & $(\exists C >
        0)(\exists N \in \nats)(\forall n \geq N)(\absa{a_{n}} < C)$
        \\
        $a_{n} = \bigoh{b_{n}}$                            &
        $\frac{a_{n}}{b_{n}} = \bigoh{1}$
        \\
        $a_{n} = \litoh{1}$                                & $(\forall C >
        0)(\exists N \in \nats)(\forall n \geq N)(\absa{a_{n}} < C)$
        \\
        $a_{n} = \litoh{b_{n}}$                            &
        $\frac{a_{n}}{b_{n}} = \litoh{1}$
        \\
        $X_{n} = \litohp{1}$                               & $(\forall
            \varepsilon > 0)(\Prb{\absa{X_{n}} \geq \varepsilon}
            \xrightarrow[]{n \rightarrow \infty} 0)$
            \\
        $X_{n} = \bigohp{1}$                               & $(\forall
            \varepsilon > 0)(\exists C > 0)(\exists N \in \nats)(\forall n \geq
            N)(\Prb{\absa{X_{n}} \geq C} \leq \varepsilon)$
            \\
        \multicolumn{2}{c}{}
        \\
        \multicolumn{2}{c}{\underline{\textbf{Functions and sets}}}
        \\
        \multicolumn{2}{c}{}
        \\
        $[n]$                                              & $\theseta{1,
        \ldots, n}$, for $n \in \nats$
        \\
        Indicator function $\Indb{\Omega}{\bfx}$           & Takes value $1$
        when $x \in \Omega$, and $0$ otherwise
        \\
        $\Pi_{K} : \reals^{d} \to K$                       &
        $\ell_{2}$-projection of any $\bfx \in \reals^{d}$ onto $K$
        \\
        $f : \Omega \to \reals$ is increasing              & If $\forall u, v
        \in \Omega$ such that $u \leq v \implies f(u) \leq f(v)$
        \\
        $f : \Omega \to \reals$ is decreasing              & If $\forall u, v
        \in \Omega$ such that $u \leq v \implies f(u) \geq f(v)$
        \\
        $\Phi : \reals \to [0, 1]$                         & Cumulative density
        function of $\distNorm(0, 1)$
        \\
        $\phi : \reals \to \reals$                         & Probability density
        function of $\distNorm(0, 1)$
        \\
        $\mclS^{\uparrow}$                                 &
        $\thesetb{\boldsymbol{\mu} \defined (\mu_1, \ldots,
        \mu_{n})^{\top}}{\mu_1 \leq \ldots \leq \mu_{n}}$
        \\
        $\mclS^{\uparrow}_{+}$                             &
        $\thesetb{\boldsymbol{\mu} \in \mclS^{\uparrow}}{\mu_{1} \geq 0}$
        \\
        $\mclS^{\uparrow}(V^{*})$                          &
        $\thesetb{\boldsymbol{\mu} \in \mclS^{\uparrow}}{V(\boldsymbol{\mu})
        \leq V^{*}}$                       \\
        $V(\boldsymbol{\mu})$                              & $\mu_{n} - \mu_{1}$
        for $\boldsymbol{\mu} \in \mclS^{\uparrow}$
        \\
        $\mclS^{\uparrow}_{k^{*}}$                         &
        $\thesetb{\boldsymbol{\mu} \in \mclS^{\uparrow}}{k(\boldsymbol{\mu})
        \leq k^{*}}$                       \\
        $\mclS^{\uparrow}(V^{*}, \eta, \univupperboundmu)$ &
        $\thesetb{\boldsymbol{\mu} \in
        \mclS^{\uparrow}(V^{*})}{\frac{1}{n}\sum_{i = 1}^{n} \mu_{i}^{2} \leq
        \univupperboundmu, \mu_{1} > \eta > 0}$
        \\
        \bottomrule
    \end{tabular}
    \label{tab:notation}
\end{table}

\subsection{Useful miscellaneous results}\label{subsec:useful-misc-results}

Here we prove some useful standard results that are used in several of the
remaining proofs. For reader convenience, we provide short proofs to ensure that
our work is self-contained.

We start with some elementary inequalities, which will be used repeatedly.
First, in \Cref{nlem:square-diff-lower-bound} we introduce a differencing
inequality we use repeatedly to construct lower bounds.

\bnlem[Difference of squares lower bound]\label{nlem:square-diff-lower-bound}
For each $a, b, l \in \reals$, such that $b, l \geq 0$ and $a - b \geq l$, the
following holds:
\begin{equation}\label{neqn:square-diff-lower-bound}
    a^{2} - b^{2} \geq l a\geq l^2
\end{equation}
\enlem

\bprfof{\Cref{nlem:square-diff-lower-bound}}\label{prf:square-diff-lower-bound}
We proceed as follows. First note that since $b, l \geq 0$ by assumption, we
have that $a - b \geq l \iff a \geq b + l \geq 0$. Now observe:
\begin{align}
    a^{2} - b^{2}
     & = (a + b)(a - b)
    \nonumber                                        \\
     & \geq l (a + b)
    \tag{since $a - b \geq l \geq 0$ by assumption.} \\
     & \geq l a
    \tag{since $b \geq 0$ by assumption.}            \\
     & \geq l^2 \tag{since $a \geq l$}
\end{align}
as required.
\eprfof

\bnlem[Lower bound via difference of
    squares]\label{nlem:lower-bound-via-square-diff} For each $a, b, C, K \in
    \reals$, such that $b \geq 0, a^{2} - b^{2} \geq C > 0, a \in [0, K]$, the
    following holds:
\begin{equation}\label{neqn:lower-bound-via-square-diff-01}
    a - b \geq \frac{C}{2 K}
\end{equation}
\enlem

\bprfof{\Cref{nlem:lower-bound-via-square-diff}}\label{prf:lower-bound-via-square-diff}
We proceed as follows. First note that since $a^{2} - b^{2} \geq C > 0$ by
assumption, we have that $a > 0$, and hence $a > b, K > 0$ since both $a, b$ are
non-negative. Now observe:
\begin{align}
    a^{2} - b^{2}
     & \geq C
    \tag{by assumption.}                             \\
    \implies a - b
     & \geq \frac{C}{a + b}
    \tag{since $a > 0, b \geq 0 \implies a +b > 0$.} \\
     & \geq \frac{C}{2 a}
    \tag{since $a \geq b$.}                          \\
     & \geq \frac{C}{2 K}
    \tag{since $a \leq K$.}
\end{align}
as required.
\eprfof

\bnlem[Maximum difference square
    inequality]\label{nlem:square-max-diff-inequality} For each $a, b, c \in
    \reals$ such that $b \leq c$ the following inequality holds:
\begin{equation}\label{neqn:square-max-diff-inequality}
    (\parens{a \vee b} - c)^{2} \leq (a - c)^{2}
\end{equation}
\enlem

\bprfof{\Cref{nlem:square-max-diff-inequality}}\label{prf:square-max-diff-inequality}
Under the assumption that $a, b, c \in \reals$ such that $b \leq c$, let $d := a
\vee b$. We then observe:

\begin{align}
    (d - c)^{2}
     & \leq (a - c)^{2}
    \nonumber                        \\
    \iff d^{2} - a^{2}
     & \leq 2dc - 2ac
    \tag{expanding and simplifying.} \\
    \iff (d + a)(d - a)
     & \leq 2c(d - a)
    \label{neqn:square-max-diff-inequality-01}
\end{align}

\nit So we need to equivalently prove that
\Cref{neqn:square-max-diff-inequality-01}. To that end we only need to consider
2 cases. Namely $a \geq b$, and $a < b$. Note that in the first case $a \geq b
\implies d \defined a \vee b = a$. In this case, both $\lhs / \rhs$ of
\Cref{neqn:square-max-diff-inequality-01} are 0, and the statement holds. Next
consider the case $a < b$. Here we have $a < b \implies d \defined a \vee b = b
> a$. We then observe the following:
\begin{align}
    a + d
     & = a + b
    \tag{since $d = b$.}               \\
     & \leq 2b
    \tag{since $a < b$ by assumption.} \\
     & \leq 2c
    \tag{since $b \leq c$ by assumption.}
\end{align}

That is, we have that $a + d \leq 2c$. Substituting back to
\Cref{neqn:square-max-diff-inequality-01} we have that $(d + a)(d - a) \leq 2c(d
- a)$, which is what we wanted to show. Which completes the proof
\Cref{neqn:square-max-diff-inequality}, as required.
\eprfof

\bnlem[Square sum inequality]\label{nlem:square-sum-inequality} For each $a, b
    \in \reals$ the following holds:
\begin{equation}\label{neqn:square-sum-inequality}
    (a + b)^{2} \leq 2 (a^{2} + b^{2})
\end{equation}
\enlem

\bprfof{\Cref{nlem:square-sum-inequality}}\label{prf:square-sum-inequality} We
proceed as follows:
\begin{align}
    (a + b)^{2}
     & = a^{2} + 2ab + b^{2}                                              \\
     & \leq a^{2} + b^{2} + 2\absa{ab}
    \tag{since $x \leq \absa{x}$ for each $x \in \reals$}                 \\
     & \leq a^{2} + b^{2} + 2 (\absa{a}^{2} + \absa{b}^{2})
    \tag{by AM-GM we have $2 \absa{ab} \leq \absa{a}^{2} + \absa{b}^{2}$} \\
     & = 2(a^{2} + b^{2})
\end{align}
as required.
\eprfof

\nit As a result of \Cref{nlem:square-sum-inequality} we obtain
\cref{ncor:square-sum-inequality-rv}.

\bncor\label{ncor:square-sum-inequality-rv} For random variables $X_{1}, X_{2}$
the following holds:
\begin{equation}\label{neqn:square-sum-inequality-rv}
    \Var{X_{1} - X_{2}} \leq 2 (\Var{X_{1}} + \Var{X_{1}})
\end{equation}
\encor

\bprfof{\Cref{ncor:square-sum-inequality-rv}} First let the centered versions of
the random variables be denoted as
\begin{align}
    \tilde{X}_{i} \defined X_{i} - \Exp{X_{i}}, \text{\, for each $i \in [2]$.}
\end{align}

It then follows that:
\begin{align}
    \Var{X_{1} - X_{2}}
     & = \Var{X_{1} - X_{2} + \Exp{X_{2}} - \Exp{X_{1}}}              \\
     & = \Var{\tilde{X}_{1} - \tilde{X}_{2}}                          \\
     & = \Exp{(\tilde{X}_{1} - \tilde{X}_{2})^{2}}
    \tag{since $\tilde{X}_{1}, \tilde{X}_{2}$ are both centered.}     \\
     & = \Exp{2 (\tilde{X}_{1}^{2} + \tilde{X}_{2}^{2})}
    \tag{using \Cref{nlem:square-sum-inequality}}                     \\
     & = 2 \parens{\Exp{\tilde{X}_{1}^{2}} + \Exp{\tilde{X}_{2}^{2}}}
    \tag{linearity of expectation.}                                   \\
     & = 2 (\Var{X_{1}} + \Var{X_{2}})
    \tag{since $\tilde{X}_{1}, \tilde{X}_{2}$ are both centered.}
\end{align}
as required.
\eprfof

\nit The following is a standard result from real analysis, which we use
repeatedly.

\bnlem[$B$-Lipschitz characterization via bounded
    derivative]\label{nlem:lipschitz-bounded-derivative} Let $f : I \to \reals$
    be continuous and once differentiable, where $I \subseteq \reals$ is an
    interval (possibly unbounded).
\begin{align}\label{neqn:lipschitz-bounded-derivative}
    \text{$f$ is $B$-Lipschitz, with $B > 0$}
     & \iff (\exists B > 0)(\forall x \in \reals) : (\absa{f^{\prime}(x)} \leq B)
\end{align}
\enlem

\bprfof{\Cref{nlem:lipschitz-bounded-derivative}}\label{prf:lipschitz-bounded-derivative}
We prove both directions. In both parts we assume that $f : I \to \reals$ be
continuous and once differentiable, where $I \subseteq \reals$ is an interval
(possibly unbounded). \\
($\implies$). Suppose that $f$ is $B$-Lipschitz, with $B > 0$. We then have
that, for some fixed (but arbitrary) $c \in I$:
\begin{align*}
    \absa{f(x) - f(c)}
     & \leq B \absa{x - c}
    \tag{by definition of $B$-Lipschitz property.} \\
    \implies \absa{\frac{f(x) - f(c)}{x - c}}
     & \leq B
    \tag{taking limits as $x \to c$.}              \\
    \implies \absa{f^{\prime}(c)}
     & \leq B
\end{align*}
Since $c \in I$ is arbitrary, indeed $\absa{f^{\prime}(x)} \leq B$, for each $x
    \in I$, as required. \\
\newline
($\impliedby$). Suppose that $\absa{f^{\prime}(x)} \leq B$, with $B > 0$.
Further let $x, y \in I$, such that $x < y$. Since $f$ is differentiable on $I$,
we have:
\begin{align*}
    \absa{f(x) - f(y)}
     & \leq \absa{f^{\prime}(c)} \absa{x - y}
    \tag{by the mean value theorem, for some $c \in (x, y)$.} \\
     & \leq B \absa{x - y}
    \tag{by assumption.}
\end{align*}
Which implies that $f$ is $B$-Lipschitz, as required.
\eprfof

\bnlem[Standard normal upper bound]\label{nlem:normal-upper-bound} Let $\phi(x),
    \Phi(x)$ respectively denote the probability density function, and
    cumulative density function of a standard normal variable. Then the
    following inequality holds:
\begin{equation}
    \frac{x \phi(x)}{2 \Phi(x) - 1} \leq \frac{1}{2}, \text{ for each $x \geq 0$}
\end{equation}
With equality if and only if $x = 0$.
\enlem

\bprfof{\Cref{nlem:normal-upper-bound}}\label{prf:normal-upper-bound} We first
note that at $x = 0$, that $\frac{x \phi(x)}{2 \Phi(x) - 1}$ is an indeterminate
form of type $\frac{0}{0}$. As such we have:
\begin{align}
    \lim_{x \to 0} \frac{x \phi(x)}{2 \Phi(x) - 1}
     & = \frac{\lim_{x \to 0} \prt{x} x \phi(x)}{\lim_{x \to 0} \prt{x} 2 \Phi(x) - 1}
    \tag{using L'Hospital's rule.}                                                     \\
     & = \frac{\lim_{x \to 0} \phi(x) + x \phi^{\prime}(x)}{\lim_{x \to 0} 2 \phi(x)}
    \nonumber                                                                          \\
     & = \frac{\lim_{x \to 0} \phi(x)}{\lim_{x \to 0} 2 \phi(x)}
    \nonumber                                                                          \\
     & = \frac{\phi(0)}{2 \phi(0)}
    \tag{by continuity of $\phi(x)$ at $x = 0$.}                                       \\
     & = \frac{1}{2}
    \label{neqn:normal-upper-bound-val-zero}
\end{align}
With our given function now defined to be $\frac{1}{2}$ at $x = 0$, we now
proceed to prove our given inequality. Observe that we can equivalently
reformulate it as:
\begin{align}
    \Phi(x) - x\phi(x) - \frac{1}{2} \geq 0
\end{align}
Setting $h(x) \defined \Phi(x) - x\phi(x) - \frac{1}{2}$, we observe that $h(0)
    = \Phi(0) - \frac{1}{2} = 0$. We need to show that $h(x) \geq 0$, for each
    $x \geq 0$, which will imply the result. We will show that $h(x)$ is
    increasing, \ie, or equivalently that $h^{\prime}(x) \geq 0$, for each $x
    \geq 0$. We then have that:
\begin{align}
    h^{\prime}(x)
     & = \phi(x) - (\phi(x) + x \phi^{\prime}(x))
    \nonumber                                                                    \\
     & = - x \phi^{\prime}(x)
    \nonumber                                                                    \\
     & = - x \parens{-x \frac{1}{\sqrt{2 \pi}} e^{-\frac{1}{2}x^{2}}}
    \tag{using $\phi(x) \defined \frac{1}{\sqrt{2 \pi}} e^{-\frac{1}{2}x^{2}}$.} \\
     & = x^{2} \phi(x)
    \nonumber                                                                    \\
     & \geq 0
    \label{neqn:nlem:normal-upper-bound-h-prime-eq}
\end{align}
Note that the inequality in \Cref{neqn:nlem:normal-upper-bound-h-prime-eq} is
strict when $x > 0$ and equality holds if and only if $x = 0$. This means the
function is strictly increasing and bounded away from 0 when for each $x > 0$,
and equal to 0 only when $x = 0$, as required.
\eprfof

\subsection{The Folded Normal
    Distribution}\label{subsec:math-folded-normal-background}

For convenience, we begin by quickly recalling the definition of the Folded
Normal distribution.

\ndefnfoldednormal*
\bnrmk
We note that \Cref{neqn:folded-normal-exp-1} can be equivalently written as
follows:
\begin{equation}
    \sigma \sqrt{2/\pi} \exp(- \mu^2/(2 \sigma^2))
    + \mu(1 - 2\Phi(-\mu/\sigma))
    \label{neqn:folded-normal-exp-2}
\end{equation}
Note that this equivalence follows from the symmetry of the standard normal CDF,
\ie, $\Phi(x) = 1 - \Phi(-x)$ for each $x \in \reals$. For our purposes we
typically use the form of \Cref{neqn:folded-normal-exp-1}.
\enrmk

\subsection{Properties of the folded normal mean:
    \texorpdfstring{$\boldsymbol{f(\mu,
    \sigma)}$}{fmusigma}}\label{subsec:prop-folded-normal-mean}

Let's start setting up some notation. First we note as previously $\absresp_{i}
    \defined \absa{\resp_{i}} = \absa{\mu_{i} + \varepsilon_{i}}$. Where we then
    have $\absresp_{i} \sim \distFoldNorm\parens{\mu_{i}, \sigma^{2}}$. Now
    denote $f(\mu_{i}, \sigma) \defined \Exp{\absresp_{i}}$, for each $i \in
    [n]$. Moreover the $\absresp_{i}$ random variables are all mutually
    independent, but not identically distributed (since their mean's, \ie,
    $f(\mu_{i}, \sigma)$ differ for each $i \in [n]$). Since we run PAVA on
    $\parens{T_{1}, \ldots, T_{n}}$ we have the resulting estimators
    $(\widehat{T}_{1}, \ldots, \widehat{T}_{n})$. We will also denote the
    population level error terms for this transformed (mean centered) response
    as $\delta_{i} \defined \absresp_{i} - f(\mu_{i}, \sigma)$. We note that the
    $\parens{\delta_{1}, \ldots, \delta_{n}}$ are all mutually independent, but
    not identically distributed.

\bnlem[Properties of the Folded Normal mean]\label{nlem:prop-folded-normal-mean}
Suppose $\resp \sim \distNorm(\mu, \sigma^{2})$. Let $T \eqas \absa{\resp}$,
then $\absresp \sim \distFoldNorm\parens{\mu, \sigma^{2}}$ per
\textnormal{\Cref{ndefn:folded-normal}}. We denote the mean of the Folded Normal
distribution by $f(\mu, \sigma) \defined \Exp{T}$. Given this setup, and fixing
$\sigma > 0$, we note the following important properties of $f(\mu, \sigma)$:
\begin{align}
    \text{$f(\mu, \sigma) \geq 0$ for each $\mu \in \reals$} \label{neqn:fold-norm-mean-nonneg}                        \\
    \text{$f(\mu, \sigma) \geq \mu$ for each $\mu \in \reals$} \label{neqn:fold-norm-mean-ge-norm}                     \\
    \text{$f(\mu, \sigma)$ is strictly increasing in $\mu \in \reals_{>0}$} \label{neqn:fold-norm-strictly-inc}        \\
    \text{$\D{f(\mu, \sigma)}{\mu} \in (0, 1)$ is for each $\mu \in \reals_{>0}$}  \label{neqn:fold-norm-mean-bounded} \\
    \text{$f(\mu, \sigma)$ is $1$-Lipschitz for each $\mu \in \reals_{>0}$}  \label{neqn:fold-norm-1-lipschitz}        \\
    \text{$f(\mu, \sigma)^{2} \leq \mu^{2} + \sigma^{2}$  for each $\mu \in \reals_{\geq 0}$}  \label{neqn:fold-norm-mean-le-sec-mom-norm}
\end{align}
Additionally for $\mu_{1} \leq \ldots \leq \mu_{n}$ we have that the
relationship holds for $V(f, \boldsymbol{\mu}, \sigma)$, \ie, the total
variation of the mean of the Folded Normal distribution:
\begin{align}\label{neqn:fold-norm-mean-total-var}
    V(f, \boldsymbol{\mu}, \sigma)
    \defined \sum_{i = 1}^{n - 1} \absa{f(\mu_{i + 1}, \sigma) - f(\mu_{i}, \sigma)}
     & \leq \mu_{n} - \mu_{1}
\end{align}
\enlem

\bprfof{\Cref{nlem:prop-folded-normal-mean}}\label{prf:prop-folded-normal-mean}
We prove each property
(\Cref{neqn:fold-norm-mean-nonneg,neqn:fold-norm-mean-ge-norm,neqn:fold-norm-strictly-inc,neqn:fold-norm-mean-bounded,neqn:fold-norm-1-lipschitz,neqn:fold-norm-mean-le-sec-mom-norm,neqn:fold-norm-mean-total-var})
in turn. As per the assumption $\sigma > 0$ is fixed and that $\resp \sim
\distNorm\parens{\mu, \sigma^{2}}$ for $\mu \in \reals$. \qedbsquare \\
\newline
(\textit{Proof of \Cref{neqn:fold-norm-mean-nonneg}}.) We have that $T \defined
    \absa{\resp} \geq 0 \text{ \as } \implies f(\mu, \sigma) \defined \Exp{T} =
    \Exp{\absa{\resp}} \geq 0$ by the monotonicity of expectation, as required.
    \qedbsquare \\
\newline
(\textit{Proof of \Cref{neqn:fold-norm-mean-ge-norm}}.) We have that $\resp \leq
    \absa{\resp} \text{\as } \implies \mu \defined \Exp{\resp} \leq
    \Exp{\absa{\resp}} = \Exp{T} \defines f(\mu, \sigma)$ again by the
    monotonicity of expectation, as required. \qedbsquare \\
\newline
(\textit{Proof of \Cref{neqn:fold-norm-strictly-inc}}.) For any $\mu > 0$ we
have that:
\begin{align*}
    f(\mu, \sigma)
     & \defined \sigma \sqrt{\frac{2}{\pi}} \exp\parens{- \frac{\mu^{2}}{2 \sigma^{2}}} -
    \mu\parens{1 - 2\Phi\parens{\frac{\mu}{\sigma}}}
    \tag{per \Cref{neqn:folded-normal-exp-1}}                                                     \\
    \implies
    \D{f(\mu, \sigma)}{\mu}
     & = -\frac{\mu}{\sigma} \sqrt{\frac{2}{\pi}} \exp\parens{- \frac{\mu^{2}}{2 \sigma^{2}}} - 1
    + 2\Phi\parens{\frac{\mu}{\sigma}}
    + \frac{2 \mu}{\sigma} \phi\parens{\frac{\mu}{\sigma}}                                        \\                               \\
     & = 2\Phi\parens{\frac{\mu}{\sigma}} - 1
    \tag{since $\frac{\mu}{\sigma} \sqrt{\frac{2}{\pi}} \exp\parens{- \frac{\mu^{2}}{2 \sigma^{2}}}
    = \frac{2 \mu}{\sigma} \phi\parens{\frac{\mu}{\sigma}}$}                                      \\
     & > 0 \tag{since $\mu, \sigma > 0$ and $\Phi\parens{\frac{\mu}{\sigma}} > \frac{1}{2}$}
\end{align*}
as required. \qedbsquare \\
\newline
(\textit{Proof of \Cref{neqn:fold-norm-mean-bounded}}.) By the previous proof,
we note that $\D{f(\mu, \sigma)}{\mu} > 0$. Also using the previous proof and
noting that $\Phi(x) > \frac{1}{2}$ for each $x > 0$, it follows that $\D{f(\mu,
\sigma)}{\mu} = - 1 + 2\Phi\parens{\frac{\mu}{\sigma}} < 1$. Combining both
parts we have that $\D{f(\mu, \sigma)}{\mu} \in (0, 1)$, as required.
\qedbsquare \\
\newline
(\textit{Proof of \Cref{neqn:fold-norm-1-lipschitz}}.) By the previous proof, we
note that $\D{f(\mu, \sigma)}{\mu} \in (0, 1) \implies \absa{\D{f(\mu,
\sigma)}{\mu}} \leq 1$ for each $\mu > 0$. It follows by the mean value theorem,
that $f(\mu, \sigma)$ is $1$-Lipschitz as required. \qedbsquare \\
\newline
(\textit{Proof of \Cref{neqn:fold-norm-mean-le-sec-mom-norm}}.) Observe that
from \Cref{neqn:folded-normal-var} we have that $g(\mu, \sigma) \defined
\Var{\absresp} = \mu^{2} + \sigma^{2} - f(\mu, \sigma)^{2}$. Since
$\Var{\absresp} \geq 0$, it follows that $f(\mu, \sigma)^{2} \leq \mu^{2} +
\sigma^{2}$  for each $\mu \in \reals$, as required. \qedbsquare \\
\newline
(\textit{Proof of \Cref{neqn:fold-norm-mean-total-var}}.) Let $i \in [n]$ be
arbitrary. Now note that by the \Cref{neqn:fold-norm-1-lipschitz} property it
follows that , we then have that:
\begin{align}
    V(f, \boldsymbol{\mu}, \sigma)
     & \defined \sum_{i = 1}^{n - 1} \absa{f(\mu_{i + 1}, \sigma) - f(\mu_{i}, \sigma)}
    \tag{by definition}                                                                                   \\
     & = \sum_{i = 1}^{n - 1} f(\mu_{i + 1}, \sigma) - f(\mu_{i}, \sigma)
    \tag{using \Cref{neqn:fold-norm-strictly-inc} and monotonicity of $\mu_{1} \leq \ldots \leq \mu_{n}$} \\
     & = f(\mu_{n}, \sigma) - f(\mu_{1}, \sigma)
    \tag{by telescoping sum}                                                                              \\
     & \leq \absa{\mu_{n} - \mu_{1}}
    \tag{using \Cref{neqn:fold-norm-1-lipschitz}}                                                         \\
     & = \mu_{n} - \mu_{1}
    \tag{by monotonicity of $\mu_{1} \leq \ldots \leq \mu_{n}$}
\end{align}
as required. \qedbsquare \\
\newline
\nit Thus all properties specified in
\Cref{neqn:fold-norm-mean-nonneg,neqn:fold-norm-mean-ge-norm,neqn:fold-norm-strictly-inc,neqn:fold-norm-mean-bounded,neqn:fold-norm-1-lipschitz,neqn:fold-norm-mean-le-sec-mom-norm,neqn:fold-norm-mean-total-var}
are now proved.
\eprfof

\subsection{Properties of the folded normal variance:
    \texorpdfstring{$\boldsymbol{g(\mu,
    \sigma)}$}{gmusigma}}\label{subsec:prop-folded-normal-var}

\bnlem[Properties of the Folded Normal
    variance]\label{nlem:prop-folded-normal-variance} Let $\absresp \sim
    \distFoldNorm\parens{\mu, \sigma^{2}}$ per \Cref{ndefn:folded-normal}, and
    let $g(\mu, \sigma) \defined \Var{T}$. Given this setup, and fixing $\sigma
    > 0$, we note the following properties of $g(\mu, \sigma)$:
\begin{align}
    g(\mu, \sigma)
     & \leq \sigma^{2} \text{, for each $\mu \in \reals$}
    \label{neqn:fold-norm-var-ge-sigma-sq}                                      \\
    g(\mu, \sigma)
     & \geq g(0, \sigma) \text{, for each $\mu \in \reals_{>0}$}
    \label{neqn:fold-norm-var-inc-in-mu}                                        \\
    \Var{T^{2}}
     & = 4 \mu^{2} \sigma^{2} + 2 \sigma^{4} \text{, for each $\mu \in \reals$}
    \label{neqn:fold-norm-sq-var}
\end{align}
\enlem

\bprfof{\Cref{nlem:prop-folded-normal-variance}}\label{prf:prop-folded-normal-variance}
We prove each properties specified in
\Cref{neqn:fold-norm-var-ge-sigma-sq,neqn:fold-norm-var-inc-in-mu,neqn:fold-norm-sq-var}
in turn. \\
\newline
(\textit{Proof of \Cref{neqn:fold-norm-var-ge-sigma-sq}}.) We have for each $\mu
    \geq 0$
\begin{align*}
    g(\mu, \sigma)
     & \defined \mu^{2} + \sigma^{2} - f(\mu, \sigma)^{2}
    \tag{per \Cref{neqn:folded-normal-var}}               \\
     & \leq \sigma^{2}
    \tag{since $f(\mu, \sigma)^{2} \geq \mu^{2}$ using \Cref{neqn:fold-norm-mean-ge-norm}}
\end{align*} \\
as required. \qedbsquare \\
\newline
(\textit{Proof of \Cref{neqn:fold-norm-var-inc-in-mu}}.) First note that $g(0,
    \sigma) = \sigma^{2} - f(0, \sigma)^{2} = \sigma^{2} -
    \parens{\frac{2}{\pi}}\sigma^{2}$. It then follows that:
\begin{align}
    g(\mu, \sigma)
     & \geq g(0, \sigma)  \nonumber                                       \\
    \iff \mu^{2} + \sigma^{2} - f(\mu, \sigma)^{2}
     & \geq \sigma^{2} - \parens{\frac{2}{\pi}}\sigma^{2} \nonumber       \\
    \iff \mu^{2} + \parens{\frac{2}{\pi}}\sigma^{2}
     & \geq f(\mu, \sigma)^{2} \label{neqn:fold-norm-var-inc-in-mu-equiv}
\end{align} \\
We will then prove the equivalent statement
\Cref{neqn:fold-norm-var-inc-in-mu-equiv}. Since $\mu, \sigma > 0$ in our case,
let $\nu \defined \frac{\mu}{\sigma} > 0$ in what follows. Then dividing both
sides of \Cref{neqn:fold-norm-var-inc-in-mu-equiv} by $\nu$ we obtain:
\begin{align*}
    \sqrt{\nu^{2}+\frac{2}{\pi}} & \geq \nu(2 \Phi(\nu) - 1)+\sqrt{\frac{2}{\pi}} e^{-\frac{\nu^{2}}{2}}
\end{align*}

\nit Let us then define
\begin{align}
    g(\nu)
     & \defined \sqrt{\nu^{2}+\frac{2}{\pi}} - \nu(2 \Phi(\nu)-1)
    - \sqrt{\frac{2}{\pi}} e^{-\frac{\nu^{2}}{2}}
\end{align}

\nit Taking the derivative of $g(\nu)$ with respect to $\nu$ we obtain:
\begin{align}
    g^{\prime}(\nu)
     & = \frac{\nu}{\sqrt{\nu^{2}+\frac{2}{\pi}}}-2 \Phi(\nu) + 1
    - \frac{2 \nu}{\sqrt{2 \pi}} e^{-\frac{\nu^{2}}{2}}
    + \nu\sqrt{\frac{2}{\pi}} e^{-\frac{\nu^{2}}{2}}              \\
     & = \frac{\nu}{\sqrt{\nu^{2}+\frac{2}{\pi}}}-2 \Phi(\nu) + 1
\end{align}

\nit As such all global and local extrema are obtained by setting
$g^{\prime}(\nu) = 0$, that is:
\begin{align}
    g^{\prime}(\nu)
     & = 0               \\
    \iff \frac{\nu}{\sqrt{\nu^{2}+\frac{2}{\pi}}}
     & = 2 \Phi(\nu) - 1
\end{align}
Now, $\nu=0$ is a clear solution, at which our function is exactly equal to $0$.
Also, we need to look at $\nu=\infty$, where we also have an identity. So we
need to take care of other possible roots to the equation
$\frac{\nu}{\sqrt{\nu^{2} + \frac{2}{\pi}}}= 2 \Phi(\nu) - 1$. Now observe that
since when $\nu \geq 0$ the function $\Phi(\nu)$ is concave and therefore $2
\Phi(\nu)-1=2(\Phi(\nu)-\Phi(0)) \leq 2 \nu \phi(0)=\sqrt{\frac{2}{\pi}} \nu$.
Thus for any non-zero solution $\bar{\nu}$ to the equation
$\frac{\nu}{\sqrt{\nu^{2}+\frac{2}{\pi}}}=2 \Phi(\nu)-1$, we must have
$\frac{\bar{\nu}}{\sqrt{\bar{\nu}^{2}+\frac{2}{\pi}}} \leq \bar{\nu}
\sqrt{\frac{2}{\pi}}$.

\nit This implies that $\bar{\nu}^{2} \geq \frac{\pi}{2}-\frac{2}{\pi}$. Now,
going back to the original function we need to show
\begin{equation*}
    \sqrt{\bar{\nu}^{2}+\frac{2}{\pi}}
    \geq \bar{\nu}(2 \Phi(\bar{\nu})-1)+\sqrt{\frac{2}{\pi}} e^{-\bar{\nu}^{2} / 2}=\frac{\bar{\nu}^{2}}{\sqrt{\bar{\nu}^{2}+\frac{2}{\pi}}}+\sqrt{\frac{2}{\pi}} e^{-\bar{\nu}^{2} / 2} .
\end{equation*}
The latter is equivalent to $\bar{\nu}^{2}+\frac{2}{\pi} \geq
    \bar{\nu}^{2}+\sqrt{\frac{2}{\pi}} e^{-\bar{\nu}^{2} / 2}
    \sqrt{\bar{\nu}^{2}+\frac{2}{\pi}}$ which is equivalent to $\frac{2}{\pi}
    e^{\bar{\nu}^{2}} \geq \bar{\nu}^{2}+\frac{2}{\pi}$. since the function $\nu
    \mapsto \frac{2}{\pi} e^{\nu^{2}}-\nu^{2}$ is increasing for positive $\nu$
    it suffices to check that $\frac{2}{\pi} e^{\bar{\nu}^{2}} \geq
    \bar{\nu}^{2}+\frac{2}{\pi}$ for $\bar{\nu}=\frac{\pi}{2}-\frac{2}{\pi}$
    (since as we know from before $\bar{\nu}$ is at least that value). This is
    true, and completes the proof, as required. \qedbsquare \\
\newline
(\textit{Proof of \Cref{neqn:fold-norm-sq-var}}.) By direct calculation we have:
\begin{align*}
    \Var{\absresp^{2}}
     & = \Exp{\absresp^{4}} - \parens{\Exp{\absresp^{2}}}^{2}                    \\
     & = \Exp{\resp^{4}} - \parens{\Exp{\resp^{2}}}^{2}
    \tag{since $\absresp \eqas \absa{\resp}$}                                    \\
     & = (\mu^{4} + 6 \mu^{2} \sigma^{2} + 3 \sigma^{4})
    - (\mu^{2} + \sigma^2)^{2}
    \tag{$2^{\text{nd}}, 4^{\text{th}}$ moments of $\distNorm(\mu, \sigma^{2})$} \\
     & = 4 \mu^{2} \sigma^{2} + 2 \sigma^{4}
\end{align*}
as required. \qedbsquare \\
\newline
\nit Thus all properties specified in
\Cref{neqn:fold-norm-var-ge-sigma-sq,neqn:fold-norm-var-inc-in-mu,neqn:fold-norm-sq-var}
are now proved.
\eprfof

\subsection{Properties of the inverse folded normal mean:
    \texorpdfstring{$\boldsymbol{f^{-1}(\mu,
    \sigma)}$}{finvmusigma}}\label{subsec:prop-folded-normal-inv-mean}

\bnlem[Properties of the Folded Normal mean
    inverse]\label{nlem:prop-folded-normal-inv-mean} Suppose $\resp \sim
    \distNorm(\mu, \sigma^{2})$. Let $T \eqas \absa{\resp}$, then $\absresp \sim
    \distFoldNorm\parens{\mu, \sigma^{2}}$ per
    \textnormal{\Cref{ndefn:folded-normal}}. We denote the mean of the Folded
    Normal distribution by $f(\mu, \sigma) \defined \Exp{T}$. Given this setup,
    and fixing $\sigma > 0$, we note the following important properties of
    $f^{-1}(u, \sigma)$ (which denotes the inverse with respect to $\mu$
    function of $f(\mu, \sigma)$ when $\sigma$ is held fixed):
\begin{align}
    \text{$f^{-1}(u, \sigma)$ exists,} \label{neqn:fold-norm-inv-mean-exists}                             \\
    \text{$\prt{\sigma}f^{-1}(u, \sigma) = -\frac{\sqrt{2/\pi} \exp(-\mu^2/(2\sigma^2))}{2
    \Phi(\mu/\sigma) - 1}$,} \label{neqn:fold-norm-inv-mean-der-sigma}                                    \\
    \text{$\prt{u}f^{-1}(u, \sigma) = 1/(2 \Phi(\mu/\sigma) - 1)$,} \label{neqn:fold-norm-inv-mean-der-u} \\
    \text{$f^{-1}(u, \sigma)$ is a Lipschitz function for each $u > f(\eta, \sigma) > 0$ for a fixed $\sigma$,}  \label{neqn:fold-norm-inv-mean-der-u-lipschitz}
\end{align}
where in the above $u = f(\mu,\sigma)$ (or in other words $\mu = f^{-1}(u,
    \sigma)$).
\enlem

\bprfof{\Cref{nlem:prop-folded-normal-inv-mean}} We prove each properties
specified in
\Cref{neqn:fold-norm-inv-mean-exists,neqn:fold-norm-inv-mean-der-sigma,neqn:fold-norm-inv-mean-der-u,neqn:fold-norm-inv-mean-der-u-lipschitz}
in turn. \\
\newline
(\textit{Proof of \Cref{neqn:fold-norm-inv-mean-exists}}.) Note that for a fixed
$\sigma > 0$ the function $f(\mu,\sigma)$ is invertible (as it is increasing,
per \Cref{nlem:prop-folded-normal-mean}), as required. \qedbsquare \\
\newline
(\textit{Proof of \Cref{neqn:fold-norm-inv-mean-der-sigma}}.) In order to find
the derivative of $\frac{\partial}{\partial \sigma} f^{-1}(\cdot, \sigma)$, we
can parametrize as follows:
\begin{align}
    u
     & = f(\mu, \sigma) \\
    v
     & = \sigma
\end{align}
We will use the inverse function theorem which says that under certain
conditions $\mu = F(u, v) = F(u, \sigma)$ and $\sigma = G(u,v) = v$, for some
functions $F$ and $G$. Note that for a fixed $\sigma > 0$ the function
$f(\mu,\sigma)$ is invertible (per \Cref{neqn:fold-norm-inv-mean-exists}). Thus
\begin{align}
    \prt{\sigma}f^{-1}(u, \sigma)
    = \D{\mu}{\sigma}
    = \D{F(u,v)}{v}
    = -\frac{\frac{\partial u}{\partial \sigma}}{J} =
    -\frac{\sqrt{2/\pi} \exp(-\mu^2/(2\sigma^2))}{J}
\end{align}
Where $J$ is the Jacobian of the transformation
\begin{align*}
    J & = \bigg| \begin{array}{cc}\frac{\partial u}{\partial \mu}
                      & \frac{\partial u}{\partial \sigma} \\ \frac{\partial v}{\partial \mu} &
                     \frac{\partial v}{\partial \sigma}\end{array}\bigg| \\
      & = 2 \Phi(\mu/\sigma) - 1                                                          \\
      & > 0 \tag{since $\mu \geq \eta > 0$}
\end{align*}
It follows that
\begin{align}
    \frac{\partial }{ \partial
        \sigma}f^{-1}(u, \sigma) = -\frac{\sqrt{2/\pi} \exp(-\mu^2/(2\sigma^2))}{2
        \Phi(\mu/\sigma) - 1}
\end{align}
As required. \qedbsquare \\
\newline
(\textit{Proof of \Cref{neqn:fold-norm-inv-mean-der-u}}.) We similarly evaluate
the derivative $\prt{u} f^{-1}(u, \sigma)$ as follows:
\begin{align}
    \prt{u} f^{-1}(u, \sigma)
     & = \prt{u} \mu                      \\
     & = \frac{\D{v}{\sigma}}{J}          \\
     & = \frac{1}{2 \Phi(\mu/\sigma) - 1}
\end{align}
As required. \qedbsquare \\
\newline
(\textit{Proof of \Cref{neqn:fold-norm-inv-mean-der-u-lipschitz}}.) We note that
\Cref{neqn:fold-norm-inv-mean-der-u} implies that
\begin{equation}
    \prt{u} f^{-1}(u, \sigma)
    \leq \frac{1}{2 \Phi(\eta /\sigma) - 1}
\end{equation}
since $\mu \geq \eta > 0$ under our setting. In this case, this holds for each
$u > f(\eta, \sigma) > 0$, for a fixed $\sigma$. Since this derivative is
bounded by this constant, it follows that $f^{-1}(u, \sigma)$ is $\frac{1}{2
\Phi(\eta /\sigma) - 1}$-Lipschitz by applying
\Cref{nlem:lipschitz-bounded-derivative}.
\qedbsquare \\
\newline
\nit Thus all properties specified in
\Cref{neqn:fold-norm-inv-mean-exists,neqn:fold-norm-inv-mean-der-sigma,neqn:fold-norm-inv-mean-der-u,neqn:fold-norm-inv-mean-der-u-lipschitz}
are now proved.
\eprfof

\subsection{Properties of:
    \texorpdfstring{$J(\sigma)$}{jsigma}}\label{subsec:prop-j-sigma}

\bnumdefn[$J(\sigma)$]\label{ndefn:j-sigma} Let $\eta > 0$ be fixed, and $\sigma
    \geq 0$ per \Cref{neqn:adversarial-iso-reg-1,neqn:adversarial-iso-reg-2},
    respectively. We define the function, $J: \reals_{\geq 0} \to \reals$, as:
\begin{equation}\label{neqn:j-sigma}
    J(\sigma)
    \defined
    \begin{cases}
        0
         & \text{if $\sigma = 0$} \\
        \sigma\parens{\frac{1}{2} - \frac{\eta/\sigma \phi(\eta/\sigma)}{2 \Phi(\eta/\sigma) - 1}}
         & \text{otherwise}
    \end{cases}
\end{equation}
\enumdefn

\nit In order to prove the key properties of $J(\sigma)$, we will first need to
prove a useful result in \Cref{nlem:j2-sigma-decreasing-aux}.

\bnlem\label{nlem:j2-sigma-decreasing-aux} We define the function, $M:
    \reals_{\geq 0} \to \reals_{\geq 0}$, as:
\begin{equation}\label{neqn:j2-sigma-decreasing-aux-01}
    M(x) \defined
    \begin{cases}
        \frac{1}{2}
         & \textnormal{if $x = 0$} \\
        \frac{x \phi(x)}{2 \Phi(x) - 1}
         & \textnormal{otherwise}.
    \end{cases}
\end{equation}
Note that $M(0) = \frac{1}{2}$, by \Cref{neqn:normal-upper-bound-val-zero}. Then
$M(x)$ is strictly decreasing for each $x > 0$.
\enlem
\bprfof{\Cref{nlem:j2-sigma-decreasing-aux}} In order to show that $M(x)$ is
decreasing for each $x > 0$, we will show that $M^{\prime}(x) < 0$ for each $x >
0$. To see this, first observe that:
\begin{align}
    M^{\prime}(x)
     & = \frac{(2 \Phi(x)-1)\parens{\phi(x)+x \phi^{\prime}(x)}-2 x \phi^{2}(x)}{(2 \Phi(x)-1)^{2}}
    \nonumber                                                                                       \\
     & = \frac{(2 \Phi(x) - 1)\parens{\phi(x) - x^{2} \phi(x)}-2 x \phi^{2}(x)}{(2 \Phi(x)-1)^{2}}
    \tag{since $\phi^{\prime}(x) + x \phi(x) = 0$}                                                  \\
     & = \frac{\phi(x)}{(2 \Phi(x)-1)^{2}}
    \parens{\parens{2 \Phi(x) - 1}\parens{1 - x^{2}} - 2 x \phi(x)}
    \nonumber                                                                                       \\
     & = \frac{\phi(x)}{(2 \Phi(x)-1)^{2}}
    \parens{2 \parens{\Phi(x) - \frac{1}{2}}\parens{1 - x^{2}} - 2 x \phi(x)}.
    \label{neqn:j2-sigma-decreasing-aux-02}
\end{align}
Now we see that:
\begin{equation}\label{neqn:j2-sigma-decreasing-aux-03}
    2 \parens{\Phi(x)-\frac{1}{2}}
    = 2 \parens{\Phi(x) - \Phi(0)} = 2 \int_{0}^{x} \frac{1}{\sqrt{2 \pi}} e^{-\frac{t^{2}}{2}} \dt
    = \frac{2}{\sqrt{2 \pi}} \int_{0}^{x} e^{-\frac{t^{2}}{2}} \dt
    \leq \frac{2 x}{\sqrt{2 \pi}},
\end{equation}
where the last inequality in \Cref{neqn:j2-sigma-decreasing-aux-03} followed
from the fact that $e^{-\frac{t^{2}}{2}} \leq 1$ for each $t \geq 0$. It then
follows that:
\begin{align}
    \parens{2 \parens{\Phi(x) - \frac{1}{2}}\parens{1 - x^{2}} - 2 x \phi(x)}
     & = \parens{\frac{2}{\sqrt{2 \pi}} \int_{0}^{x} e^{-\frac{t^{2}}{2}} \dt}\parens{1 - x^{2}}
    - \frac{2 x}{\sqrt{2 \pi}} e^{-\frac{x^{2}}{2}}
    \nonumber                                                                                    \\
     & \leq \frac{2 x}{\sqrt{2 \pi}} \parens{1 - x^{2}}
    - \frac{2 x}{\sqrt{2 \pi}} e^{-\frac{x^{2}}{2}}
    \tag{using \Cref{neqn:j2-sigma-decreasing-aux-03}}                                           \\
     & < 0,
    \label{neqn:j2-sigma-decreasing-aux-04}
\end{align}
Where \Cref{neqn:j2-sigma-decreasing-aux-04} followed by observing that since $1
    - x^{2} < 1 - \frac{x^{2}}{2} < e^{-\frac{x^{2}}{2}}$ for each $x > 0$. Now
    since $\frac{\phi(x)}{(2 \Phi(x)-1)^{2}} > 0$ for each $x > 0$, we have by
    applying \Cref{neqn:j2-sigma-decreasing-aux-04} to
    \Cref{neqn:j2-sigma-decreasing-aux-03} that $M^{\prime}(x) < 0$, for each $x
    > 0$, as required.
\eprfof

\bnlem[Properties of $J(\sigma)$]\label{nlem:j-sigma-properties} Let $J(\sigma)$
be defined as per \Cref{neqn:j-sigma}. Then $J(\sigma)$ satisfies the following
properties:
\begin{align}
    \text{$J(\sigma) > 0$ for each $\sigma \in \reals_{> 0}$ and 0 if and only if $\sigma = 0$} \label{neqn:j-sigma-properties-01} \\
    \text{$J(\sigma)$ is continuous for each $\sigma \in \reals_{> 0}$} \label{neqn:j-sigma-properties-02}                         \\
    \text{For any $0 < \sigma_{1} < \sigma_{2}$,
        $\min_{\sigma \in [\sigma_{1}, \sigma_{2}]} J(\sigma) \geq
            \sigma_{1} \parens{\frac{1}{2} - \frac{\eta/\sigma_{2} \phi(\eta/\sigma_{2})}{2 \Phi(\eta/\sigma_{2}) - 1}}
            > 0$} \label{neqn:j-sigma-properties-03}
\end{align}
\enlem

\bprfof{\Cref{nlem:j-sigma-properties}}\label{prf:unif-sigma-lower-bound} We
prove each property
(\Cref{neqn:j-sigma-properties-01,neqn:j-sigma-properties-02,neqn:j-sigma-properties-03})
in turn. Throughout these proofs, we write:
\begin{equation}
    J(\sigma)
    \defined J_{1}(\sigma) J_{2}(\sigma) \text{, where $J_{1}(\sigma) \defined \sigma$, and $J_{2}(\sigma)
            \defined \frac{1}{2} - \frac{\eta/\sigma \phi(\eta/\sigma)}{2 \Phi(\eta/\sigma)
                - 1}$}
\end{equation}
\newline
(\textit{Proof of \Cref{neqn:j-sigma-properties-01}}.) Observe that both
$J_{1}(\sigma), J_{2}(\sigma)$ are zero if and only if $\sigma = 0$. In the case
of $J_{2}(\sigma)$ this follows from \Cref{nlem:normal-upper-bound}. Now for
$\sigma > 0$, $J_{1}(\sigma) \defined \sigma > 0$, by assumption. And the fact
that $J_{2}(\sigma) > 0$, for $\sigma > 0$ again follows directly from
\Cref{nlem:normal-upper-bound}. As such, $J(\sigma) > 0$ for each $\sigma > 0$,
since it is the product of two strictly positive functions over this support, as
required. \qedbsquare \\
\newline
(\textit{Proof of \Cref{neqn:j-sigma-properties-02}}.) $J_{1}(\sigma)$ is
continuous for $\sigma > 0$. Moreover since $\phi(x), \Phi(x)$ for a standard
normal are continuous over their support, $\reals$, it follows that
$J_{2}(\sigma)$ is also continuous for $\sigma > 0$. As such, $J(\sigma)$ is
continuous for each $\sigma > 0$, since it is the product of two continuous
functions, as required. \qedbsquare \\
\newline
(\textit{Proof of \Cref{neqn:j-sigma-properties-03}}.) Note that for any two
fixed $\sigma_{1}, \sigma_{2}$, such that $0 < \sigma_{1} < \sigma_{2}$, the
interval $[\sigma_{1}, \sigma_{2}]$ is compact. From
\Cref{neqn:j-sigma-properties-02}, we know that $J(\sigma)$ is continuous for
$\sigma > 0$, and so it attains it's minimum (and maximum) on this interval.
Moreover, from \Cref{neqn:j-sigma-properties-01}, it follows that $\min_{\sigma
\in [\sigma_{1}, \sigma_{2}]} J(\sigma) > 0$. Now we note that $\sigma \mapsto
J_{1}(\sigma) \defined \sigma$, is increasing in $\sigma$. Moreover for each
$\sigma > 0$, we have that $J_{2}(\sigma) \defined \frac{1}{2} -
M(\frac{\eta}{\sigma})$, where the function $M$ is as defined in
\Cref{neqn:j2-sigma-decreasing-aux-01}. Moreover it follows from
\Cref{nlem:j2-sigma-decreasing-aux} that $J_{2}(\sigma)$ is strictly decreasing
for each $\sigma > 0$. By the non-negativity of $J(\sigma)$ over it's domain, we
have that $J(\sigma) > J_{1}(\sigma_{1}) J_{2}(\sigma_{2})$ for each $\sigma \in
[\sigma_{1}, \sigma_{2}]$. From this we have that $\min_{\sigma \in [\sigma_{1},
\sigma_{2}]} J(\sigma) \geq J_{1}(\sigma_{1}) J_{2}(\sigma_{2}) = \sigma_{1}
\parens{\frac{1}{2} - \frac{\eta/\sigma_{2} \phi(\eta/\sigma_{2})}{2
\Phi(\eta/\sigma_{2}) - 1}} > 0$, as required. \qedbsquare \\
\newline
\nit Thus all properties specified in
\Cref{neqn:j-sigma-properties-01,neqn:j-sigma-properties-02,neqn:j-sigma-properties-03}
are now proved.
\eprfof

\subsection{Properties of:
    \texorpdfstring{$G(\sigma)$}{gsigma}}\label{subsec:prop-g-sigma}

\bnumdefn[$G(\sigma)$]\label{ndefn:g-sigma} Under the setup of \texttt{ASCI}
generating process per \Cref{ndefn:isotonic-adversarial-sign-model}, and per the
\texttt{ASCIFIT} model we define the function, $G: \reals_{\geq 0} \to \reals$,
as:
\begin{equation}\label{neqn:g-sigma}
    G(\sigma)
    \defined
    \sigma^{2} +
    \frac{1}{n} \sum_{i = 1}^{n} (f^{-1}(\widehat{\absresp}_{i} \vee f(\eta, \sigma), \sigma))^2
\end{equation}
\enumdefn

\bnlem[Properties of $G(\sigma)$]\label{nlem:prop-g-sigma} Under the setup of
\texttt{ASCI} generating process per
\Cref{ndefn:isotonic-adversarial-sign-model}, and with $G(\sigma)$ defined as
per \Cref{ndefn:g-sigma}, we note the following important properties of
$G(\sigma)$:
\begin{align}
    \text{$\prt{\sigma} G(\sigma)
    = \frac{4}{n} \sum_{i = 1}^{n} \sigma {\parens{\frac{1}{2} - \frac{ f^{-1}(\widehat{\absresp}_{i}, \sigma)/\sigma\phi(f^{-1}(\widehat{\absresp}_{i}, \sigma)/\sigma)\mathbbm{1}(\widehat{\absresp}_{i} \geq f(\eta, \sigma)) }{2 \Phi(f^{-1}(\widehat{\absresp}_{i}, \sigma)/\sigma) - 1}}}$}
    \label{neqn:prop-g-sigma-01}, \\
    \text{$G(\sigma)$ is increasing for $\sigma \geq 0$, and strictly increasing for $\sigma > 0$.}
    \label{neqn:prop-g-sigma-02}
\end{align}
\enlem

\bprfof{\Cref{nlem:prop-g-sigma}}\label{prf:prop-g-sigma} We prove each property
(\Cref{neqn:prop-g-sigma-01,neqn:prop-g-sigma-02}) in turn. Throughout these
proofs, $J(\sigma)$ is as defined in \Cref{ndefn:j-sigma}, and $G(\sigma)$ is as
defined in \Cref{ndefn:g-sigma}. \\
\newline
(\textit{Proof of \Cref{neqn:prop-g-sigma-01}}.). Using the definition, we have:
\begin{align}
     & \prt{\sigma}G(\sigma)
    \nonumber                                                                                                                                                                                                                                                                                          \\
     & = 2 \sigma - \frac{2}{n} \sum_{i = 1}^{n} \frac{\sqrt{2/\pi} f^{-1}(\widehat{\absresp}_{i}, \sigma) \exp(-f^{-1}(\widehat{\absresp}_{i}, \sigma)^2/(2\sigma^2)) \mathbbm{1}(\widehat{\absresp}_{i} \geq f(\eta, \sigma))}{2 \Phi(f^{-1}(\widehat{\absresp}_{i} , \sigma)/\sigma) - 1} \nonumber
    \tag{using \Cref{neqn:fold-norm-inv-mean-der-sigma}}                                                                                                                                                                                                                                               \\
     & = 2\sigma - \frac{4\sigma}{ n} \sum_{i = 1}^{n} \frac{ f^{-1}(\widehat{\absresp}_{i}, \sigma)/\sigma\phi(f^{-1}(\widehat{\absresp}_{i}, \sigma)/\sigma)\mathbbm{1}(\widehat{\absresp}_{i} \geq f(\eta, \sigma)) }{2 \Phi(f^{-1}(\widehat{\absresp}_{i}, \sigma)/\sigma) - 1}
    \label{neqn:prop-g-sigma-03}                                                                                                                                                                                                                                                                       \\
     & = \frac{4}{n} \sum_{i = 1}^{n} \sigma {\parens{\frac{1}{2} - \frac{ f^{-1}(\widehat{\absresp}_{i}, \sigma)/\sigma\phi(f^{-1}(\widehat{\absresp}_{i}, \sigma)/\sigma)\mathbbm{1}(\widehat{\absresp}_{i} \geq f(\eta, \sigma)) }{2 \Phi(f^{-1}(\widehat{\absresp}_{i}, \sigma)/\sigma) - 1}}}
    \nonumber
\end{align}
as required. \qedbsquare \\
\newline
(\textit{Proof of \Cref{neqn:prop-g-sigma-02}}.). Now per
\Cref{nlem:normal-upper-bound}, we have that $x \mapsto x \phi(x)/(2\Phi(x)-1)
\leq 1/2$ for all $x \geq 0$, and moreover it is decreasing for $x > 0$.
Therefore the derivative, $\prt{\sigma}G(\sigma)$, is bounded from below by $0$.
As such $G(\sigma)$ is increasing in $\sigma$, for $\sigma \geq 0$. In fact,
since $\eta, \sigma > 0$, it follows that $\frac{\eta}{\sigma} > 0$. In turn, we
have that $\prt{\sigma} G(\sigma)$ is bounded from below by $J(\sigma) \defined
\sigma\parens{\frac{1}{2} - \frac{\eta/\sigma \phi(\eta/\sigma)}{2
\Phi(\eta/\sigma) - 1}} > 0$, for each $\sigma > 0$, using
\Cref{nlem:j-sigma-properties}. It follows that $G(\sigma)$ is \textit{strictly}
increasing in $\sigma$, for $\sigma > 0$, as required. \qedbsquare \\
\newline
\nit Thus all properties specified in
\Cref{neqn:prop-g-sigma-01,neqn:prop-g-sigma-02} are now proved.
\eprfof

\newpage
\section{Proofs of \Cref{sec:intro}}\label{sec:intro-and-prior-work-proofs}

\subsection{Mathematical
    Preliminaries}\label{subsec:intro-and-prior-work-prelims}

\bnlem[Symmetrization with Rademacher random
    variables]\label{nlem:indep-symmetric-rademacher} Suppose that $\varepsilon$
    is a symmetric distribution \ie $\varepsilon \eqd - \varepsilon$, $\xi \sim
    \distRademacher{\parens{\alpha}}$, with $\alpha \in [0, 1]$. If $\xi \indep
    \varepsilon$ then $\xi \varepsilon \eqd \varepsilon$.
\enlem

\bprfof{\Cref{nlem:indep-symmetric-rademacher}}\label{prf:indep-symmetric-rademacher}
Let us define $Q \defined \xi \varepsilon$. We then have the following:
\begin{align*}
    \Prb{Q \geq q}
     & \defined \Prb{\xi \varepsilon \geq q}
    \tag{since $Q \defined \xi \varepsilon$.}                                    \\
     & = \Prbb{\xi \varepsilon \geq q}{\xi = - 1}\Prb{\xi = -1} +
    \Prbb{\xi \varepsilon \geq q}{\xi = 1}\Prb{\xi = 1}
    \tag{since $\xi \sim \distRademacher{\parens{\alpha}}$.}                     \\
     & = \Prb{-\varepsilon \geq q}(1 - \alpha) +
    \Prb{\varepsilon \geq q}(\alpha)
    \tag{since $\xi \indep \varepsilon$.}                                        \\
     & = \Prb{\varepsilon \geq q}(1 - \alpha) + \Prb{\varepsilon \geq q}(\alpha)
    \tag{since $\varepsilon \eqd - \varepsilon$.}                                \\
     & = \Prb{\varepsilon \geq q}\parens{1 - \alpha + \alpha}                    \\
     & = \Prb{\varepsilon \geq q}
\end{align*}
So we have that $Q \defined \xi \varepsilon \eqd \varepsilon$, as required.
\eprfof

The setting can be simplified if the adversary chooses the sign-corruptions
independent of the error terms. To see this, first note that $\varepsilon_{i}$
are centered (\ie symmetric) Gaussian random variables. Now, if the $(\xi_1,
\ldots, \xi_n)$ are picked independently from $(\varepsilon_1, \ldots,
\varepsilon_n)$, the \texttt{ASCI} generating process response reduces to $R_{i}
= \xi_{i} \mu_{i} + \varepsilon_{i}$. That is our setting encompasses this more
simplified setting, and is shown formally in
\Cref{ncor:reduced-adversarial-gen-process}. Further, we note that in the case
where $\xi_{i} \eqas 1$ then and $\mu_1 \leq \mu_{2} \leq \ldots \leq \mu_{n}$,
then this is equivalent to the standard univariate isotonic regression setup.

\bncor\label{ncor:reduced-adversarial-gen-process} In the case where $\xi_{i}
    \indep \varepsilon_{i}$, for each $i \in [n]$ we have that the \texttt{ASCI}
    generating process simplifies to $R_{i} = \xi_{i}\mu_{i} + \varepsilon_{i}$.
\encor

\bprfof{\Cref{ncor:reduced-adversarial-gen-process}}\label{prf:reduced-adversarial-gen-process}
We note that the underlying adversarial generating process is given by $R_{i} =
\xi_{i}(\mu_{i} + \varepsilon_{i}) = \xi_{i}\mu_{i} + \xi_{i}\varepsilon_{i}$,
for each $i \in [n]$. Now since $\xi_{i} \indep \varepsilon_{i}$ we have by
applying \Cref{nlem:indep-symmetric-rademacher} for each $i \in [n]$ that
$\xi_{i} \varepsilon_{i} \distiid \varepsilon_{i}$. And so the required
adversarial model can be written as $R_{i} = \xi_{i}\mu_{i} + \varepsilon_{i}$,
as required.
\eprfof

\subsection{Important Model Definitions}\label{subsec:equiv-conds-model-defns}

First, we formally (redefine) the generating model described in
\Cref{nexa:asci-rademacher-mixture}.

\begin{restatable}[Two-component Gaussian mixture \texttt{ASCI} special case
        from
        \Cref{nexa:asci-rademacher-mixture}]{ndefn}{ndefnascirademachermixture}\label{ndefn:asci-rademacher-mixture}
        We consider $n$ observations, $\thesetb{\resp_{i}}{i \in [n]}$, where
        each observation $\resp_{i}$ is generated from the following model:
    \begin{align}
        \resp_{i}
         & = \xi_{i}\mu_{i} + \varepsilon_{i}
        \label{neqn:asci-rademacher-mixture-01}                                 \\
        \text{s.t. } 0
         & < \eta \leq \mu_{1} \leq \mu_{2} \leq \ldots \leq \mu_{n}
        \label{neqn:asci-rademacher-mixture-02}                                 \\
        \text{ and } \varepsilon_{i}
         & \distiid \distNorm\parens{0, \sigma^2}
        \label{neqn:asci-rademacher-mixture-03}                                 \\
        \text{ and } \xi_{i}
         & \distiid \distRademacher{(p)}, \, p \in (0, 1), \text{ and } \xi_{i}
        \indep \varepsilon_{i}
        \label{neqn:asci-rademacher-mixture-04}
    \end{align}
\end{restatable}

\nit Second, we formally define the generating model described in
\Cref{nexa:equiv-conds}.

\begin{restatable}[Non-convex generating model from
        \Cref{nexa:equiv-conds}]{ndefn}{ndefnnonconvexisotonicadversarialsignmodel}\label{ndefn:non-convex-isotonic-adversarial-sign-model}
        We consider $n$ observations, $\thesetb{\resp_{i}}{i \in [n]}$, where
        each observation $\resp_{i}$ is generated from the following model:
    \begin{align}
        \resp_{i}
         & = \gamma_{i} + \varepsilon_{i}
        \label{neqn:non-convex-adversarial-iso-reg-1}                                              \\
        \text{s.t. } 0
         & < \eta \leq \absa{\gamma_{1}} \leq \absa{\gamma_{2}} \leq \ldots \leq \absa{\gamma_{n}}
        \label{neqn:non-convex-adversarial-iso-reg-2}                                              \\
        \text{ and } \varepsilon_{i}
         & \distiid \distNorm\parens{0,\sigma^2}
        \label{neqn:non-convex-adversarial-iso-reg-3}
    \end{align}
\end{restatable}

\bnrmk\label{nrmk:non-convex-isotonic-adversarial-sign-model} From a simulation
perspective, each $\gamma_{i}$ is generated first subject to
\Cref{neqn:non-convex-adversarial-iso-reg-2}, then $\xi_{i}$ is sampled
independently, and both are added to give each response $R_{i}$.
\enrmk

\nit Third, we introduce an alternative model as per
\Cref{ndefn:alternative-non-convex-model}.

\begin{restatable}[Alternative non-convex
        model]{ndefn}{ndefnalternativenonconvexmodel}\label{ndefn:alternative-non-convex-model}
    \begin{align}
        \resp_{i}
         & = \xi_{i} a_{i} + \varepsilon_{i}
        \label{neqn:non-convex-equiv-adversarial-iso-reg-1}    \\
        \text{s.t. } 0
         & < \eta \leq a_{1} \leq a_{2} \leq \ldots \leq a_{n}
        \label{neqn:non-convex-equiv-adversarial-iso-reg-2}    \\
        \text{ and } \varepsilon_{i}
         & \distiid \distNorm\parens{0,\sigma^2}
        \label{neqn:non-convex-equiv-adversarial-iso-reg-3}    \\
        \text{ and } \xi_{i}
         & = \sgn{(\gamma_{i})}
        \label{neqn:non-convex-equiv-adversarial-iso-reg-4}    \\
        \text{ and } \xi_{i}
         & \indep \varepsilon_{i}
        \label{neqn:non-convex-equiv-adversarial-iso-reg-5}    \\
        \text{ and } a_{i}
         & = \absa{\gamma_{i}}
        \label{neqn:non-convex-equiv-adversarial-iso-reg-6}
    \end{align}
\end{restatable}

\nit Finally, for convenience we recall
\Cref{ndefn:isotonic-adversarial-sign-model} as follows.

\ndefnisotonicadversarialsignmodel*

\subsection{Proof justification for
    \Cref{nexa:asci-rademacher-mixture}}\label{subsec:asci-rademacher-mixture-proofs}

\bnprop[Justification for
    \Cref{nexa:asci-rademacher-mixture}]\label{nprop:asci-rademacher-mixture-proofs}
    Under the model generating process described in
    \Cref{nexa:asci-rademacher-mixture} (\ie, per
    \Cref{ndefn:asci-rademacher-mixture}), the following model definition
    inclusion holds.
\begin{equation}\label{asci-rademacher-mixture-proofs-1}
    \textnormal{\Cref{ndefn:asci-rademacher-mixture}}
    \subseteq
    \textnormal{\Cref{ndefn:isotonic-adversarial-sign-model}}
\end{equation}
\enprop

\bnrmk\label{nrmk:asci-rademacher-mixture-proofs} Here, each definitional
inclusion is to be read as the former generating model definition being a
special case of the latter generating model definition.
\enrmk

\bprfof{\Cref{nexa:asci-rademacher-mixture}}\label{prf:asci-rademacher-mixture}
Our basic strategy is to show each model inclusion in turn.\\
\newline
(\Cref{ndefn:asci-rademacher-mixture} $\subseteq$
\Cref{ndefn:isotonic-adversarial-sign-model}). Observe that
\Cref{neqn:asci-rademacher-mixture-02,neqn:asci-rademacher-mixture-03} are
definitionally equivalent to
\Cref{neqn:adversarial-iso-reg-2,neqn:adversarial-iso-reg-3}, respectively.
Moreover we have that \Cref{neqn:asci-rademacher-mixture-04} is a special case
of \Cref{neqn:adversarial-iso-reg-4}. Finally, from
\Cref{neqn:asci-rademacher-mixture-04} we have that $\xi_{i} \distiid
\distRademacher{(p)}, \, p \in (0, 1)$, and $\xi_{i} \indep \varepsilon_{i}$.
Thus from \Cref{ncor:reduced-adversarial-gen-process}, it follows that
\Cref{neqn:asci-rademacher-mixture-01} is a special case of
\Cref{neqn:asci-rademacher-mixture-01}. \\
\newline
In summary we have shown that \Cref{asci-rademacher-mixture-proofs-1} holds,
from which it follows that \Cref{nexa:asci-rademacher-mixture} (or equivalently
\Cref{ndefn:asci-rademacher-mixture}) is a special case of
\Cref{ndefn:isotonic-adversarial-sign-model}, as required.
\eprfof

\subsection{Proof justification for
    \Cref{nexa:equiv-conds}}\label{subsec:equiv-conds-proofs}

We now provide a formal proof justification that \Cref{nexa:equiv-conds} is a
special case of the generating process described in
\Cref{ndefn:isotonic-adversarial-sign-model}.

\bnprop[Justification for
    \Cref{nexa:equiv-conds}]\label{nprop:equiv-conds-proofs} Under the model
    generating process described in \Cref{nexa:equiv-conds}, the following model
    definition inclusion holds.
\begin{equation}\label{equiv-conds-proofs-1}
    \textnormal{\Cref{ndefn:non-convex-isotonic-adversarial-sign-model}}
    =
    \textnormal{\Cref{ndefn:alternative-non-convex-model}}
    \subseteq
    \textnormal{\Cref{ndefn:isotonic-adversarial-sign-model}}
\end{equation}
\enprop

\bnrmk\label{nrmk:equiv-conds-proofs} As with
\Cref{nprop:asci-rademacher-mixture-proofs}, each definitional inclusion is to
be read as the former generating model definition being a special case of the
latter generating model definition. In the case of equivalence, we note that
that both inclusions hold between the model definitions.
\enrmk

\bprfof{\Cref{nexa:equiv-conds}}\label{prf:equiv-conds}  Our basic strategy is
to show each model inclusion in turn.\\
\newline
(\Cref{ndefn:non-convex-isotonic-adversarial-sign-model} $=$
\Cref{ndefn:alternative-non-convex-model}) This follows by construction. Observe
that
\Cref{neqn:non-convex-equiv-adversarial-iso-reg-4,neqn:non-convex-equiv-adversarial-iso-reg-6}
imply that $\xi_{i} a_{i} = \sgn{(\gamma_{i})} \absa{\gamma_{i}} = \gamma_{i}$,
so that
\Cref{neqn:non-convex-adversarial-iso-reg-1,neqn:non-convex-equiv-adversarial-iso-reg-1}
are equivalent. In addition from
\Cref{neqn:non-convex-equiv-adversarial-iso-reg-6}, we have that $a_{i} =
\absa{\gamma_{i}}$ and thus
\Cref{neqn:non-convex-adversarial-iso-reg-2,neqn:non-convex-equiv-adversarial-iso-reg-2}
are equivalent, as are
\Cref{neqn:non-convex-adversarial-iso-reg-3,neqn:non-convex-equiv-adversarial-iso-reg-3}.
As such the equality is established between the two generating model
definitions.\\
\newline
(\Cref{ndefn:alternative-non-convex-model} $\subseteq$
\Cref{ndefn:isotonic-adversarial-sign-model}). Observe that by
\Cref{neqn:adversarial-iso-reg-3,neqn:non-convex-equiv-adversarial-iso-reg-3}
are definitionally equivalent. Observe from
\Cref{neqn:non-convex-equiv-adversarial-iso-reg-4} that $\xi_{i} =
\sgn{(\gamma_{i})} \in \theseta{-1, 1}$ which is a special case of
\Cref{neqn:adversarial-iso-reg-4}. For each observation $i \in [n]$ using
\Cref{neqn:non-convex-equiv-adversarial-iso-reg-6} that by setting $a_{i}
\defined \mu_{i}$ that
\Cref{neqn:adversarial-iso-reg-2,neqn:non-convex-equiv-adversarial-iso-reg-2}
are equivalent. Finally since $\xi_{i} \indep \varepsilon_{i}$ from
\Cref{neqn:non-convex-equiv-adversarial-iso-reg-5}, we note that
\Cref{neqn:non-convex-equiv-adversarial-iso-reg-1} is a special case of
\Cref{neqn:adversarial-iso-reg-1} by applying
\Cref{ncor:reduced-adversarial-gen-process} to the observation $R_{i}$, for each
$i \in [n]$.\\
\newline
In summary we have shown that \Cref{equiv-conds-proofs-1} holds, from which it
follows that \Cref{ndefn:non-convex-isotonic-adversarial-sign-model} is a
special case of \Cref{ndefn:isotonic-adversarial-sign-model}, as required.
\eprfof

\newpage
\section{Proofs of
  \Cref{sec:ascifit-three-step-est-proc}}\label{sec:ascifit-three-step-est-proc-proofs}

\subsection{Mathematical
    Preliminaries}\label{subsec:ascifit-three-step-est-proc-prelims}

\bnthm[Projection onto the nonnegative monotone
    cone]\label{nthm:nemeth-projection-nonneg-cone} Suppose that
    $\mclS^{\uparrow} \subseteq \reals^{n}$ is the monotone cone, that is,
\begin{equation}
    \mclS^{\uparrow}
    \defined \thesetb{\boldsymbol{\mu} \defined (\mu_1, \ldots,
        \mu_{n})^{\top} \in \reals^{n}}{\mu_1 \leq \ldots \leq \mu_{n}}. \nonumber
\end{equation}
and $\mclS^{\uparrow}_{+} \subseteq \reals^{n}$ is the nonnegative monotone
cone, that is,
\begin{equation}
    \mclS^{\uparrow}_{+} \defined \thesetb{\boldsymbol{\mu} \defined
        \parens{\mu_1, \ldots, \mu_{n}}^{\top} \in \mclS^{\uparrow}}{\mu_{1} \geq 0}. \nonumber
\end{equation}
Then for an arbitrary $\bfv \in \reals^{n}$ it holds that
\begin{equation}
    \projn{\mclS^{\uparrow}_{+}}{\bfv} = \parens{\projn{\mclS^{\uparrow}}{\bfv}}^{+}, \nonumber
\end{equation}
where for any $\bfz \in \reals^{n}$, $\bfz^{+} \in \reals^{n}$ stands for the
lattice operation defined by the order induced by the nonnegative orthant in
$\reals^{n}$. That is, we define the operation componentwise as
$\parens{\bfz^{+}}_{i} \defined \parens{\bfz}_{i} \vee 0$ for each component
index $i \in [n]$.
\enthm

\bprfof{\Cref{nthm:nemeth-projection-nonneg-cone}}\label{prf:nemeth-projection-nonneg-cone}
See \citet[Corollary~1]{nemeth2012projectmonotonenonnegconepava} for details.
\eprfof

\bnrmk\label{nrmk:nemeth-projection-nonneg-cone} In effect,
\Cref{nthm:nemeth-projection-nonneg-cone} basically states that in order to
project onto the nonegative monotone cone, $K$, one can instead first project
onto the monotone cone, $W$, first, and then take the non-negative part along
each component. This is useful, since one can leverage algorithms like the PAVA
which already efficiently handle projection onto the unrestricted monotone cone,
$W$.
\enrmk

\subsection{Proof of
    \Cref{nprop:intuition-for-correction}}\label{subsec:intuition-for-correction-proofs}

\npropintuitionforcorrection*
\bprfof{\Cref{nprop:intuition-for-correction}}\label{prf:intuition-for-correction}
This follows along the following lines. First subtract $f(\eta, \sigma)$ from
all $\tilde \absresp_{i}$ to bring the first problem to

\begin{align}
    \argmin_{\bar T_i} \sum_{i = 1}^{n} ((\absresp_{i} - f(\eta, \sigma)) - \bar \absresp_{i})^{2} \mbox{ s.t. } 0 \leq \bar T_1 \leq \ldots \leq \bar T_n,
\end{align}

where $\bar \absresp_{i} = \tilde{\absresp}_{i} - f(\eta, \sigma)$. Now the
solution to the unrestricted problem

\begin{align*}
    \argmin_{T_i^*} \sum_{i = 1}^{n} ((\absresp_{i} - f(\eta, \sigma)) -  \absresp^*_{i})^{2} \mbox{ s.t. } T^*_1 \leq \ldots \leq T^*_n,
\end{align*}
is $T^{*}_{i} = \hat {\absresp}_{i} - f(\eta, \sigma)$.
Next we apply \Cref{nthm:nemeth-projection-nonneg-cone}, we see that $\bar
    \absresp_{i} = T^{*}_{i}  \vee 0$, so that $ \tilde{\absresp}_{i} = \bar
    \absresp_{i} + f(\eta, \sigma) =T^{*}_{i} \vee 0 + f(\eta,\sigma) =
    (T^{*}_{i} + f(\eta, \sigma)) \vee  f(\eta, \sigma) = \hat{\absresp}_{i}
    \vee f(\eta, \sigma)$ which is what we wanted to show.
\eprfof

\newpage
\section{Proofs of
  \Cref{sec:step-two-unique-root}}\label{sec:step-two-unique-root-proofs}

\subsection{Mathematical
    Preliminaries}\label{subsec:step-two-unique-root-prelims}

\nit The key idea to prove this theorem here is to apply
\citep[Theorem~2.2(ii)]{zhang2002riskboundsisotonicreg} to our specific setting.
To ensure our work is self-contained, we translate this result into the notation
of our paper:

\bnthm[Theorem~2.2 (ii)
    \citep{zhang2002riskboundsisotonicreg}]\label{nthm:zhang-thm-22} Let $R_{n,
    p}(f, \boldsymbol{\mu}, \sigma, \sigma_{p}) \defined \parens{\frac{1}{n}
    \sum_{i = 1}^{n}\Exp{\absa{\widehat{\absresp}_{i} - f(\mu_{i},
    \sigma)}^{p}}}^{\frac{1}{p}}$. Let $\delta_{i} \defined \absresp_{i} -
    f(\mu_{i}, \sigma)$ be independent, with $\Exp{\delta_{i}} = 0$ and
    $\Exp{\absa{\delta_{i}}^{p \vee 2}} \leq \sigma_{p}^{p \vee 2}$, $p \geq 1$
    then:
\begin{align}
    R_{n, p}(f, \boldsymbol{\mu}, \sigma, \sigma_{p})
     & \leq 2^{\frac{1}{p}} \sigma_{p} C_{p}
    \min \brackets{1,
        \frac{3}{2}\braces{\frac{3}{(3-p)_{+}}\parens{\frac{V(f, \boldsymbol{\mu}, \sigma)}{n \sigma_{p} C_{p}}}^{\frac{p}{3}}
            + \frac{1}{n} \int_{0}^{n} \frac{d x}{(x \vee 1)^{\frac{p}{2}}}}^{\frac{1}{p}}
    }
\end{align}
where $C_{p}$ are constants depending on p only in general.
\enthm

\bprfof{\Cref{nthm:zhang-thm-22}}\label{prf:zhang-thm-22} See
\citet[Theorem~2.2(ii)]{zhang2002riskboundsisotonicreg} for details. Note that
to translate between our notation and theirs respectively, we have $\absresp_{i}
\equiv y_{i}, \widehat{T}_{i} \equiv \widehat{f}_{n}(t_{i}), f(\mu_{i}, \sigma)
\equiv f(t_{i}), \delta_{i} \equiv \varepsilon_{i}$ for each $i \in [n]$.
\eprfof

\bncor[Upper bound for $R_{n,2}^{2}(f, \boldsymbol{\mu}, \sigma,
        \sigma_{2})$]\label{ncor:zhang-thm-22} In our setting, define $X
        \defined \frac{1}{n} \sum_{i = 1}^{n} \parens{\widehat{\absresp}_{i} -
        f(\mu_{i}, \sigma)}^2$. We then have:
\begin{align}
    \Exp{X}
     & \leq \min \brackets{2 \sigma^{2} C_{2}^{2},
        \frac{27}{4}\parens{\frac{\mu_{n} - \mu_{1}}{n}}^{\frac{2}{3}} (\sigma C_{2})^{\frac{4}{3}}
    + \frac{2 \sigma^{2} C_{2}^{2}}{n} \parens{1 + \log{n}}} \label{neqn:upp-bound-Rn2-squared} \\
     & \defines r_{n,2}(\mu_{n}, \mu_{1}, \sigma)       \nonumber
\end{align}
where $C_{2}$ is a constant.
\encor

\bprfof{\Cref{ncor:zhang-thm-22}}\label{prf:ncor-zhang-thm-22} Since $X \defined
    \frac{1}{n} \sum_{i = 1}^{n} \parens{\widehat{\absresp}_{i} - f(\mu_{i},
    \sigma)}^2$, then $X = R_{n,2}^{2}(f, \boldsymbol{\mu}, \sigma,
    \sigma_{2})^{2}$, by definition in the setting of \Cref{nthm:zhang-thm-22},
    assuming the relevant sufficient conditions are met. We now need to check
    the sufficient condition for \Cref{nthm:zhang-thm-22}. Here we have, for
    each $i \in [n]$, that $\delta_{i} \defined \absresp_{i} - f(\mu_{i},
    \sigma)$. Note that by definition $\Exp{\delta_{i}} = \Exp{\absresp_{i}} -
    f(\mu_{i}, \sigma) = 0$. We observe that $(\delta_{1}, \ldots, \delta_{n})$
    are independent since the original responses, \ie $(\resp_{1}, \ldots,
    \resp_{n})$ are independent by assumption. And taking absolute values and
    centering are measurable transformations which preserve their independence.
    We note that as per \citet[Theorem~2.2(ii)]{zhang2002riskboundsisotonicreg},
    we are required to check the sufficient condition $\Exp{\absa{\delta_{i}}^{p
    \vee 2}} \leq \sigma_{p}^{p \vee 2}$. In our case, with $p = 2$, this is
    equivalent to showing that $\Exp{\delta_{i}^{2}} \leq \sigma_{2}^{2}$. Then
    for each $i \in [n]$ we have:
\begin{align}
    \Exp{\absa{\delta_{i}}^{p \vee 2}}
     & = \Exp{\delta_{i}^{2}}
    \tag{since $p = 2$.}                                              \\
     & = \Var{\absresp_{i}}
    \tag{since $\delta_{i}$ are mean centered $\absresp_{i}$ values.} \\
     & \defines g(\mu_{i}, \sigma)
    \tag{by definition.}                                              \\
     & \leq \sigma^{2}
    \tag{using \Cref{neqn:fold-norm-var-ge-sigma-sq}}                 \\
     & \defines \sigma_{2}^{2}
    \label{neqn:sigma-sub-2-zhang}
\end{align}
As required, by defining $\sigma_{2} \defined \sigma$. So we meet this
sufficient condition. Additionally observe that
\begin{align}
    \int_{0}^{n} \frac{d x}{(x \vee 1)}
     & = \int_{0}^{1} \frac{d x}{(x \vee 1)} + \int_{1}^{n} \frac{d x}{(x \vee 1)}
    \tag{by truncation} \nonumber                                                  \\
     & = \int_{0}^{1} d x + \int_{1}^{n} \frac{d x}{x} \nonumber                   \\
     & = 1 + \log{n} \label{neqn:upp-bound-Rn2-squared-integral}
\end{align}

\nit Now, in our setting note that $V(f, \boldsymbol{\mu}, \sigma) \leq \mu_{n}
    - \mu_{1}$ using \Cref{neqn:fold-norm-mean-total-var}, it follows that:
\begin{align}
    \Exp{X}
     & \defined R_{n,2}^{2}(f, \boldsymbol{\mu}, \sigma, \sigma_{2})
    \tag{by definition.}                                                       \\
     & \leq \min \theseta{2 \sigma_{2}^{2} C_{2}^{2},
        \frac{27}{4}\parens{\frac{\mu_{n} - \mu_{1}}{n}}^{\frac{2}{3}} (\sigma_{2} C_{2})^{\frac{4}{3}}
        + \frac{2 \sigma_{2}^{2} C_{2}^{2}}{n} \int_{0}^{n} \frac{d x}{(x \vee 1)}}
    \tag{setting $p = 2$ in \Cref{nthm:zhang-thm-22}.}                         \\
     & = \min \theseta{2 \sigma_{2}^{2} C_{2}^{2},
        \frac{27}{4}\parens{\frac{\mu_{n} - \mu_{1}}{n}}^{\frac{2}{3}} (\sigma_{2} C_{2})^{\frac{4}{3}}
        + \frac{2 \sigma_{2}^{2} C_{2}^{2}}{n} \parens{1 + \log{n}}}
    \tag{using \Cref{neqn:upp-bound-Rn2-squared-integral}}                     \\
     & = \min \theseta{2 \sigma^{2} C_{2}^{2},
        \frac{27}{4}\parens{\frac{\mu_{n} - \mu_{1}}{n}}^{\frac{2}{3}} (\sigma C_{2})^{\frac{4}{3}}
        + \frac{2 \sigma^{2} C_{2}^{2}}{n} \parens{1 + \log{n}}}
    \tag{since $\sigma_{2} \defined \sigma$ per \Cref{neqn:sigma-sub-2-zhang}} \\
     & \defines r_{n,2}(\mu_{n}, \mu_{1}, \sigma)
\end{align}
\newline
as required.
\eprfof

\bnlem[Concentration for mean Folded Normal]\label{nlem:conc-folded-normal} In
our setting we assume that $\frac{1}{n} \sum_{i = 1}^{n} \mu_{i}^{2} \leq
\univupperboundmu$, for each $n \in \nats$. Define $X \defined \frac{1}{n}
\sum_{i = 1}^{n} \parens{\absresp_{i} - f(\mu_{i}, \sigma)}^{2}$. We then have:
\begin{equation}\label{neqn:conc-folded-normal-01}
    \absa{X - \Exp{X}}
    \leq  2 \gamma \sigma \sqrt{\frac{5 \sigma^{2} + 4 \univupperboundmu}{n}}
\end{equation}
with probability at least $1 - \gamma^{-2}$, where $\Exp{X} = \frac{1}{n}
    \sum_{i = 1}^{n} g(\mu_{i}, \sigma) = \frac{1}{n}\sum_{i = 1}^{n}
    \parens{\mu_{i}^{2} + \sigma^{2} - f(\mu_{i}, \sigma)^2}$.
\enlem

\bprfof{\Cref{nlem:conc-folded-normal}}\label{prf:conc-folded-normal} First we
determine $\Exp{X}$ as follows:
\begin{align*}
    \Exp{X}
     & \defined \Exp{\frac{1}{n} \sum_{i = 1}^{n} \parens{\absresp_{i} -
            f(\mu_{i}, \sigma)}^{2}}
    \tag{by definition of $X$}                                                                 \\
     & = \frac{1}{n} \sum_{i = 1}^{n} \Exp{\parens{\absresp_{i} -
    f(\mu_{i}, \sigma)}^{2}}                                                                   \\
     & = \frac{1}{n} \sum_{i = 1}^{n} \Var{\absresp_{i}}
    \tag{by since $f(\mu_{i}, \sigma) \defined \Exp{\absresp_{i}}$.}                           \\
     & = \frac{1}{n} \sum_{i = 1}^{n} g(\mu_{i}, \sigma)
    \tag{by since $g(\mu_{i}, \sigma) \defined \Var{\absresp_{i}}$.}                           \\
     & = \frac{1}{n} \sum_{i = 1}^{n} \parens{\mu_{i}^{2} + \sigma^{2} - f(\mu_{i}, \sigma)^2}
    \tag{using \Cref{neqn:folded-normal-var}}
\end{align*}
as required. Next we determine $\Var{X}$ as follows:
\begin{align*}
    \Var{X}
     & \defined \Var{\frac{1}{n} \sum_{i = 1}^{n} \parens{\absresp_{i} -
            f(\mu_{i}, \sigma)}^{2}}
    \tag{by definition of $X$}                                           \\
     & = \frac{1}{n^{2}} \sum_{i = 1}^{n} \Var{\parens{\absresp_{i} -
            f(\mu_{i}, \sigma)}^{2}}
    \tag{by the independence of $\absresp_{i}$}
\end{align*}
Now note that for each $i \in [n]$ we have:
\begin{align}
    \Var{\parens{\absresp_{i} - f(\mu_{i}, \sigma)}^{2}}
     & = \Var{\absresp_{i}^{2} + f(\mu_{i}, \sigma)^{2} - 2 f(\mu_{i}, \sigma) \absresp_{i}}
    \nonumber                                                                                 \\
     & = \Var{\absresp_{i}^{2} - 2 f(\mu_{i}, \sigma) \absresp_{i}}
    \tag{by translation invariance.}                                                          \\
     & \leq 2 \parens{\Var{\absresp_{i}^{2}} + \Var{2 f(\mu_{i}, \sigma) \absresp_{i}}}
    \tag{using \Cref{ncor:square-sum-inequality-rv}}                                          \\
     & = 2 \Var{\absresp_{i}^{2}} + 8 f(\mu_{i}, \sigma)^{2} \Var{\absresp_{i}}
    \nonumber                                                                                 \\
     & = 2 \Var{\absresp_{i}^{2}} + 8 f(\mu_{i}, \sigma)^{2} g(\mu_{i}, \sigma)
    \tag{since $g(\mu_{i}, \sigma) \defined \Var{\absresp_{i}}$}                              \\
     & \leq 2 (4 \mu_{i}^{2} \sigma^{2} + 2 \sigma^{4}) + 8 f(\mu_{i}, \sigma)^{2} \sigma^{2}
    \tag{using \Cref{neqn:fold-norm-var-ge-sigma-sq,neqn:fold-norm-sq-var}}                   \\
     & = 8 \mu_{i}^{2} \sigma^{2} + 8 f(\mu_{i}, \sigma)^{2} \sigma^{2} + 4 \sigma^{4}
    \nonumber                                                                                 \\
     & \leq 16 f(\mu_{i}, \sigma)^{2} \sigma^{2} + 4 \sigma^{4}
    \tag{using \Cref{neqn:fold-norm-mean-ge-norm}}                                            \\
     & \leq 16 (\mu_{i}^{2} + \sigma^{2}) \sigma^{2} + 4 \sigma^{4}
    \tag{using \Cref{neqn:fold-norm-mean-le-sec-mom-norm}}                                    \\
     & = 16 \mu_{i}^{2}\sigma^{2} + 20 \sigma^{4}
    \label{neqn:var-cent-fold-norm-squared}
\end{align}
Therefore we have that
\begin{align*}
    \Var{X}
     & = \frac{1}{n^{2}} \sum_{i = 1}^{n} \Var{\parens{\absresp_{i} -
    f(\mu_{i}, \sigma)}^{2}}                                                                              \\
     & \leq \frac{1}{n^{2}} \sum_{i = 1}^{n} \parens{16 f(\mu_{i}, \sigma)^{2} \sigma^{2} + 4 \sigma^{4}}
    \tag{using \Cref{neqn:fold-norm-mean-ge-norm}}                                                        \\
     & \leq \frac{1}{n^{2}} \sum_{i = 1}^{n} \parens{16 \mu_{i}^{2}\sigma^{2} + 20 \sigma^{4}}
    \tag{using \Cref{neqn:var-cent-fold-norm-squared}}                                                    \\
     & \leq \frac{16 C \sigma^{2} + 20 \sigma^{4}}{n}
    \tag{assuming $\frac{1}{n} \sum_{i = 1}^{n} \mu_{i}^{2} \leq C$, for each $n \in \nats$.}
\end{align*}

\nit From this it follows that:
\begin{align}
    \Prb{\absa{X - \Exp{X}} \geq t}
     & \leq \frac{\Var{X}}{t^{2}}
    \text{, \; for each $t > 0$}
    \tag{using Chebychev's inequality}                      \\
     & \leq \frac{16 C \sigma^{2} + 20 \sigma^{4}}{n t^{2}}
    \text{, \; for each $t > 0$}
\end{align}

\nit It then follows that by setting the upper bound ($\rhs$) to $\gamma^{-2}
    \in (0, 1)$, that
\begin{equation*}
    \frac{16 C \sigma^{2} + 20 \sigma^{4}}{n t^{2}}
    = \frac{1}{\gamma^{2}}
    \implies t
    = \gamma \sigma \sqrt{\frac{5 \sigma^{2} + 4 \univupperboundmu}{n}}
\end{equation*}

\nit We then have that $\absa{X - \Exp{X}} \leq 2 \gamma \sigma \sqrt{\frac{5
            \sigma^{2} + 4 \univupperboundmu}{n}}$, with probability at least $1
            - \gamma^{-2}$, as required.
\eprfof

\nit Our end goal is to to show a the following high probability result
described in \Cref{nthm:conc-fitted-folded-normal}.

\bnthm[Concentration of fitted Folded
    Normal]\label{nthm:conc-fitted-folded-normal}
\begin{equation}\label{neqn:conc-fitted-folded-normal}
    \frac{1}{n} \sum_{i = 1}^{n} \parens{\widehat{\absresp}_{i}
        - f(\mu_{i}, \sigma)}^2
    \leq \delta r_{n,2}(\mu_n, \mu_1, \sigma)
\end{equation}
with probability at least $1 - \delta^{-1}$.
\enthm

\bprfof{\Cref{nthm:conc-fitted-folded-normal}}\label{prf:conc-fitted-folded-normal}
First to simplify notation we let $X \defined \frac{1}{n} \sum_{i = 1}^{n}
\parens{\widehat{\absresp}_{i} - f(\mu_{i}, \sigma)}^2$ represent the quantity
of interest. Observe that $X \geq 0$ \as by definition, so that $\absa{X} \eqas
X$. Then for any $t > 0$ we have:
\begin{align}
    \Prb{X \geq t}
     & \leq \frac{\Exp{X}}{t}
    \tag{by Markov's inequality}                          \\
     & \leq    \frac{R_{n, 2}^{2}(f(\mu_{i}, \sigma))}{t}
    \tag{by definition, per \Cref{ncor:zhang-thm-22}}     \\
     & \leq   \frac{r_{n,2}(\mu_{n}, \mu_{1}, \sigma)}{t}
    \tag{using \Cref{ncor:zhang-thm-22}}
\end{align}
It then follows that by setting the upper bound ($\rhs$) to $\delta^{-1} \in (0,
    1)$, that
\begin{equation*}
    \frac{r_{n,2}(\mu_{n}, \mu_{1}, \sigma)}{t}
    = \frac{1}{\delta}
    \implies t
    = \delta r_{n,2}(\mu_{n}, \mu_{1}, \sigma)
\end{equation*}
We then have that $\absa{X} \eqas X \leq \delta r_{n,2}(\mu_{n}, \mu_{1},
    \sigma)$, with probability at least $1 - \delta^{-1}$, as required.
\eprfof

\bnlem[Concentration of $\frac{1}{n} \sum_{i = 1}^{n}
        \absresp_{i}^{2}$]\label{nlem:conc-samp-mean-ti-sq} In our setting,
        define $X \defined \frac{1}{n} \sum_{i = 1}^{n} \absresp_{i}^{2}$. We
        then have:
\begin{align}
    \absa{X - \Exp{X}}\label{neqn:conc-samp-mean-ti-sq-01}
    \leq 2 \gamma \sigma \sqrt{\frac{2 \sigma^{2} + 4 \univupperboundmu}{n}}
\end{align}
with probability at least $1 - \gamma^{-2}$, where $\Exp{X} = \frac{1}{n}
    \sum_{i = 1}^{n} \parens{\mu_{i}^{2} + \sigma^{2}}$.
\enlem

\bprfof{\Cref{nlem:conc-samp-mean-ti-sq}} Let $X \defined \frac{1}{n} \sum_{i =
        1}^{n} \absresp_{i}^{2}$. First we determine $\Exp{X}$ as follows:
\begin{align}
    \Exp{X}
     & \defined \Exp{\frac{1}{n} \sum_{i = 1}^{n} \absresp_{i}^{2}}
    \tag{by definition of $X$}                                                                                        \\
     & = \frac{1}{n} \sum_{i = 1}^{n} \Exp{\absresp_{i}^{2}}
    \tag{by linearity of expectation.}                                                                                \\
     & = \frac{1}{n} \sum_{i = 1}^{n} \parens{\Var{\absresp_{i}} + \parens{\Exp{\absresp_{i}}}^{2}}
    \nonumber                                                                                                         \\
     & = \frac{1}{n} \sum_{i = 1}^{n} \parens{\mu_{i}^{2} + \sigma^{2} - f(\mu_{i}, \sigma)^2 + f(\mu_{i}, \sigma)^2}
    \tag{using \Cref{neqn:folded-normal-exp-1,neqn:folded-normal-var}}                                                \\
     & = \frac{1}{n} \sum_{i = 1}^{n} \parens{\mu_{i}^{2} + \sigma^{2}}
    \nonumber
\end{align}
as required. Next we determine $\Var{X}$ as follows:
\begin{align}
    \Var{X}
     & \defined \Var{\frac{1}{n} \sum_{i = 1}^{n} \absresp_{i}^{2}}
    \tag{by definition of $X$.}                                                            \\
     & = \frac{1}{n^{2}} \sum_{i = 1}^{n} \Var{\absresp_{i}^{2}}
    \tag{by the independence of $\absresp_{i}$.}                                           \\
     & = \frac{1}{n^{2}} \sum_{i = 1}^{n} \parens{4 \mu_{i}^{2} \sigma^{2} + 2 \sigma^{4}}
    \tag{using \Cref{neqn:fold-norm-sq-var}}                                               \\
     & \leq \frac{4 \univupperboundmu \sigma^{2} + 2 \sigma^{4}}{n}
    \tag{assuming $\frac{1}{n} \sum_{i = 1}^{n} \mu_{i}^{2} \leq C$, for each $n \in \nats$.}
\end{align}

\nit From this it follows that:
\begin{align}
    \Prb{\absa{X - \Exp{X}} \geq t}
     & \leq \frac{\Var{X}}{t^{2}}
    \text{, \; for each $t > 0$}
    \tag{using Chebychev's inequality}                                    \\
     & \leq \frac{4 \univupperboundmu \sigma^{2} + 2 \sigma^{4}}{n t^{2}}
    \text{, \; for each $t > 0$}
\end{align}

\nit It then follows that by setting the upper bound ($\rhs$) to $\gamma^{-2}
    \in (0, 1)$, that
\begin{equation*}
    \frac{4 \univupperboundmu \sigma^{2} + 2 \sigma^{4}}{n t^{2}}
    = \frac{1}{\gamma^{2}}
    \implies t
    = \gamma \sigma \sqrt{\frac{2 \sigma^{2} + 4 \univupperboundmu}{n}}
\end{equation*}

\nit We then have that $\absa{X - \Exp{X}} \leq 2 \gamma \sigma \sqrt{\frac{2
            \sigma^{2} + 4 \univupperboundmu}{n}}$, with probability at least $1
            - \gamma^{-2}$, as required.
\eprfof

\subsection{Proof of
    \Cref{nthm:eqn-has-unique-root}}\label{subsec:eqn-has-unique-root-proofs}

\nthmeqnhasuniqueroot*
\bprfof{\Cref{nthm:eqn-has-unique-root}}\label{prf:eqn-has-unique-root}

First, under the \texttt{ASCIFIT} setup, we can rewrite
\Cref{neqn:second-moment-matching-step} as $H(\sigma) = 0$, where:

\begin{align}
    H(\sigma)
     & \defined G(\sigma)
    - \frac{1}{n} \sum_{i = 1}^{n} \absresp_{i}^2.
    \label{neqn:h-sigma-function-ascifit} \\
    G(\sigma)
     & \defined
    \sigma^{2} +
    \frac{1}{n} \sum_{i = 1}^{n} (f^{-1}(\widehat{\absresp}_{i} \vee f(\eta, \sigma), \sigma))^2
    \label{neqn:g-sigma-function-ascifit}
\end{align}

\nit Our goal in this proof is to show that $H(\sigma) = 0$ has a solution
$\sigma^{*} \in \brackets{0,\sqrt{\frac{1}{n} \sum_{i = 1}^{n}
\absresp_{i}^2}}$, which occurs with high probability. We note that per
\Cref{nlem:prop-g-sigma} that $G(\sigma)$ is increasing for $\sigma \geq 0$ and
\textit{strictly} increasing for $\sigma > 0$ (per \Cref{neqn:prop-g-sigma-02}).
Moreover to see that the equation $H(\sigma) = 0$ has a unique root we appeal to
the Intermediate Value Theorem. Specifically we are required to find two values
for $\sigma$, \ie $\theset{\sigma_{1}, \sigma_{2}}$, such that the following
conditions hold:
\begin{align}
    G(\sigma_{2})
     & \geq \frac{1}{n} \sum_{i = 1}^n \absresp_{i}^{2}
    \label{neqn:g-sigma-condn-higher}                   \\
    G(\sigma_{1})
     & \leq \frac{1}{n} \sum_{i = 1}^n \absresp_{i}^{2}
    \label{neqn:g-sigma-condn-lower}
\end{align}

By taking $\sigma_{2} \defined \sqrt{\frac{1}{n}\sum_{i = 1}^{n}
        \absresp_{i}^{2}}$, we observe that \as:
\begin{align}
    G(\sigma_{2})
     & = \frac{1}{n}\sum_{i = 1}^{n} \absresp_{i}^{2} +
    \ubrace{\frac{1}{n} \sum_{i = 1}^{n}
        \parens{f^{-1}\parens{\widehat{\absresp}_{i} \vee f\parens{\eta, \frac{1}{n}\sum_{i = 1}^{n} \absresp_{i}^{2}},
    \frac{1}{n}\sum_{i = 1}^{n} \absresp_{i}^{2}}}^2}{\geq 0 \text{ \as}}{} \\
     & \geq \sum_{i = 1}^{n} \absresp_{i}^{2}
\end{align}

So indeed $\sigma_{2} \defined \sqrt{\frac{1}{n}\sum_{i = 1}^{n}
        \absresp_{i}^{2}}$ satisfies the required condition in
        \Cref{neqn:g-sigma-condn-higher}. We now claim that $\sigma_{1} \defined
        0$ will satisfy \Cref{neqn:g-sigma-condn-lower}. First observe that:
\begin{equation}\label{neqn:g-sigma-value-zero}
    G(0)
    = \frac{1}{n} \sum_{i = 1}^{n} f^{-1}(\widehat{\absresp}_{i} \vee f(\eta, 0), 0)^{2}
    = \frac{1}{n} \sum_{i= 1}^{n} (\widehat{\absresp}_{i} \vee \eta)^2,
\end{equation}
we then want to show that
\begin{equation}\label{neqn:g-sigma-value-zero-goal}
    \frac{1}{n} \sum_{i= 1}^{n} (\widehat{\absresp}_{i} \vee \eta)^2
    \leq \frac{1}{n} \sum_{i = 1}^{n} \absresp^{2}_{i},
\end{equation}
holds with high probability, to be specified later.

Furthermore, since $\widehat{\absresp}_{i} \vee \eta$ is the solution to an
optimization problem we have that \as:
\begin{equation}\label{neqn:eqn-has-unique-root-proj}
    \sum_{i = 1}^{n}(\widehat{\absresp}_{i} \vee \eta - \eta)(\absresp_{i} - \eta)
    = \sum_{i = 1}^{n}(\widehat{\absresp}_{i} \vee \eta - \eta)^2.
\end{equation}

We see that \Cref{neqn:eqn-has-unique-root-proj} holds since when you project
any vector $\bfv \in \reals^{n}$ on a monotone cone $K \subseteq \reals^{n}$,
then $\projn{K}{\bfv}^{\top} \bfv = \norma{\projn{K}{\bfv}}_{2}^{2}$ per
\citet[Equation~1.16]{bellec2018sharporacleineqsshprestreg}. Specifically, in
our case we have that $K = \mclS^{\uparrow}_{+} \defined
\thesetc{\boldsymbol{\mu} \defined \parens{\mu_1, \ldots, \mu_{n}}^{\top} \in
\reals^{n}}{0 \leq \mu_{1} \leq \mu_{2} \leq \ldots \leq \mu_{n}}$, and $\bfv
\defined \parens{T_1 - \eta, \ldots, T_n - \eta}^{\top}$. We further observe
that \Cref{neqn:eqn-has-unique-root-proj} can be rewritten as follows \as:
\begin{align}
    \sum_{i = 1}^{n}(\widehat{\absresp}_{i} \vee \eta - \eta)^2
     & =  \sum_{i = 1}^{n}(\widehat{\absresp}_{i} \vee \eta - \eta)(\absresp_{i} - \eta)
    \tag{\Cref{neqn:eqn-has-unique-root-proj}}                                           \\
    \iff \sum_{i = 1}^{n}(\widehat{\absresp}_{i} \vee \eta)^2
    + 2 \eta \sum_{i = 1}^{n}(\widehat{\absresp}_{i} \vee \eta - \eta)
    + \eta^{2}
     & =  \sum_{i = 1}^{n}(\widehat{\absresp}_{i} \vee \eta)\absresp_{i}
    - \eta \sum_{i = 1}^{n}(\widehat{\absresp}_{i} \vee \eta - \eta)
    \nonumber                                                                            \\
     & - \eta \sum_{i = 1}^{n}\widehat{\absresp}_{i}
    + \eta^{2}
    \tag{expanding $\lhs / \rhs$.}                                                       \\
    \iff \sum_{i = 1}^{n}(\widehat{\absresp}_{i} \vee \eta)^2
     & = \sum_{i = 1}^{n}(\widehat{\absresp}_{i} \vee \eta)\absresp_{i}
    - \eta \parens{\sum_{i = 1}^{n}\absresp_{i}
        - \sum_{i = 1}^{n}\widehat{\absresp}_{i} \vee \eta}
    \label{neqn:eqn-has-unique-root-proj-01}                                             \\
    \iff \sum_{i = 1}^{n}(\widehat{\absresp}_{i} \vee \eta)^2
     & = \sum_{i = 1}^{n}(\widehat{\absresp}_{i} \vee \eta)\absresp_{i}
    - \eta \parens{\sum_{i = 1}^{n}\widehat{\absresp}_{i}
        - \sum_{i = 1}^{n}\widehat{\absresp}_{i} \vee \eta},
    \label{neqn:eqn-has-unique-root-proj-02}
\end{align}
where to go from \Cref{neqn:eqn-has-unique-root-proj-01} to
\Cref{neqn:eqn-has-unique-root-proj-02} we used the fact that $\sum_{i =
1}^{n}\absresp_{i} = \sum_{i =1}^{n}\widehat{\absresp}_{i}$. This holds since we
know that $\widehat{\absresp}_{i}$ are the \texttt{PAVA} solutions. Now we
derive the following upper bound \as:
\begin{align}
     & \frac{1}{n}\parens{\sum_{i = 1}^{n}\widehat{\absresp}_{i} - \sum_{i = 1}^{n}\widehat{\absresp}_{i} \vee \eta}
    \nonumber                                                                                                        \\
     & = \frac{1}{n} \sum_{i = 1}^{n}\widehat{\absresp}_{i}
    - \frac{1}{n} \sum_{i = 1}^{n}f(\mu_{i},\sigma)
    + \frac{1}{n} \sum_{i = 1}^{n}f(\mu_{i}, \sigma)
    - \frac{1}{n} \sum_{i = 1}^{n}\widehat{\absresp}_{i} \vee \eta
    \nonumber                                                                                                        \\
     & = \frac{1}{n}\sum_{i = 1}^{n} \parens{\widehat{\absresp}_{i} - f(\mu_{i},\sigma)}
    + \frac{1}{n} \sum_{i = 1}^{n} \parens{f(\mu_{i}, \sigma) - \widehat{\absresp}_{i} \vee \eta}
    \label{neqn:eqn-has-unique-root-proj-03}                                                                         \\
     & \leq \sqrt{\frac{1}{n}\sum_{i = 1}^{n}\parens{\widehat{\absresp}_{i} - f(\mu_{i},\sigma)}^{2}}
    + \sqrt{\frac{1}{n}\sum_{i = 1}^{n}\parens{\widehat{\absresp}_{i} \vee \eta - f(\mu_{i},\sigma)}^{2}}
    \label{neqn:eqn-has-unique-root-proj-04},
\end{align}
where the transition between
\Cref{neqn:eqn-has-unique-root-proj-03,neqn:eqn-has-unique-root-proj-04} was by
applying the Cauchy-Schwartz inequality to each summand. Note that for each $i
\in [n]$, we have that $f(\mu_{i}, \sigma) \geq \mu_{i} \geq \eta$ per
\Cref{neqn:adversarial-iso-reg-2,neqn:fold-norm-mean-ge-norm}. Then using
\Cref{nlem:square-max-diff-inequality} we have \as:

\begin{equation}\label{neqn:eqn-has-unique-root-proj-05}
    \frac{1}{n} \sum_{i= 1}^{n} \parens{(\widehat{\absresp}_{i} \vee \eta) - f(\mu_{i}, \sigma)}^{2}
    \leq \frac{1}{n} \sum_{i= 1}^{n} (\widehat{\absresp}_{i} - f(\mu_{i}, \sigma))^2.
\end{equation}

Hence by first using the inequality in \Cref{neqn:eqn-has-unique-root-proj-05}
to upper bound \Cref{neqn:eqn-has-unique-root-proj-04}, we can in turn upper
bound the $\lhs$ of \Cref{neqn:eqn-has-unique-root-proj-02} as follows \as:
\begin{equation}\label{neqn:eqn-has-unique-root-proj-06}
    \frac{1}{n} \sum_{i = 1}^{n}(\widehat{\absresp}_{i} \vee \eta)^2
    \leq \frac{1}{n}\sum_{i = 1}^{n}(\widehat{\absresp}_{i} \vee \eta)\absresp_{i}
    + 2\eta \sqrt{\frac{1}{n}\sum_{i = 1}^{n}(\widehat{\absresp}_{i} - f(\mu_{i},\sigma))^{2}}.
\end{equation}

\nit On the other hand we have by \Cref{nthm:conc-fitted-folded-normal} that:
\begin{equation}\label{neqn:conc-fitted-folded-normal-01}
    \frac{1}{n} \sum_{i = 1}^{n} \parens{\widehat{\absresp}_{i}
        - f(\mu_{i}, \sigma)}^2
    \leq \delta r_{n,2}(\mu_n, \mu_1, \sigma),
\end{equation}
with probability at least $1 - \delta^{-1}$, for $\delta^{-1} \in (0, 1)$. Thus
from \Cref{neqn:eqn-has-unique-root-proj-06}, we have:
\begin{align}
    \frac{1}{n} \sum_{i = 1}^{n}(\widehat{\absresp}_{i} \vee \eta)^2
     & \leq \frac{1}{n}\sum_{i = 1}^{n}(\widehat{\absresp}_{i} \vee \eta)\absresp_{i}
    + 2 \eta \parens{\delta r_{n,2}(\mu_n, \mu_1, \sigma)}^{\frac{1}{2}}
    \tag{using \Cref{neqn:conc-fitted-folded-normal-01}}                                     \\
     & = \frac{1}{n}\sum_{i = 1}^{n}(\widehat{\absresp}_{i} \vee \eta)\widehat{\absresp}_{i}
    + 2 \eta \parens{\delta r_{n,2}(\mu_n, \mu_1, \sigma)}^{\frac{1}{2}}.
    \label{neqn:conc-fitted-folded-normal-02-02}
\end{align}

Note that the final equality in \Cref{neqn:conc-fitted-folded-normal-02-02}
holds, since $\sum_{i = 1}^{n}(\widehat{\absresp}_{i} \vee \eta)\absresp_{i} =
\sum_{i = 1}^{n}(\widehat{\absresp}_{i} \vee \eta)\widehat{\absresp}_{i}$, by
again since we know that $\widehat{\absresp}_{i}$ are the \texttt{PAVA}
solutions. We then apply Cauchy-Schwartz to this summand of
\Cref{neqn:conc-fitted-folded-normal-02-02} to obtain the following upper bound
with probability at least $1 - \delta^{-1}$, for $\delta^{-1} \in (0, 1)$.
\begin{equation}\label{neqn:conc-fitted-folded-normal-03}
    \frac{1}{n} \sum_{i = 1}^{n}(\widehat{\absresp}_{i} \vee \eta)^2
    \leq \frac{1}{n}\sqrt{\sum_{i = 1}^{n}(\widehat{\absresp}_{i} \vee \eta)^2}
    \sqrt{\sum_{i = 1}^{n}\widehat{\absresp}_{i}^2}
    + 2 \eta \parens{\delta r_{n,2}(\mu_n, \mu_1, \sigma)}^{\frac{1}{2}},
\end{equation}

\nit Now observe that since $\eta > 0$, the following holds \as:
\begin{equation}\label{neqn:eta-upper-bound-t-hat}
    \eta
    = \absa{\eta}
    = \sqrt{\frac{1}{n} \sum_{i = 1}^{n}\eta^{2}}
    \leq \sqrt{\frac{1}{n} \sum_{i = 1}^{n}(\widehat{\absresp}_{i} \vee \eta)^2}
\end{equation}

\nit Then using \Cref{neqn:eta-upper-bound-t-hat} we have the following:
\begin{align}
     & \eta \parens{\sqrt{\frac{1}{n} \sum_{i = 1}^{n}(\widehat{\absresp}_{i} \vee \eta)^2} - \sqrt{\frac{1}{n} \sum_{i = 1}^{n}\widehat{\absresp}_{i}^2}}
    \nonumber                                                                                                                                                                                                                      \\
     & \leq \sqrt{\frac{1}{n} \sum_{i = 1}^{n}(\widehat{\absresp}_{i} \vee \eta)^2} \parens{\sqrt{\frac{1}{n} \sum_{i = 1}^{n}(\widehat{\absresp}_{i} \vee \eta)^2} - \sqrt{\frac{1}{n} \sum_{i = 1}^{n}\widehat{\absresp}_{i}^2}}
    \tag{using \Cref{neqn:eta-upper-bound-t-hat}}                                                                                                                                                                                  \\
     & = \frac{1}{n} \sum_{i = 1}^{n}(\widehat{\absresp}_{i} \vee \eta)^2  - \frac{1}{n}\sqrt{\sum_{i = 1}^{n}(\widehat{\absresp}_{i} \vee \eta)^2}
    \sqrt{\sum_{i = 1}^{n}\widehat{\absresp}_{i}^2}.
    \label{neqn:eta-upper-bound-t-hat-01}
\end{align}

By applying the upper bound derived in \Cref{neqn:conc-fitted-folded-normal-03}
to \Cref{neqn:eta-upper-bound-t-hat-01} we obtain the following:
\begin{equation}
    \sqrt{\frac{1}{n} \sum_{i = 1}^{n}(\widehat{\absresp}_{i} \vee \eta)^2} - \sqrt{\frac{1}{n} \sum_{i = 1}^{n}\widehat{\absresp}_{i}^2}
    \leq  2 \parens{\delta r_{n,2}(\mu_n, \mu_1, \sigma)}^{\frac{1}{2}},
    \label{neqn:eta-upper-bound-t-hat-02}
\end{equation}
with probability at least $1 - \delta^{-1}$, for $\delta^{-1} \in (0, 1)$.

Now we will show that $\sqrt{\frac{1}{n} \sum_{i =
            1}^{n}\widehat{\absresp}_{i}^2}$ is a constant distance away from
            $\sqrt{\frac{1}{n} \sum_{i = 1}^{n}\absresp_{i}^{2}}$, which will
            imply that for large $n$ the value at $0$ is smaller than the target
            value, \ie, $G(\sigma_{1}) \defined G(0) \leq \frac{1}{n} \sum_{i =
            1}^n \absresp_{i}^{2}$ as required per
            \Cref{neqn:g-sigma-condn-lower}.

On the other hand using \Cref{nlem:conc-folded-normal}, we have:
\begin{equation}\label{neqn:conc-folded-normal-02}
    \absa{\frac{1}{n}\sum_{i= 1}^{n} ( \absresp_{i} - f(\mu_{i}, \sigma))^{2} -
    \frac{1}{n}\sum_{i = 1}^{n} \parens{\mu_{i}^{2} + \sigma^{2} - f(\mu_{i}, \sigma)^2}}
    \leq l(\gamma, \univupperboundmu, \sigma),
\end{equation}
with probability at least $1 - \gamma^{-2}$, where $l(\gamma, \univupperboundmu,
    \sigma) \defined \gamma \sigma \sqrt{\frac{5 \sigma^{2} + 4
    \univupperboundmu}{n}}$. Subtracting the inequalities in
    \Cref{neqn:conc-fitted-folded-normal-01,neqn:conc-folded-normal-02} we then
    obtain:
\begin{align}
         & \frac{1}{n}\sum_{i = 1}^{n} \absresp_{i}^{2}
    - \frac{1}{n}\sum_{i = 1}^{n} \widehat{\absresp}_{i}^{2}
    + \frac{2}{n}\sum_{i = 1}^{n} (\widehat{\absresp}_{i} - \absresp_{i}) f(\mu_{i},\sigma)
    \nonumber                                                                                           \\
         & \geq \frac{1}{n} \sum_{i = 1}^{n} \parens{\mu_{i}^{2} + \sigma^{2} - f(\mu_{i}, \sigma)^{2}}
    - l(\gamma, \univupperboundmu, \sigma) - \delta r_{n,2}(\mu_n, \mu_1, \sigma)
    \label{neqn:samp-mean-diff-t-that-01}                                                               \\
    \iff & \frac{1}{n}\sum_{i = 1}^{n} \absresp_{i}^{2}
    - \frac{1}{n}\sum_{i = 1}^{n} \widehat{\absresp}_{i}^{2}
    \nonumber                                                                                           \\
         & \geq \frac{1}{n}\sum_{i = 1}^{n} \parens{\mu_{i}^{2} + \sigma^{2} - f(\mu_{i}, \sigma)^{2}}
    \nonumber                                                                                           \\
         & - \parens{l(\gamma, \univupperboundmu, \sigma)
        + \delta r_{n,2}(\mu_n, \mu_1, \sigma)
        + \frac{2}{n}\sum_{i = 1}^{n} (\widehat{\absresp}_{i} - \absresp_{i}) f(\mu_{i},\sigma)}.
    \label{neqn:samp-mean-diff-t-that-02}
\end{align}
Now, in order sharpen the lower bound in \Cref{neqn:samp-mean-diff-t-that-02},
we upper bound the term $\frac{2}{n}\sum_{i = 1}^{n} (\widehat{\absresp}_{i} -
\absresp_{i}) f(\mu_{i},\sigma)$ as follows:
\begin{align}
     & \frac{2}{n}\sum_{i = 1}^{n} (\widehat{\absresp}_{i} - \absresp_{i}) f(\mu_{i},\sigma)
    \nonumber                                                                                                                                                   \\
     & = \frac{2}{n}\sum_{i = 1}^{n} (\widehat{\absresp}_{i} - \absresp_{i})(f(\mu_{i},\sigma) - \widehat{\absresp}_{i})
    \tag{since $\widehat{\absresp}_{i}$ are the \texttt{PAVA} solutions.}                                                                                       \\
     & \leq \frac{2}{n} \sqrt{\sum_{i = 1}^{n} (\widehat{\absresp}_{i} - \absresp_{i})^2} \sqrt{\sum_{i = 1}^{n}(f(\mu_{i},\sigma) - \widehat{\absresp}_{i})^2}
    \tag{by Cauchy-Schwartz.}                                                                                                                                   \\
     & \leq  \frac{2}{n}\sqrt{\sum_{i = 1}^{n} (f(\mu_{i}, \sigma) - \absresp_{i})^2} \sqrt{\sum_{i = 1}^{n}(f(\mu_{i},\sigma) - \widehat{\absresp}_{i})^2}
    \tag{since $\widehat{\absresp}_{i}$ are \texttt{PAVA}, \ie, LSE solutions.}                                                                                 \\
     & = 2 \sqrt{\frac{1}{n} \sum_{i = 1}^{n} \parens{f(\mu_{i}, \sigma) - \absresp_{i}}^2}
    \sqrt{\frac{1}{n} \sum_{i = 1}^{n}\parens{f(\mu_{i},\sigma) - \widehat{\absresp}_{i}}^2}
    \label{neqn:samp-mean-diff-t-that-03}                                                                                                                       \\
     & \leq 2\parens{l(\gamma, \univupperboundmu, \sigma)
        \delta r_{n,2}(\mu_n, \mu_1, \sigma)}^{\frac{1}{2}},
    \label{neqn:samp-mean-diff-t-that-04}
\end{align}
with probability at least with probability at least $1 - \delta^{-1} -
    \gamma^{-2}$, by the union bound. Note that to obtain
    \Cref{neqn:samp-mean-diff-t-that-04} we applied the bounds in
    \Cref{neqn:conc-fitted-folded-normal-04,neqn:conc-folded-normal-02} to
    \Cref{neqn:samp-mean-diff-t-that-03}. Now using the bound in
    \Cref{neqn:samp-mean-diff-t-that-04} in \Cref{neqn:samp-mean-diff-t-that-02}
    we conclude that:
\begin{align}
     & \frac{1}{n}\sum_{i = 1}^{n} \absresp_{i}^{2}
    - \frac{1}{n}\sum_{i = 1}^{n} \widehat{\absresp}_{i}^{2}
    \nonumber                                              \\
     & \geq \frac{1}{n}\sum_{i = 1}^{n} g(\mu_{i}, \sigma)
    \nonumber                                              \\
     & - \parens{l(\gamma, \univupperboundmu, \sigma)
        + \delta r_{n,2}(\mu_n, \mu_1, \sigma)
        + 2\parens{l(\gamma, \univupperboundmu, \sigma)
            \delta r_{n,2}(\mu_n, \mu_1, \sigma)}^{\frac{1}{2}}}
    \label{neqn:samp-mean-diff-t-that-05}
\end{align}
Now from \Cref{neqn:fold-norm-var-inc-in-mu} we have that $g(\mu, \sigma) \geq
    g(0, \sigma) = \sigma^{2}\parens{1 - \frac{2}{\pi}}$, for each $\mu > 0$.
    Hence if $\sigma \geq \univlowerboundsigma > 0$, then
\begin{equation}\label{neqn:samp-mean-diff-t-that-06}
    \frac{1}{n}\sum_{i = 1}^{n} \absresp_{i}^{2} -
    \frac{1}{n}\sum_{i = 1}^{n} \widehat{\absresp}_{i}^{2}
    \geq \univlowerboundsigma^{2} \parens{1 - \frac{2}{\pi}}
    - \parens{l(\gamma, \univupperboundmu, \sigma)
        + \delta r_{n,2}(\mu_n, \mu_1, \sigma)
        + 2\parens{l(\gamma, \univupperboundmu, \sigma)
            \delta r_{n,2}(\mu_n, \mu_1, \sigma)}^{\frac{1}{2}}}
    > 0
\end{equation}
$\frac{1}{n}\sum_{i = 1}^{n} \absresp_{i}^{2} - \frac{1}{n}\sum_{i = 1}^{n}
    \widehat{\absresp}_{i}^{2} \geq \univlowerboundsigma^{2} \parens{1 -
    \frac{2}{\pi}} > 0$, with probability at least $1 - \delta^{-1} -
    \gamma^{-2}$. Hence under the assumption that $r_{n,2}(f, \mu_n, \mu_1,
    \sigma) = o(1)$, for sufficiently large $n$ the above will be bigger than a
    constant. Now by Lemma \ref{nlem:conc-samp-mean-ti-sq} we have
    $\frac{1}{n}\sum_{i = 1}^n T_i^2 \leq \frac{1}{n}\sum_{i = 1}^n (\mu_i^2 +
    \sigma^2) + 2\gamma \sigma \sqrt{\frac{2\sigma^2 + 4 \univupperboundmu}{n}}$
    which is upper bounded by some constant for sufficiently large $n$ given our
    assumption that $\frac{1}{n} \sum_{i = 1}^n \mu_i^2 \leq C$, for each $n \in
    \nats$, and $\sigma \leq \univupperboundsigma$ for some constants
    $\univupperboundmu, \univupperboundsigma > 0$. It follows that by applying
    \Cref{nlem:lower-bound-via-square-diff} to
    \Cref{neqn:samp-mean-diff-t-that-06} we have:
\begin{equation}
    \sqrt{\frac{1}{n}\sum_{i = 1}^{n}
        \absresp_{i}^{2}} - \sqrt{\frac{1}{n}\sum_{i = 1}^{n} \widehat{\absresp}_{i}^{2}} \geq \kappa > 0,
\end{equation}
for sufficiently large $n$, where $\kappa$ is some positive constant, with
probability at least $1-  \delta^{-1} - 2\gamma^{-2}$.

Going back to equation \eqref{neqn:eta-upper-bound-t-hat-02}, it follows that
required equation will have a solution between $\brackets{0, \sqrt{\frac{1}{n}
\sum \absresp_{i}^{2}}}$, with probability at least with probability at least $1
- \delta^{-1} - 2\gamma^{-2}$, as required.
\eprfof


\subsection{Proof of
    \Cref{nthm:sigma-hat-close-sigma}}\label{subsec:sigma-hat-close-sigma-proofs}

\bnlem[Upper and lower bounds for
    $\widehat{\sigma}$]\label{nlem:bounds-sigma-hat} Assume that there exist
    constants $\univlowerboundsigma, \univupperboundsigma, \univupperboundmu >
    0$ such that $\univlowerboundsigma \leq \sigma \leq \univupperboundsigma$
    and $\frac{1}{n}\sum_{i = 1}^n \mu_i^2 \leq \univupperboundmu$, for each $n
    \in \nats$. Then under the \texttt{ASCI} setting per
    \Cref{ndefn:isotonic-adversarial-sign-model}, the following hold:
\begin{align}
    \text{$\widehat{\sigma} \geq K_1$ with probability at least $1 - \delta^{-1}- 2\gamma^{-2}$, for sufficiently large $n$.}
    \label{neqn:lower-bounds-sigma-hat} \\
    \text{$\widehat{\sigma} \leq K_{2}$ with probability at least $1 - \delta^{-1}- 2\gamma^{-2}$, for sufficiently large $n$,}
    \label{neqn:upper-bounds-sigma-hat}
\end{align}
where $K_{1}, K_{2} > 0$ are fixed constants and $\gamma^{-2}, \delta^{-1} \in
    (0,1)$  are as in the proof of Theorem \ref{nthm:eqn-has-unique-root}.
\enlem

\bprfof{\Cref{nlem:bounds-sigma-hat}}\label{prf:nlem:bounds-sigma-hat} We prove
each property (\Cref{neqn:lower-bounds-sigma-hat,neqn:upper-bounds-sigma-hat})
in turn. \\
\newline
(\textit{Proof of \Cref{neqn:lower-bounds-sigma-hat}}.) We note that by
assumption we have $0 < \univlowerboundsigma \leq \sigma \leq
\univupperboundsigma$. We now want to show that $\widehat{\sigma}$ is positively
bounded away from 0, with high probability. First, observe that per
\Cref{nthm:eqn-has-unique-root} that $\widehat{\sigma}$ uniquely solves
$G(\widehat{\sigma}) = \frac{1}{n} \sum_{i= 1}^{n} \absresp_{i}^{2}$, with high
probability. Per \Cref{neqn:g-sigma-value-zero}, we then have that:
\begin{equation}\label{neqn:lower-bounds-sigma-hat-01}
    G(\widehat{\sigma}) - G(0)
    = \frac{1}{n} \sum_{i= 1}^{n} \widehat{\absresp}_{i}^{2}
    - \frac{1}{n} \sum_{i= 1}^{n} (\widehat{\absresp}_{i} \vee \eta)^2.
\end{equation}
We then have by the Mean Value Theorem, and the fact that $G(\sigma)$ is
increasing for each $\sigma \geq 0$ (per \Cref{nlem:prop-g-sigma}), that there
exists a $\widetilde{\sigma} \in [0, \widehat{\sigma}]$ such that
\begin{equation}\label{neqn:lower-bounds-sigma-hat-02}
    G(\widehat{\sigma}) - G(0)
    = G^{\prime}(\widetilde{\sigma}) \widehat{\sigma}.
\end{equation}
Now we have $\widetilde{\sigma} \leq \widehat{\sigma}$, or equivalently that $2
    \widetilde{\sigma} \leq 2 \widehat{\sigma}$. Since $G^{\prime}(\sigma) \leq
    2 \sigma$ using \Cref{neqn:prop-g-sigma-03}, it follows that
    $G^{\prime}(\sigma) \leq 2 \widetilde{\sigma} \leq 2 \widehat{\sigma}$.
    Using this and \Cref{neqn:lower-bounds-sigma-hat-02}, we see that:
\begin{equation}\label{neqn:lower-bounds-sigma-hat-03}
    G(\widehat{\sigma}) - G(0)
    = G^{\prime}(\widetilde{\sigma}) \widehat{\sigma}
    \leq \parens{2 \widehat{\sigma}} \widehat{\sigma}
    \leq 2 \widehat{\sigma}^{2},
\end{equation}
Now using \Cref{neqn:lower-bounds-sigma-hat-03} and the proof of
\Cref{nthm:eqn-has-unique-root} we have that $G(\widehat{\sigma}) - G(0) =
\frac{1}{n} \sum_{i = 1}^{n} \absresp^{2}_{i} - \frac{1}{n} \sum_{i= 1}^{n}
(\widehat{\absresp}_{i} \vee \eta)^2$ is positively bounded away from 0 with
high probability. So it follows that $\widehat{\sigma} \geq
\sqrt{\frac{G(\widehat{\sigma}) - G(0)}{2}} > 0$, with high probability, as
required. \qedbsquare \\
\newline
(\textit{Proof of \Cref{neqn:upper-bounds-sigma-hat}}.) First, observe that per
\Cref{nthm:eqn-has-unique-root} that $\widehat{\sigma}$ uniquely solves
$G(\widehat{\sigma}) = \frac{1}{n} \sum_{i= 1}^{n} \absresp_{i}^{2}$, with high
probability. By \Cref{ndefn:g-sigma} this implies that $\widehat{\sigma} \leq
\frac{1}{n} \sum_{i= 1}^{n} \absresp_{i}^{2}$, with high probability. Moreover
by \Cref{nlem:conc-samp-mean-ti-sq} we have $\frac{1}{n}\sum_{i = 1}^n
\absresp_{i}^{2} \leq \frac{1}{n} \sum_{i = 1}^n (\mu_{i}^{2} + \sigma^{2}) +
2\gamma \sigma \sqrt{\frac{2\sigma^2 + 4 \univupperboundmu}{n}}$ with
probability at least $\gamma^{-1}$, where $1 - \gamma^{-2}$ for $\gamma \in (0,
1)$. This in turn is bounded, in high probability, by some constant, $K_{2} > 0$
for sufficiently large $n$ given our assumptions $\frac{1}{n} \sum_{i = 1}^n
\mu_i^2 \leq C$, for each $n \in \nats$, and $\sigma \leq \univupperboundsigma$,
as required. \qedbsquare \\
\newline
\nit Thus all properties specified in
\Cref{neqn:lower-bounds-sigma-hat,neqn:upper-bounds-sigma-hat} are now proved.
\eprfof

\nthmsigmahatclosesigma*
\bprfof{\Cref{nthm:sigma-hat-close-sigma}}\label{prf:sigma-hat-close-sigma}
Recall our map $G(\sigma) \defined \sigma^{2} + \frac{1}{n} \sum_{i = 1}^{n}
(f^{-1}(\widehat{\absresp}_{i} \vee f(\eta, \sigma), \sigma))^2$ as originally
defined in \Cref{neqn:g-sigma-function-ascifit}. We will first try to show that
$G(\sigma)$ is close to $\frac{1}{n}\sum_{i = 1}^{n}\absresp_{i}^{2}$. First
note that $f^{-1}(\cdot \vee f(\eta, \sigma), \sigma)$ is a $L \defined
\frac{1}{2\Phi(\eta/\sigma)- 1}$-Lipschitz function per
\Cref{nlem:prop-folded-normal-inv-mean} and the fact that $\sigma$ is a (both
upper and lower) bounded quantity by assumption. Thus it follows that

\begin{equation}\label{neqn:sigma-hat-close-sigma-01}
    \absa{f^{-1}(\widehat{\absresp}_{i} \vee f(\eta, \sigma), \sigma) - f^{-1}(f(\mu_{i}, \sigma), \sigma)}
    \leq L \absa{\widehat{\absresp}_{i} \vee f(\eta, \sigma) - f(\mu_{i}, \sigma)},
\end{equation}
and therefore
\begin{align}
     & \frac{1}{n}\sum_{i = 1}^{n} (f^{-1}(\widehat{\absresp}_{i} \vee f(\eta, \sigma), \sigma) - \mu_{i})^{2}                                   \\
     & = \frac{1}{n}\sum_{i = 1}^{n} \absa{f^{-1}(\widehat{\absresp}_{i} \vee f(\eta, \sigma), \sigma) - f^{-1}(f(\mu_{i}, \sigma), \sigma)}^{2}
    \nonumber                                                                                                                                    \\
     & \leq \frac{L^{2}}{n} \sum_{i = 1}^{n}(\widehat{\absresp}_{i} \vee f(\eta, \sigma) - f(\mu_{i}, \sigma))^{2}
    \tag{using \Cref{neqn:sigma-hat-close-sigma-01}}                                                                                             \\
     & \leq \frac{L^{2}}{n} \sum_{i = 1}^{n}(\widehat{\absresp}_{i} - f(\mu_{i}, \sigma))^{2}
    \tag{using \Cref{neqn:fold-norm-strictly-inc,nlem:square-max-diff-inequality}.}
\end{align}

\nit In sum, we have established:
\begin{equation}\label{neqn:sigma-hat-close-sigma-02}
    \frac{1}{n}\sum_{i = 1}^{n} (f^{-1}(\widehat{\absresp}_{i} \vee f(\eta, \sigma), \sigma) - \mu_{i})^{2}
    \leq \frac{L^{2}}{n} \sum_{i = 1}^{n}(\widehat{\absresp}_{i} - f(\mu_{i}, \sigma))^{2},
\end{equation}

\nit We saw earlier by \Cref{nthm:conc-fitted-folded-normal} we have that

\begin{equation}\label{neqn:conc-fitted-folded-normal-04}
    \frac{1}{n} \sum_{i = 1}^{n} \parens{\widehat{\absresp}_{i}
        - f(\mu_{i}, \sigma)}^2
    \leq \delta r_{n,2}(\mu_n, \mu_1, \sigma),
\end{equation}
with probability at least $1 - \delta^{-1}$, for $\delta^{-1} \in (0, 1)$.
Combining \Cref{neqn:sigma-hat-close-sigma-02,neqn:conc-fitted-folded-normal-04}
we have that
\begin{equation}\label{neqn:sigma-hat-close-sigma-03}
    \frac{1}{n}\sum_{i = 1}^{n} (f^{-1}(\widehat{\absresp}_{i} \vee f(\eta, \sigma), \sigma) - \mu_{i})^{2}
    \leq L^{2} \delta r_{n,2}(\mu_n, \mu_1, \sigma)
\end{equation}
with probability at least $1 - \delta^{-1}$, for $\delta^{-1} \in (0, 1)$. Thus
by the triangle inequality, and reverse triangle inequality we have
\begin{align}\label{neqn:sigma-hat-close-sigma-04}
    \frac{1}{n} \sum_{i = 1}^{n} (f^{-1}(\widehat{\absresp}_{i} \vee f(\eta, \sigma), \sigma))^{2}
    \in \brackets{\frac{1}{n} \sum_{i = 1}^{n}\mu_{i}^{2} - h_{n},
        \frac{1}{n} \sum_{i = 1}^{n}\mu_{i}^{2} + L^{2} \delta r_{n,2}(\mu_n, \mu_1, \sigma)
        + h_{n}}
\end{align}
where $h_{n} \defined 2 \sqrt{\frac{1}{n} \sum_{i = 1}^{n}\mu_{i}^2}
    \sqrt{\frac{1}{n} \sum_{i = 1}^{n}(f^{-1}(\widehat{\absresp}_{i} \vee
    f(\eta, \sigma), \sigma)- \mu_{i})^2}$. Given our assumption that
    $\frac{1}{n}\sum_{i = 1}^{n}\mu_{i}^{2} \leq \univupperboundmu$, for each $n
    \in \nats$, we have that:
\begin{align}\label{neqn:sigma-hat-close-sigma-05}
    h_{n}
     & \leq 2 L \parens{\univupperboundmu \delta r_{n,2}(\mu_n, \mu_1, \sigma)}^{\frac{1}{2}}
\end{align}
with probability at least $1 - \delta^{-1}$, for $\delta^{-1} \in (0, 1)$. Then
combining \Cref{neqn:sigma-hat-close-sigma-04,neqn:sigma-hat-close-sigma-05}, we
have that there exists some $l_{1} \in [-2,2]$ for sufficiently large $n$ such
that
\begin{align}\label{neqn:sigma-hat-close-sigma-06}
    \frac{1}{n} \sum_{i = 1}^{n} (f^{-1}(\widehat{\absresp}_{i} \vee f(\eta, \sigma), \sigma))^{2}
    = \frac{1}{n} \sum_{i = 1}^{n}\mu_{i}^{2}
    + l_{1} L \parens{2 \univupperboundmu \delta r_{n,2}(\mu_n, \mu_1, \sigma)}^{\frac{1}{2}}.
\end{align}
with probability at least $1 - 2 \delta^{-1}$, for $\delta^{-1} \in (0, 1)$,
using the union bound.

Similarly, using \Cref{nlem:conc-samp-mean-ti-sq} we have that there exists some
$l_{2}(\sigma, \univupperboundmu, \gamma) \in \reals$ such that
\begin{align}\label{neqn:sigma-hat-close-sigma-07}
    \frac{1}{n} \sum_{i = 1}^{n}\absresp_{i}^{2}
    = \sigma^{2} + \frac{1}{n} \sum_{i = 1}^{n}\mu_{i}^{2} + l_{2}(\sigma,
    \univupperboundmu, \gamma) n^{-1/2}.
\end{align}
with probability at least $1 - \gamma^{-2}$, for $\gamma^{-2} \in (0, 1)$.
Moreover per \Cref{nthm:eqn-has-unique-root} we have that $\widehat{\sigma}$
uniquely solves $G(\widehat{\sigma}) = \frac{1}{n} \sum_{i= 1}^{n}
\absresp_{i}^{2}$, with high probability. That is:
\begin{align}\label{neqn:sigma-hat-close-sigma-08}
    G(\widehat{\sigma})
    = \widehat{\sigma}^{2} + \frac{1}{n} \sum_{i = 1}^{n} (f^{-1}(\widehat{\absresp}_{i} \vee f(\eta, \widehat{\sigma}), \widehat{\sigma}))^{2}
    = \frac{1}{n} \sum_{i = 1}^{n}\absresp_{i}^{2}.
\end{align}
Then combining
\Cref{neqn:sigma-hat-close-sigma-07,neqn:sigma-hat-close-sigma-08}, we conclude
that
\begin{align}\label{Gsigma-Ghatsigma}
    |G(\sigma) - G(\hat \sigma)| \leq l_{1} L \parens{2 \univupperboundmu \delta r_{n,2}(\mu_n, \mu_1, \sigma)}^{\frac{1}{2}} + l_2(\sigma, C ,\gamma)n^{-1/2}.
\end{align}

We now consider two cases, namely $\widehat{\sigma} > \sigma$ and $\sigma >
    \widehat{\sigma}$. In the first case, with $\widehat{\sigma} > \sigma$, we
    seek to show that $G^{\prime}(\xi) \geq K_{1} > 0$, in high probability, for
    each $\xi \in (\sigma, \widehat{\sigma})$.  Here $K_{1}$ represents a
    positive constant. By Lemma \ref{nlem:bounds-sigma-hat} both $\sigma$ and
    $\hat \sigma$ are upper and lower bounded by some constants which implies
    that $\xi$ is also upper and lower bounded by some constants call them $C_1$
    and $C_2$, i.e., $C_1\leq \xi \leq C_2$. Since $G'(\xi) \geq J(\xi) = \xi
    \bigg(\frac{1}{2}- \frac{\eta/\xi \phi(\eta/\xi)}{2\Phi(\eta/\xi) -
    1}\bigg)$. As we argued earlier $J(\xi)$ is positive and since it is a
    continuous function and the set $[C_1,C_2]$ is compact it achieves its
    minimum, which is strictly positive. Hence $G'(\xi) \geq K_1 > 0$.

Similarly, in the second case, with $\sigma > \widehat{\sigma}$, we can also
show that $G^{\prime}(\xi) \geq K_{2} > 0$, in high probability, for each $\xi
\in (\sigma, \widehat{\sigma})$. Where again, $K_{2}$ represents a positive
constant.

Then by using the Mean Value Theorem we have that there exists some $\xi \in
    (\sigma, \widehat{\sigma})$ such that $\absa{G(\sigma) -
    G(\widehat{\sigma})} = G^{\prime}(\xi) \absa{\sigma - \widehat{\sigma}} >
    \min(K_1,K_2)\absa{\sigma - \widehat{\sigma}}$. Thus from equation
    \eqref{Gsigma-Ghatsigma} we have
\begin{equation}\label{neqn:sigma-hat-positive-04}
    l_{1} L \parens{2 \univupperboundmu \delta r_{n,2}(\mu_n, \mu_1, \sigma)}^{\frac{1}{2}} + l_2(\sigma, C ,\gamma)n^{-1/2} = \absa{G(\sigma) - G(\widehat{\sigma})}
    \geq \min(K_1,K_2) \absa{\sigma - \widehat{\sigma}},
\end{equation}
and hence $\absa{\sigma - \widehat{\sigma}} \lesssim l_{1} L \parens{2
        \univupperboundmu \delta r_{n,2}(\mu_n, \mu_1, \sigma)}^{\frac{1}{2}} +
        l_2(\sigma, C ,\gamma)n^{-1/2}$.
\eprfof


\subsection{Proof of
    \Cref{nthm:mu-hat-close-mu}}\label{subsec:mu-hat-close-mu-proofs}

\nthmmuhatclosemu*
\bprfof{\Cref{nthm:mu-hat-close-mu}}\label{prf:mu-hat-close-mu} We will now
consider $\frac{1}{n} \sum_{i = 1}^{n}(f^{-1} (\widehat{\absresp}_{i} \vee
f(\eta, \widehat{\sigma}), \widehat{\sigma}) - \mu_{i})^2$. We observe that \as:
\begin{align}
     & \frac{1}{n} \sum_{i = 1}^{n}(f^{-1} (\widehat{\absresp}_{i} \vee f(\eta, \widehat{\sigma}), \widehat{\sigma}) - \mu_{i})^{2}
    \label{neqn:mu-hat-close-mu-01}                                                                                                                                                             \\
     & \leq \frac{1}{n} \sum_{i = 1}^{n}(f^{-1} (\widehat{\absresp}_{i} \vee f(\eta, \widehat{\sigma}), \widehat{\sigma}) - f^{-1} (\widehat{\absresp}_{i} \vee f(\eta,  \sigma),  \sigma))^{2}
    \nonumber                                                                                                                                                                                   \\
     & + \frac{1}{n} \sum_{i = 1}^{n}(f^{-1} (\widehat{\absresp}_{i} \vee f(\eta,  \sigma),  \sigma) - \mu_{i})^{2}
    \nonumber                                                                                                                                                                                   \\
     & + 2 \sqrt{\frac{1}{n} \sum_{i = 1}^{n}(f^{-1} (\widehat{\absresp}_{i} \vee f(\eta,  \sigma),  \sigma) - \mu_{i})^2}
    \sqrt{\frac{1}{n} \sum_{i = 1}^{n}(f^{-1} (\widehat{\absresp}_{i} \vee f(\eta, \widehat{\sigma}), \widehat{\sigma}) - f^{-1} (\widehat{\absresp}_{i} \vee f(\eta,  \sigma),  \sigma))^2},
    \label{neqn:mu-hat-close-mu-02}
\end{align}
where the transition between
\Cref{neqn:mu-hat-close-mu-01,neqn:mu-hat-close-mu-02} was by applying adding
and subtracting $f^{-1} (\widehat{\absresp}_{i} \vee f(\eta, \sigma),  \sigma)$,
then applying the triangle inequality, and finally applying the Cauchy-Schwartz
inequality to the cross product summand.

We now set to upper bound the \Cref{neqn:mu-hat-close-mu-02} further. First, we
saw in \Cref{neqn:sigma-hat-close-sigma-03} that $\frac{1}{n}\sum_{i = 1}^{n}
(f^{-1}(\widehat{\absresp}_{i} \vee f(\eta, \sigma), \sigma) - \mu_{i})^{2} \leq
L^{2} \delta r_{n,2}(\mu_n, \mu_1, \sigma)$, with probability at least $1 -
\delta^{-1}$, for $\delta^{-1} \in (0, 1)$. Next, we will tackle the term
\begin{equation}\label{neqn:mu-hat-close-mu-03}
    \frac{1}{n} \sum_{i = 1}^{n}(f^{-1} (\widehat{\absresp}_{i} \vee f(\eta, \widehat{\sigma}), \widehat{\sigma})
    - f^{-1} (\widehat{\absresp}_{i} \vee f(\eta,  \sigma),  \sigma) )^2.
\end{equation}
Note that map $\sigma \mapsto f^{-1}(\widehat{\absresp}_{i} \vee f(\eta,
    \sigma), \sigma)$ is a $L := \frac{\sqrt{2/\pi} \exp(-\mu^2/\sigma^2/2)}{2
    \Phi(\mu/\sigma) - 1} \leq \frac{\sqrt{2/\pi}}{2 \Phi(\eta/\sigma) -
    1}$-Lipschitz per
    \Cref{nlem:lipschitz-bounded-derivative,nlem:prop-folded-normal-inv-mean},
    and in addition both $\sigma, \hat \sigma$ are upper and lower bounded by
    constants. It follows that
\begin{align}
     & \frac{1}{n} \sum_{i = 1}^{n}(f^{-1} (\widehat{\absresp}_{i} \vee f(\eta, \widehat{\sigma}), \widehat{\sigma})
    - f^{-1} (\widehat{\absresp}_{i} \vee f(\eta,  \sigma),  \sigma) )^{2}
    \nonumber                                                                                                        \\
     & \leq \frac{1}{n} \sum_{i = 1}^{n} \absa{L (\sigma - \widehat{\sigma})}^{2}
    \tag{using $L$-Lipschitz property.}                                                                              \\
     & = L^{2} (\sigma - \widehat{\sigma})^{2}
    \nonumber                                                                                                        \\
     & \lesssim 2 \univupperboundmu \delta r_{n,2}(\mu_n, \mu_1, \sigma) + \gamma^2 n^{-1},
    \label{neqn:mu-hat-close-mu-04}
\end{align}
where \Cref{neqn:mu-hat-close-mu-04} follows from
\Cref{nthm:sigma-hat-close-sigma} with probability at least $1 - \delta^{-1} -
2\gamma^{-2}$.

Then applying the upper bounds in
\Cref{neqn:mu-hat-close-mu-03,neqn:mu-hat-close-mu-04} appropriately to each
corresponding summand of \Cref{neqn:mu-hat-close-mu-02}, we conclude that
\begin{equation}\label{neqn:mu-hat-close-mu-05}
    \frac{1}{n} \sum_{i = 1}^{n}(f^{-1} (\widehat{\absresp}_{i} \vee f(\eta, \widehat{\sigma}), \widehat{\sigma}) - \mu_{i})^{2}
    \lesssim \delta r_{n,2}(\mu_n, \mu_1, \sigma) + \gamma^2 n^{-1},
\end{equation}
with probability at least $1 - \delta^{-1} - 2\gamma^{-2}$.
\eprfof

\newpage
\section{Proofs of \Cref{sec:lower-bounds}}\label{sec:lower-bounds-proofs}

\subsection{Mathematical
    Preliminaries}\label{subsec:preliminarieslower-bounds-proofs}

Since we adapt the lower bound construction from
\citet{bellec2015sharporaclemonoconvexreg} for our \texttt{ASCI} setting, we
first introduce the relevant related notation and definitions here first for
classes of monotonic sequences. We denote $\mclS^{\uparrow} \defined
\thesetb{\boldsymbol{\mu} \defined (\mu_1, \ldots, \mu_{n})^{\top}}{\mu_1 \leq
\ldots \leq \mu_{n}}$ to be the set of all non-decreasing sequences. We define
$k(\boldsymbol{\mu}) \geq 1$, for $\boldsymbol{\mu} \in \mclS^{\uparrow}$, to be
the integer such that $k(\boldsymbol{u})-1$ is the number of inequalities
$\mu_{i} \leq \mu_{i+1}$ that are strict for $i \in [n-1]$ (\ie, number of jumps
of $\boldsymbol{\mu}$). The class of monotone functions we will consider are
$\mclS^{\uparrow}(V^{*}) \defined \thesetb{\boldsymbol{\mu} \in
\mclS^{\uparrow}}{V(\boldsymbol{\mu}) \leq V^{*}}$, for some fixed $V^{*} \in
\reals$, and $V(\boldsymbol{\mu}) = \mu_{n} - \mu_{1}$, is the total variation
of any $\boldsymbol{\mu} \in \mclS^{\uparrow}$. We also consider the restricted
class of monotone sequences, $\mclS^{\uparrow}_{k^{*}} \defined
\thesetb{\boldsymbol{\mu} \in \mclS^{\uparrow}}{k(\boldsymbol{\mu}) \leq
k^{*}}$, and $\mclS^{\uparrow}(V^{*}, \eta, \univupperboundmu) \defined
\thesetb{\boldsymbol{\mu} \in \mclS^{\uparrow}(V^{*})}{\frac{1}{n}\sum_{i =
1}^{n} \mu_{i}^{2} \leq \univupperboundmu, \mu_{1} > \eta > 0}$.

\subsection{Proof of
    \Cref{nprop:minimax-lower-bounds}}\label{subsec:nthm-minimax-lower-bounds-proofs}

We follow directly the proof technique and construction from
\citet[Proposition~4]{bellec2015sharporaclemonoconvexreg}, but make suitable
adaptations for our \texttt{ASCI} setup. We largely follow their notation to
help readers align the commonalities and differences in the underlying
constructions used. Our first lower bound result is stated in
\Cref{nprop:minimax-lower-bounds}.

\begin{restatable}[Minimax lower
        bounds]{nprop}{nthmminimaxlowerbounds}\label{nprop:minimax-lower-bounds}
        Let $n \geq 2, V^{*} > 0$ and $\sigma > 0$. There exist absolute
        constants $c, c^{\prime}>0$ such that for any positive integer $k^{*}
        \leq n$ satisfying $\parens{k^{*}}^{3} \leq \frac{16 n
        \parens{V^{*}}^{2}}{\sigma^{2}}$ we have
    \begin{equation}\label{neqn:minimax-lower-bounds-01}
        \inf_{\hat{\boldsymbol{\mu}}}
        \sup_{\mclS^{\uparrow}_{k^{*}} \, \cap \, \mclS^{\uparrow}(V^{*}, \eta, \univupperboundmu)}
        \Prba{\boldsymbol{\mu}}{\frac{1}{n} \norma{\widehat{\boldsymbol{\mu}} -
                     \boldsymbol{\mu}}^{2} \geq \frac{c \sigma^{2} k}{n}} > c^{\prime}
    \end{equation}
    where $\eta > 0$ is a fixed positive constant per
    \Cref{neqn:adversarial-iso-reg-2}, $\univupperboundmu \geq
    \parens{V^{*}}^{2} + 4 \gamma^{2} + 2 \eta^{2}$, $\gamma \defined
    \frac{1}{8} \sqrt{\frac{\sigma^{2} k^{*}}{n}}$, $\mathop{\bf P
    \/}_{\boldsymbol{\mu}}$ denotes the distribution of $\parens{R_{1}, \ldots,
    R_{n}}^{\top}$ satisfying \Cref{neqn:adversarial-iso-reg-1}, and
    $\inf_{\hat{\boldsymbol{\mu}}}$ is the infimum over all estimators.
\end{restatable}

\bprfof{\Cref{nprop:minimax-lower-bounds}} Let $n$ be a multiple of $k^{*} \in
    \nats$. Then for any $\boldsymbol{\omega}, \boldsymbol{\omega^{\prime}} \in
    \theset{0, 1}^{k^{*}}$, using the Varshamov-Gilbert bound
    \citep[Lemma~2.9]{tsybakov2009intrononparmestimation}, there exists a set
    $\Omega \in \theset{0, 1}^{k^{*}}$ such that:

\begin{equation}
    \bfzero = (0, \ldots, 0)^{\top} \in \Omega,
    \quad \log (\absa{\Omega} - 1) \geq \frac{k^{*}}{8},
    \quad \text { and }
    \quad \hammd{\boldsymbol{\omega}}{\boldsymbol{\omega}^{\prime}} > \frac{k^{*}}{8}
\end{equation}

\nit for any two distinct $\boldsymbol{\omega}, \boldsymbol{\omega}^{\prime} \in
    \Omega$. For each $\boldsymbol{\omega} \in \Omega$, define a vector
    $\bfu^{\boldsymbol{\omega}} \in \mathbb{R}^{n}$ componentwise, for each
    component index $i \in [n]$ as follows:

\begin{align}
    \bfu_{i}^{\boldsymbol{\omega}}
     & \defined \frac{\floor{(i-1) \frac{k^{*}}{n}} V^{*}}{2 k^{*}}
    + \gamma \boldsymbol{\omega}_{\floor{(i-1) \frac{k^{*}}{n}} + 1}.
    \label{neqn:minimax-lower-bounds-bellec-ui}                     \\
    \bar{\bfu}_{i}^{\boldsymbol{\omega}}
     & \defined \bfu_{i}^{\boldsymbol{\omega}} + \eta.
    \label{neqn:minimax-lower-bounds-bellec-ui-eta}
\end{align}

\nit where $\gamma \defined \frac{1}{8} \sqrt{\frac{\sigma^{2} k^{*}}{n}}$ and
$\eta > 0$ is a fixed positive constant per \Cref{neqn:adversarial-iso-reg-2}.
Importantly we note that $\bfu_{i}^{\boldsymbol{\omega}}$ per
\Cref{neqn:minimax-lower-bounds-bellec-ui} is precisely as constructed in
\citet[Proposition~4]{bellec2015sharporaclemonoconvexreg}. However, critically
the construction in \Cref{neqn:minimax-lower-bounds-bellec-ui-eta} is adapted to
our \texttt{ASCI} setting, by componentwise translation by $\eta > 0$. More
compactly, it is also convenient to represent this construction as
$\bar{\bfu}^{\boldsymbol{\omega}} \defined \bfu^{\boldsymbol{\omega}} +
\boldsymbol{\eta}$, where $\boldsymbol{\eta} \defined (\eta, \ldots,
\eta)^{\top} \in \reals^{n}$.

As per \citet[Proposition~4]{bellec2015sharporaclemonoconvexreg} we first note
the following properties for $\bfu_{i}^{\boldsymbol{\omega}}$, for each $i \in
[n]$. For any $\boldsymbol{\omega} \in \Omega, \bfu^{\boldsymbol{\omega}}$ is a
piecewise constant sequence with $k\parens{\bfu^{\boldsymbol{\omega}}} \leq
k^{*}, \bfu^{\boldsymbol{\omega}}$ is a non-decreasing sequence because $\gamma
\leq \frac{V^{*}}{2 k^{*}}$, and by construction
$V\parens{\bfu^{\boldsymbol{\omega}}} \leq V^{*}$. Thus,
$\bfu^{\boldsymbol{\omega}} \in \mathcal{S}_{k^{*}}^{\uparrow} \cap
\mathcal{S}^{\uparrow}(V)$ for all $\boldsymbol{\omega} \in \Omega$.

Now we observe the following corresponding properties of the
$\boldsymbol{\eta}$-translated sequence $\bar{\bfu}^{\boldsymbol{\omega}}$.
First note that since for any $\boldsymbol{\omega} \in \Omega,
\bfu^{\boldsymbol{\omega}}$ is a piecewise constant non-decreasing sequence, so
is $\bar{\bfu}_{j}^{\boldsymbol{\omega}}$, by translation invariance. Next,
consider any arbitrary index $j \in [n]$ relating to a `jump' in
$\bfu^{\boldsymbol{\omega}}$, \ie, $\bfu_{j}^{\boldsymbol{\omega}} < \bfu_{j +
1}^{\boldsymbol{\omega}}$ (note the \textit{strict} inequality). We then have
that:
\begin{align}
    \bfu_{j}^{\boldsymbol{\omega}}
     & < \bfu_{j + 1}^{\boldsymbol{\omega}}
    \tag{by assumption.}                           \\
    \iff \bfu_{j}^{\boldsymbol{\omega}} + \eta
     & < \bfu_{j + 1}^{\boldsymbol{\omega}} + \eta
    \nonumber                                      \\
    \iff \bar{\bfu}_{j}^{\boldsymbol{\omega}}
     & < \bar{\bfu}_{j + 1}^{\boldsymbol{\omega}}
    \tag{using \Cref{neqn:minimax-lower-bounds-bellec-ui-eta}}
\end{align}
So any `jump' in the original sequence $\bfu^{\boldsymbol{\omega}}$ corresponds
to a jump in the $\boldsymbol{\eta}$-translated sequence
$\bar{\bfu}^{\boldsymbol{\omega}}$. That is, we have
$k\parens{\bar{\bfu}^{\boldsymbol{\omega}}} =
k\parens{\bfu^{\boldsymbol{\omega}}} \leq k^{*}$. In addition, we note that
\begin{align*}
    V\parens{\bar{\bfu}^{\boldsymbol{\omega}}}
     & = \bar{\bfu}_{n}^{\boldsymbol{\omega}} - \bar{\bfu}_{1}^{\boldsymbol{\omega}}                     \\
     & = \parens{\bfu_{n}^{\boldsymbol{\omega}} + \eta} - \parens{\bfu_{1}^{\boldsymbol{\omega}} + \eta} \\
     & = \bfu_{n}^{\boldsymbol{\omega}} - \bfu_{1}^{\boldsymbol{\omega}}                                 \\
     & = V\parens{\bfu^{\boldsymbol{\omega}}}                                                            \\
     & \leq V^{*}
    \tag{by construction of $\bfu^{\boldsymbol{\omega}}$.}
\end{align*}

By construction we also have that $\bar{\bfu}_{1}^{\boldsymbol{\omega}} \defined
    \bfu_{1}^{\boldsymbol{\omega}} + \eta \geq \eta > 0$, since
    $\bfu_{1}^{\boldsymbol{\omega}} \geq 0$ by construction (in fact each
    component is non-negative). Finally, per our \texttt{ASCI} setting, we want
    to check if there exists a $\univupperboundmu > 0$, such that $\frac{1}{n}
    \sum_{i = 1}^{n} \parens{\bar{\bfu}_{i}^{\boldsymbol{\omega}}}^{2} \leq
    \univupperboundmu$, for each $n \in \nats$. Given
    $\bar{\bfu}^{\boldsymbol{\omega}}$ we observe the following for each
    component index $i \in [n]$:
\begin{align}
    \parens{\bar{\bfu}_{i}^{\boldsymbol{\omega}}}^{2}
     & \defined \parens{\bfu_{i}^{\boldsymbol{\omega}} + \eta}^{2}
    \tag{using \Cref{neqn:minimax-lower-bounds-bellec-ui}}                                 \\
     & \leq 2 \parens{\parens{\bfu_{i}^{\boldsymbol{\omega}}}^{2} + \eta^{2}}
    \tag{using \Cref{nlem:square-sum-inequality}}                                          \\
     & = 2 \parens{\parens{\frac{\floor{(i-1) \frac{k}{n}} V^{*}}{2 k}
    + \gamma \boldsymbol{\omega}_{\floor{(i-1) \frac{k}{n}} + 1}}^{2} + \eta^{2}}
    \tag{using \Cref{neqn:minimax-lower-bounds-bellec-ui}}                                 \\
     & \leq 2 \parens{2 \brackets{\parens{\frac{\floor{(i-1) \frac{k}{n}} V^{*}}{2 k}}^{2}
    + \parens{\gamma \boldsymbol{\omega}_{\floor{(i-1) \frac{k}{n}} + 1}}^{2}} + \eta^{2}}
    \tag{using \Cref{nlem:square-sum-inequality}}                                          \\
     & \leq 2 \parens{2 \brackets{\parens{\frac{V^{*} k}{2 k}}^{2}
            + \gamma^{2}} + \eta^{2}}
    \tag{since $\frac{i-1}{n} \leq 1$ for each $i \in [n]$.}                               \\
     & \leq \parens{V^{*}}^{2} + 4 \gamma^{2} + 2 \eta^{2}
    \label{neqn:minimax-lower-bound-sec-mom-cond}
\end{align}

\nit So indeed it follows from \Cref{neqn:minimax-lower-bound-sec-mom-cond}
that:
\begin{equation}
    \frac{1}{n} \sum_{i = 1}^{n} \parens{\bar{\bfu}_{i}^{\boldsymbol{\omega}}}^{2}
    \leq \parens{V^{*}}^{2} + 4 \gamma^{2} + 2 \eta^{2}
    \defines \univupperboundmu
\end{equation}

\nit So that we have $\bar{\bfu}^{\boldsymbol{\omega}} \in
    \mclS^{\uparrow}_{k^{*}} \, \cap \, \mclS^{\uparrow}(V^{*}, \eta,
    \univupperboundmu)$. Moreover, for any $\boldsymbol{\omega},
    \boldsymbol{\omega}^{\prime} \in \Omega$, we observe that:
\begin{align}
    \absa{\bar{\bfu}^{\boldsymbol{\omega}}-\bar{\bfu}^{\boldsymbol{\omega^{\prime}}}}^{2}
     & \defined \absa{(\bfu^{\boldsymbol{\omega}} + \boldsymbol{\eta}) - (\bfu^{\boldsymbol{\omega^{\prime}}} + \boldsymbol{\eta})}^{2}
    \tag{by construction $\bar{\bfu}^{\boldsymbol{\omega}}
        \defined \bfu^{\boldsymbol{\omega}} + \boldsymbol{\eta}$.}
    \nonumber                                                                                                                            \\
     & = \norma{\bfu^{\boldsymbol{\omega}} - \bfu^{\boldsymbol{\omega^{\prime}}}}^{2}
    \nonumber                                                                                                                            \\
     & = \frac{\gamma^{2}}{k^{*}} \hammd{\boldsymbol{\omega}}{\boldsymbol{\omega}^{\prime}}
    \nonumber                                                                                                                            \\
     & \geq \frac{\gamma^{2}}{8}
    \nonumber                                                                                                                            \\
     & = \frac{\sigma^{2} k^{*}}{512 n}
    \nonumber
\end{align}

\nit Set for brevity $\mathop{\bf P\/}_{\boldsymbol{\omega}} = \mathop{\bf
        P\/}_{\bar{\bfu}^{\boldsymbol{\omega}}}$. The Kullback-Leibler
        divergence $\kl{\mathop{\bf P\/}_{\boldsymbol{\omega}}}{\mathop{\bf
        P\/}_{\boldsymbol{\omega}^{\prime}}}$, between $\mathop{\bf
        P\/}_{\boldsymbol{\omega}}$ and $\mathop{\bf
        P\/}_{\boldsymbol{\omega}^{\prime}}$, is equal to $\frac{n}{2
        \sigma^{2}} \norma{\bfu^{\boldsymbol{\omega}} -
        \bfu^{\boldsymbol{\omega^{\prime}}}}^{2}$ for all $\boldsymbol{\omega},
        \boldsymbol{\omega}^{\prime} \in \Omega$. Thus,
\begin{align}
    \kl{\stackrel{}{{\bf P\/}_{\boldsymbol{\omega}}}}{\stackrel{}{{\bf P\/}_{\bfzero}}}
    = \frac{\gamma^{2} n \hammd{\mathbf{0}}{\boldsymbol{\omega}}}{2 k^{*} \sigma^{2}}
    \leq \frac{k^{*}}{128}
    \leq \frac{\log (\absa{\Omega} - 1)}{16}
\end{align}
Applying \citet[Theorem~2.7]{tsybakov2009intrononparmestimation} with $\alpha=1
    / 16$ completes the proof.
\eprfof

\subsection{Proof of
    \Cref{nprop:minimax-lower-bounds-final}}\label{subsec:ncor-minimax-lower-bounds-proofs}

From \Cref{nprop:minimax-lower-bounds}, in line with
\citet[Corollary~5]{bellec2015sharporaclemonoconvexreg}, we immediately obtain
the following result in \Cref{nprop:minimax-lower-bounds-final}. Once again, we
utilize the technique of \citet[Corollary~5]{bellec2015sharporaclemonoconvexreg}
to obtain the following corollary. The important changes to ensure that we adapt
to our \texttt{ASCI} setting are captured in \Cref{nprop:minimax-lower-bounds}
and our proof thereof.

\ncorminimaxlowerbounds*
\bprfof{\Cref{nprop:minimax-lower-bounds}} As per
\citet[Corollary~5]{bellec2015sharporaclemonoconvexreg}, to prove this corollary
it is enough to note that if $\frac{16 n \parens{V^{*}}^{2}}{\sigma^{2}} \geq
1$, by choosing $k^{*}$ in \Cref{nprop:minimax-lower-bounds} as the integer part
of $\parens{\frac{16 n \parens{V^{*}}^{2}}{\sigma^{2}}}^{\frac{1}{3}}$, we
obtain the lower bound corresponding to $\parens{\frac{\sigma^{2}
V^{*}}{n}}^{\frac{2}{3}}$ under the maximum in
\Cref{neqn:minimax-lower-bounds-02}. On the other hand, if $\frac{16 n
\parens{V^{*}}^{2}}{\sigma^{2}} < 1$ the term $\frac{\sigma^{2}}{n}$ is
dominant, so that we need to have the lower bound of the order
$\frac{\sigma^{2}}{n}$, which is trivial (it follows from a reduction to the
bound for the class composed of two constant functions).
\eprfof

\end{document}